\documentclass[11pt]{article}
\UseRawInputEncoding
\usepackage{currfile}
\usepackage{color}
\usepackage{mathrsfs}
\usepackage{caption}
\usepackage{subfigure}
\usepackage{setspace}
\usepackage{epstopdf}

\usepackage{amscd}
\usepackage{amssymb,amsfonts,psfrag,cite}
\usepackage{amsmath}
  \usepackage{paralist}
  \usepackage[colorlinks=true]{hyperref}
  \usepackage[all]{xy}
\usepackage{graphicx}
\newtheorem{theo}{Theorem}
\newtheorem{prop}[theo]{Proposition}

\newtheorem{defi}[theo]{Definition}
\newtheorem{lemm}[theo]{Lemma}
\newtheorem{rema}[theo]{Remark}

\makeatletter \@addtoreset{equation}{section}
\@addtoreset{theo}{section} \makeatother

\begin{document}
\date{}
\title{Nonholonomic Controlled Hamiltonian System: Symmetric Reduction
and Hamilton-Jacobi Equations}
\author{Hong Wang  \\
School of Mathematical Sciences and LPMC,\\
Nankai University,  Tianjin 300071, P.R.China \\
E-mail: hongwang@nankai.edu.cn\\\\
May 5, 2022} \maketitle

{\bf Abstract:} In order to describe the impact of nonholonomic constraints
for the dynamics of a regular controlled Hamiltonian (RCH) system,
in this paper, for an RCH system with nonholonomic constraint, we first derive its
distributional RCH system, by analyzing carefully
the structure of dynamical vector field of the nonholonomic RCH system.
Secondly, we derive precisely
the geometric constraint conditions of the induced distributional two-form
for the dynamical vector field of the distributional RCH system,
which are called the Type I and Type II of Hamilton-Jacobi equations.
Thirdly, we generalize the above results for the nonholonomic
reducible RCH system with symmetry,
and prove two types of Hamilton-Jacobi theorems
for the nonholonomic reduced distributional RCH system.
Moreover, we consider the nonholonomic
reducible RCH system with momentum map,
by combining with the regular point and regular orbit
reduction theory and the analysis of dynamics of RCH system,
we give the geometric formulations of the nonholonomic
regular point reduced and orbit reduced distributional RCH systems,
and prove two types of Hamilton-Jacobi theorems for these
reduced distributional RCH systems.
These researches reveal the deeply internal
relationships of the nonholonomic constraints,
the induced (resp. reduced) distributional two-forms,
the dynamical vector fields and controls of the nonholonomic RCH system
and its (reduced) distributional RCH systems.\\



\tableofcontents

\section{Introduction}

It is well known that the theory of controlled mechanical systems
became an important subject in recent years. Its research
gathers together some separate areas of research such as mechanics,
differential geometry and nonlinear control theory, etc., and the
emphasis of this research on geometry is motivated by the aim of
understanding the structure of equations of motion of the system, in
a way that helps both for analysis and design. Thus, it is natural to
study controlled mechanical systems by combining with the analysis
of dynamic systems and the geometric reduction theory of Hamiltonian
and Lagrangian systems. In particular,
a regular controlled Hamiltonian (RCH) system is a
Hamiltonian system with external force and control,
which is defined in Marsden et al.\cite{mawazh10},
from the viewpoint of completeness of Marsden-Weinstein symplectic
reduction. In general,
under the actions of external force and control,
an RCH system is not Hamiltonian,
however, it is a dynamical system closely related to a
Hamiltonian system, and it can be explored and studied by extending
the methods for external force and control in the study of Hamiltonian systems.\\

On the other hand, the Hamilton-Jacobi theory is an important research subject
in mathematics and analytical mechanics.
see Abraham and Marsden \cite{abma78}, Arnold
\cite{ar89} and Marsden and Ratiu \cite{mara99},
and the Hamilton-Jacobi equation is
also fundamental in the study of the quantum-classical relationship
in quantization, and it also plays an important role
in the study of stochastic dynamical systems, see
Woodhouse \cite{wo92}, Ge and Marsden \cite{gema88},
and L\'{a}zaro-Cam\'{i} and Ortega \cite{laor09}.
Hamilton-Jacobi theory from the variational
point of view is originally developed by Jacobi in 1866, which states
that the integral of Lagrangian of a system along the solution of
its Euler-Lagrange equation satisfies the Hamilton-Jacobi equation.
The classical description of this problem from the generating function and the geometrical point
of view is given by Abraham and Marsden in \cite{abma78} as follows:
Let $Q$ be a smooth manifold and $TQ$ the tangent bundle, $T^* Q$
the cotangent bundle with the canonical symplectic form $\omega$,
and the projection $\pi_Q: T^* Q \rightarrow Q $ induces the map $
T\pi_{Q}: TT^* Q \rightarrow TQ. $
\begin{theo}
Assume that the triple $(T^*Q,\omega,H)$ is a Hamiltonian system
with Hamiltonian vector field $X_H$, and $W: Q\rightarrow
\mathbb{R}$ is a given generating function. Then the following two assertions
are equivalent:\\
\noindent $(\mathrm{i})$ For every curve $\sigma: \mathbb{R}
\rightarrow Q $ satisfying $\dot{\sigma}(t)= T\pi_Q
(X_H(\mathbf{d}W(\sigma(t))))$, $\forall t\in \mathbb{R}$, then
$\mathbf{d}W \cdot \sigma $ is an integral curve of the Hamiltonian
vector field $X_H$.\\
\noindent $(\mathrm{ii})$ $W$ satisfies the Hamilton-Jacobi equation
$H(q^i,\frac{\partial W}{\partial q^i})=E, $ where $E$ is a
constant.
\end{theo}

From the proof of the above theorem given in
Abraham and Marsden \cite{abma78}, we know that
the assertion $(\mathrm{i})$ with equivalent to
Hamilton-Jacobi equation $(\mathrm{ii})$ by the generating function,
gives a geometric constraint condition of the canonical symplectic form
on the cotangent bundle $T^*Q$
for Hamiltonian vector field of the system.
Thus, the Hamilton-Jacobi equation reveals the deeply internal relationships of
the generating function, the canonical symplectic form
and the dynamical vector field of a Hamiltonian system.\\

However, note that the RCH system defined on the cotangent bundle
$T^*Q$, under the actions of external force and control, in general,
may not be a Hamiltonian system,
and it has yet no generating function,
we cannot give the Hamilton-Jacobi equations for the RCH system
and its regular reduced systems just like same as the above Theorem 1.1.
We have to look for a new way. It is worthy of noting that,
in Wang \cite{wa13d}, the author first proved
an important lemma, also see Wang \cite{wa17}, which is
a modification for the corresponding result of Abraham and Marsden
in \cite{abma78}, such that we can
derive precisely the geometric constraint conditions of
canonical symplectic form and regular reduced symplectic forms for the
dynamical vector fields of the RCH system and its regular reduced RCH systems.
which are called the Type I and Type II of Hamilton-Jacobi equations,
and moreover, can prove that the RCH-equivalence for the RCH system,
leave the solutions of corresponding Hamilton-Jacobi equations invariant.
These researches reveal the deeply internal
relationships of the geometrical structures of phase spaces, the dynamical
vector fields and controls of the RCH system.\\

In mechanics, it is very often that many systems have constraints.
A nonholonomic RCH system is an RCH system with nonholonomic constraint.
Usually, under the restriction given by nonholonomic constraint,
in general, the dynamical vector field of a nonholonomic RCH system
may not be Hamiltonian. Thus, we can not
describe the Hamilton-Jacobi equations for nonholonomic RCH system
from the viewpoint of generating function
as in the classical Hamiltonian case.
Since the Hamilton-Jacobi theory is developed based on the
Hamiltonian picture of dynamics, it is natural idea to extend the
Hamilton-Jacobi theory to the nonholonomic RCH systems,
and they with symmetry and momentum map. Note that
for a nonholonomic Hamiltonian system, from Bates and $\acute{S}$niatycki \cite{basn93},
one can derive a distributional Hamiltonian system, which is
very important, and it is called a semi-Hamiltonian
system in Patrick \cite{pa07}. In addition, there are
some important research for nonholonomic reduction and
Hamilton-Jacobi theory of nonholonomic systems, see
Bates and $\acute{S}$niatycki \cite{basn93},
Cantrijn et al. \cite{calemama99},
Cari$\tilde{n}$ena et al. \cite{cagrmamamuro06, cagrmamamuro10},
Cendra et al.\cite{cemara01},
Cushman et al. \cite{cudusn10, cukesnba95}, Gotay et al. \cite{gone79},
Koiller \cite{ko92}, Koon and Marsden in \cite{koma97, koma98},
Le\'{o}n et al. \cite{lero89}, Montgomery \cite{mo02},
and $\acute{S}$niatycki \cite{sn98, sn01} and so on.
It has been becoming one of the most active subjects in the study of
modern applied mathematics and analytical mechanics.\\

In order to describe the impact of nonholonomic constraints
for the dynamical vector fields of an RCH system and its regular
reduced RCH systems, by using the method in
Le\'{o}n and Wang \cite{lewa15} and
Bates and $\acute{S}$niatycki \cite{basn93}, and
analyzing carefully
the structure of dynamical vector field of the nonholonomic RCH system,
we derive precisely the distributional RCH system and the regular
reduced distributional RCH systems, which are determined by the non-degenerate
distributional two-forms induced from the canonical symplectic form
and reduced symplectic forms. Now, it is a natural problem how to
derive precisely the geometric constraint conditions of
the distributional two-form and the regular reduced distributional two-form for the
dynamical vector fields of the distributional RCH system and its regular reduced
distributional RCH systems, that is, two types of Hamilton-Jacobi equations.
These research are our goal in this paper.\\

The paper is organized as follows. In section 2 we first
give some definitions and basic facts about the RCH system,
the nonholonomic constraint, the nonholonomic RCH system
and the distributional RCH system,
which will be used in subsequent sections.
In section 3, we first prove an important lemma, which is a tool for our research.
Then derive precisely
the geometric constraint conditions of the non-degenerate distributional two-form
for the dynamical vector field of distributional RCH system,
that is, the two types of Hamilton-Jacobi equation.
The nonholonomic reducible RCH systems with
symmetries, as well as momentum maps, are considered respectively in
section 4 and section 5, and
by combining with the regular point and regular orbit
reduction theory and the analysis of dynamics of nonholonomic RCH system,
we can derive precisely
the geometric constraint conditions of a variety of nonholonomic reduced
distributional two-forms
for the nonholonomic reducible dynamical vector fields,
that is, the two types of Hamilton-Jacobi
equations for a variety of nonholonomic reduced distributional RCH systems.
These research work develop the nonholonomic
reduction and Hamilton-Jacobi theory of the nonholonomic RCH
systems with symmetries, as well as momentum maps,
and make us have much deeper understanding
and recognition for the structures of the nonholonomic constraint,
dynamics and control of RCH system.

\section{Nonholonomic RCH System and Distributional RCH System}

In this section, we first give some definitions and basic facts about
the RCH system, the nonholonomic constraint and the nonholonomic RCH system.
Moreover, by analyzing carefully the structure for the nonholonomic dynamical
vector field, we give a geometric formulation of distributional RCH system,
which is determined by a non-degenerate distributional two-form induced
from the canonical symplectic form.
All of them will be used in subsequent sections.
We shall follow the notations and
conventions introduced in Marsden et al. \cite{mawazh10},
Bates and $\acute{S}$niatycki in \cite{basn93}, Cushman et al.
\cite{cudusn10} and \cite{cukesnba95}, Montgomery \cite{mo02},
 Le\'{o}n and Wang \cite{lewa15}, Wang \cite{wa21c}.\\

From Marsden et al.\cite{mawazh10}, we know that
the symplectic reduced space of a Hamiltonian system
defined on the cotangent bundle of a configuration manifold may not be a
cotangent bundle, and hence the set of Hamiltonian systems with
symmetries on the cotangent bundle is not complete under the
Marsden-Weinstein reduction. This is a serious problem. If we define
directly a controlled Hamiltonian system with symmetry on a
cotangent bundle, then it is possible that the Marsden-Weinstein
reduced RCH system may not have definition.
In order to describe uniformly RCH systems defined on a cotangent
bundle and on the regular reduced spaces, we
first define an RCH system on a symplectic fiber bundle,
then we can obtain the RCH system on the cotangent bundle of a configuration
manifold as a special case.\\

Let $(\mathbb{E},M,\pi)$ be a fiber bundle, and for each point $x \in M$,
assume that the fiber $\mathbb{E}_x=\pi^{-1}(x)$ is a smooth submanifold of $\mathbb{E}$
and with a symplectic form $\omega_{\mathbb{E}}(x)$, that is,
$(\mathbb{E}, \omega_{\mathbb{E}})$ is a
symplectic fiber bundle. If for any Hamiltonian function $H: \mathbb{E} \rightarrow
\mathbb{R}$, we have a Hamiltonian vector field $X_H$,
which satisfies the Hamilton's equation, that is,
$\mathbf{i}_{X_H}\omega_{\mathbb{E}}=\mathbf{d}H$,
then $(\mathbb{E}, \omega_{\mathbb{E}}, H )$ is a
Hamiltonian system. Moreover, if considering the external force and
control, we can define a kind of regular controlled Hamiltonian
(RCH) system on the symplectic fiber bundle $\mathbb{E}$ as follows.

\begin{defi} (RCH System) An RCH system on $\mathbb{E}$ is a 5-tuple
$(\mathbb{E}, \omega_{\mathbb{E}}, H, F, W)$, where $(\mathbb{E}, \omega_{\mathbb{E}}, H )$ is a
Hamiltonian system, and the function $H: \mathbb{E} \rightarrow \mathbb{R}$
is called the Hamiltonian, a fiber-preserving map $F: \mathbb{E}\rightarrow
\mathbb{E}$ is called the (external) force map, and a fiber sub-manifold
$W$ of $\mathbb{E}$ is called the control subset.
\end{defi}
Sometimes, $W$ is also denoted the set of fiber-preserving
maps from $\mathbb{E}$ to $W$. When a feedback control law $u:
\mathbb{E}\rightarrow W$ is chosen, the 5-tuple $(\mathbb{E}, \omega_{\mathbb{E}}, H,
F, u)$ is a closed-loop dynamical system. In particular, when $Q$
is a smooth manifold, and $T^\ast Q$ its cotangent bundle with a
symplectic form $\omega$ (not necessarily canonical symplectic
form), then $(T^\ast Q, \omega )$ is a symplectic vector bundle. If
we take that $\mathbb{E}= T^* Q$, from above definition we can obtain an RCH
system on the cotangent bundle $T^\ast Q$, that is, 5-tuple $(T^\ast
Q, \omega, H, F, W)$. Here for convenience, we assume that all
controls appearing in this paper are the admissible controls.\\

In order to describe the dynamics of an RCH system, we have to give a
good expression of the dynamical vector field of the RCH system, by using
the notations of vertical lifted maps of a vector along a fiber, see
Marsden et al. \cite{mawazh10}. For a given RCH System $(T^\ast Q, \omega, H, F, W)$,
the dynamical vector field $X_H$ of the associated Hamiltonian system $(T^\ast Q,
\omega, H) $ satisfies $\mathbf{i}_{X_H}\omega=\mathbf{d}H $.
If considering the external force $F: T^*Q \rightarrow T^*Q, $
which is a fiber-preserving map, by using the
notation of vertical lift map of a vector along a fiber, the
change of $X_H$ under the action of $F$ is that
$$\textnormal{vlift}(F)X_H(\alpha_x)= \textnormal{vlift}((TFX_H)(F(\alpha_x)), \alpha_x)
= (TFX_H)^v_\gamma(\alpha_x),$$
where $\alpha_x \in T^*_x Q, \; x\in Q $ and the geodesic $\gamma$ is a straight
line in $T^*_x Q$ connecting $F_x(\alpha_x)$ and $\alpha_x$. In the
same way, when a feedback control law $u: T^\ast Q \rightarrow W,$
which is a fiber-preserving map, is chosen,
the change of $X_H$ under the action of $u$ is that
$$\textnormal{vlift}(u)X_H(\alpha_x)= \textnormal{vlift}((TuX_H)(F(\alpha_x)), \alpha_x)
= (TuX_H)^v_\gamma(\alpha_x).$$
In consequence, we can give an expression of the dynamical vector
field of the RCH system as follows.
\begin{theo}
The dynamical vector field of an RCH system $(T^\ast
Q,\omega,H,F,W)$ with a control law $u$ is the synthetic
of Hamiltonian vector field $X_H$ and its changes under the actions
of the external force $F$ and control $u$, that is,
$$X_{(T^\ast Q,\omega,H,F,u)}(\alpha_x)= X_H(\alpha_x)+ \textnormal{vlift}(F)X_H(\alpha_x)
+ \textnormal{vlift}(u)X_H(\alpha_x),$$ for any $\alpha_x \in T^*_x
Q, \; x\in Q $. For convenience, it is simply written as
\begin{equation}X_{(T^\ast Q,\omega,H,F,u)}
=X_H +\textnormal{vlift}(F) +\textnormal{vlift}(u). \label{2.1}
\end{equation}
\end{theo}
Where $\textnormal{vlift}(F)=\textnormal{vlift}(F)X_H$,
and $\textnormal{vlift}(u)=\textnormal{vlift}(u)X_H.$ are the
changes of $X_H$ under the actions of $F$ and $u$.
We also denote that $\textnormal{vlift}(W)= \bigcup\{\textnormal{vlift}(u)X_H |
\; u\in W\}$. \\

From the expression (2.1) of the dynamical vector
field of an RCH system, we know that under the actions of the external force $F$
and control $u$, in general, the dynamical vector
field is not Hamiltonian, and hence the RCH system is not
yet a Hamiltonian system. However,
it is a dynamical system closed relative to a
Hamiltonian system, and it can be explored and studied by extending
the methods for external force and control
in the study of Hamiltonian system.\\

In order to describe the nonholonomic RCH system,
in the following we first give the completeness and regularity
conditions for nonholonomic constraint of a mechanical system,
see Le\'{o}n and Wang in \cite{lewa15}. In fact,
in order to describe the dynamics of a nonholonomic mechanical system,
we need some restriction conditions for nonholonomic constraints of
the system. At first, we note that the set of Hamiltonian vector fields
forms a Lie algebra with respect to the Lie bracket, since
$X_{\{f,g\}}=-[X_f, X_g]. $ But, the Lie bracket operator, in
general, may not be closed on the restriction of a nonholonomic
constraint. Thus, we have to give the following completeness
condition for nonholonomic constraints of a system.\\

{\bf $\mathcal{D}$-completeness } Let $Q$ be a smooth manifold and
$TQ$ its tangent bundle. A distribution $\mathcal{D} \subset TQ$ is
said to be {\bf completely nonholonomic} (or bracket-generating) if
$\mathcal{D}$ along with all of its iterated Lie brackets
$[\mathcal{D},\mathcal{D}], [\mathcal{D}, [\mathcal{D},\mathcal{D}]],
\cdots ,$ spans the tangent bundle $TQ$. Moreover, we consider a
nonholonomic mechanical system on $Q$, which is
given by a Lagrangian function $L: TQ \rightarrow \mathbb{R}$
subject to constraint determined by a nonholonomic
distribution $\mathcal{D}\subset TQ$ on the configuration manifold $Q$.
Then the nonholonomic system is said to be {\bf completely nonholonomic},
if the distribution $\mathcal{D} \subset TQ$ determined by the nonholonomic
constraint is completely nonholonomic.\\

{\bf $\mathcal{D}$-regularity } In the following we always assume
that $Q$ is an $n$-dimensional smooth manifold with coordinates $(q^i)$, and $TQ$ its
tangent bundle with coordinates $(q^i,\dot{q}^i)$, and $T^\ast Q$
its cotangent bundle with coordinates $(q^i,p_j)$, which are the
canonical cotangent coordinates of $T^\ast Q$ and $\omega=
dq^{i}\wedge dp_{i}$ is canonical symplectic form on $T^{\ast}Q$. If
the Lagrangian $L: TQ \rightarrow \mathbb{R}$ is hyperregular, that
is, the Hessian matrix
$(\partial^2L/\partial\dot{q}^i\partial\dot{q}^j)$ is nondegenerate
everywhere, then the Legendre transformation $FL: TQ \rightarrow T^*
Q$ is a diffeomorphism. In this case the Hamiltonian $H: T^* Q
\rightarrow \mathbb{R}$ is given by $H(q,p)=\dot{q}\cdot
p-L(q,\dot{q}) $ with Hamiltonian vector field $X_H$,
which is defined by the Hamilton's equation
$\mathbf{i}_{X_H}\omega=\mathbf{d}H$, and
$\mathcal{M}=\mathcal{F}L(\mathcal{D})$ is a constraint submanifold
in $T^* Q$. In particular, for the nonholonomic constraint
$\mathcal{D}\subset TQ$, the Lagrangian $L$ is said to be {\bf
$\mathcal{D}$-regular}, if the restriction of Hessian matrix
$(\partial^2L/\partial\dot{q}^i\partial\dot{q}^j)$ on $\mathcal{D}$
is nondegenerate everywhere. Moreover, a nonholonomic system is said
to be {\bf $\mathcal{D}$-regular}, if its Lagrangian $L$ is {\bf
$\mathcal{D}$-regular}. Note that the restriction of a positive
definite symmetric bilinear form to a subspace is also positive
definite, and hence nondegenerate. Thus, for a simple nonholonomic
mechanical system, that is, whose Lagrangian is the total kinetic
energy minus potential energy, it is {\bf $\mathcal{D}$-regular}
automatically.\\

A nonholonomic RCH system is a 6-tuple $(T^\ast Q,\omega,\mathcal{D},H,F,W)$,
which is an RCH system with a
$\mathcal{D}$-completely and $\mathcal{D}$-regularly nonholonomic
constraint $\mathcal{D} \subset TQ$. Under
the restriction given by constraints, in general, the dynamical
vector field of a nonholonomic RCH system may not be
Hamiltonian vector field. However the nonholonomic RCH system is
also a dynamical system closely related to a Hamiltonian system.
In the following we shall derive a distributional RCH system of the nonholonomic
RCH system $(T^*Q,\omega,\mathcal{D},H,F,W)$,
by analyzing carefully the structure for the nonholonomic
dynamical vector field similar to the method used in Le\'{o}n and Wang \cite{lewa15}.\\

We consider that the constraint submanifold
$\mathcal{M}=\mathcal{F}L(\mathcal{D})\subset T^*Q$ and
$i_{\mathcal{M}}: \mathcal{M}\rightarrow T^*Q $ is the inclusion,
the symplectic form $\omega_{\mathcal{M}}= i_{\mathcal{M}}^* \omega $,
is induced from the canonical symplectic form $\omega$ on $T^* Q$.
We define the distribution $\mathcal{F}$ as the pre-image of the nonholonomic
constraint $\mathcal{D}$ for the map $T\pi_Q: TT^* Q \rightarrow TQ$,
that is, $\mathcal{F}=(T\pi_Q)^{-1}(\mathcal{D})\subset TT^*Q,
$ which is a distribution along $\mathcal{M}$, and
$\mathcal{F}^\circ:=\{\alpha \in T^*T^*Q | <\alpha,v>=0, \; \forall
v\in TT^*Q \}$ is the annihilator of $\mathcal{F}$ in
$T^*T^*Q_{|\mathcal{M}}$. We consider the following nonholonomic
constraints condition
\begin{align} (\mathbf{i}_X \omega -\mathbf{d}H) \in \mathcal{F}^\circ,
\;\;\;\;\;\; X \in T \mathcal{M},
\label{2.2} \end{align} from Cantrijn et al.
\cite{calemama99}, we know that there exists an unique nonholonomic
vector field $X_n$ satisfying the above condition $(2.2)$, if the
admissibility condition $\mathrm{dim}\mathcal{M}=
\mathrm{rank}\mathcal{F}$ and the compatibility condition
$T\mathcal{M}\cap \mathcal{F}^\bot= \{0\}$ hold, where
$\mathcal{F}^\bot$ denotes the symplectic orthogonal of
$\mathcal{F}$ with respect to the canonical symplectic form
$\omega$ on $T^*Q$. In particular, when we consider the Whitney sum
decomposition $T(T^*Q)_{|\mathcal{M}}=T\mathcal{M}\oplus
\mathcal{F}^\bot$ and the canonical projection $P:
T(T^*Q)_{|\mathcal{M}} \rightarrow T\mathcal{M}$,
then we have that $X_n= P(X_H)$.\\

From the condition (2.2) we know that the nonholonomic vector field,
in general, may not be Hamiltonian, because of the restriction
of nonholonomic constraint. But, we hope to study the dynamical
vector field of nonholonomic RCH system by using the similar
method of studying Hamiltonian vector field.
From Le\'{o}n and Wang \cite{lewa15} and
Bates and $\acute{S}$niatycki \cite{basn93}, we can define the
distribution $ \mathcal{K}=\mathcal {F}\cap T\mathcal{M}.$ and
$\mathcal{K}^\bot=\mathcal {F}^\bot\cap T\mathcal{M}, $ where
$\mathcal{K}^\bot$ denotes the symplectic orthogonal of
$\mathcal{K}$ with respect to the canonical symplectic form
$\omega$. If the admissibility condition $\mathrm{dim}\mathcal{M}=
\mathrm{rank}\mathcal{F}$ and the compatibility condition
$T\mathcal{M}\cap \mathcal{F}^\bot= \{0\}$ hold, then we know that the
restriction of the symplectic form $\omega_{\mathcal{M}}$ on
$T^*\mathcal{M}$ fibrewise to the distribution $\mathcal{K}$, that
is, $\omega_\mathcal{K}= \tau_{\mathcal{K}}\cdot
\omega_{\mathcal{M}}$ is non-degenerate, where $\tau_{\mathcal{K}}$
is the restriction map to distribution $\mathcal{K}$. It is worthy
of noting that $\omega_\mathcal{K}$ is not a true two-form on a
manifold, so it does not make sense to speak about it being closed.
We call $\omega_\mathcal{K}$ as a distributional two-form to avoid
any confusion. Because $\omega_\mathcal{K}$ is non-degenerate as a
bilinear form on each fibre of $\mathcal{K}$, there exists a vector
field $X_{\mathcal{K}}$ on $\mathcal{M}$ which takes values in the
constraint distribution $\mathcal{K}$,
such that the distributional Hamiltonian equation holds, that is,
\begin{align}
\mathbf{i}_{X_\mathcal{K}}\omega_{\mathcal{K}}=\mathbf{d}H_\mathcal {K}
\label{2.3} \end{align}
where $\mathbf{d}H_\mathcal{K}$ is the restriction of
$\mathbf{d}H_\mathcal{M}$ to $\mathcal{K}$.
and the function $H_{\mathcal{K}}$ satisfies
$\mathbf{d}H_{\mathcal{K}}= \tau_{\mathcal{K}}\cdot \mathbf{d}H_{\mathcal {M}}$,
and $H_\mathcal{M}= \tau_{\mathcal{M}}\cdot H$ is the restriction of $H$ to
$\mathcal{M}$. Moreover, from the distributional Hamiltonian equation (2.3),
we have that $X_{\mathcal{K}}= \tau_{\mathcal{K}}\cdot X_H.$\\

Moreover, if considering the external force $F$ and control subset $W$,
and define $F_\mathcal{K}=\tau_{\mathcal{K}}\cdot \textnormal{vlift}(F_{\mathcal{M}})X_H,$
and for a control law $u\in W$,
$u_\mathcal{K}= \tau_{\mathcal{K}}\cdot  \textnormal{vlift}(u_{\mathcal{M}})X_H,$
where $F_\mathcal{M}= \tau_{\mathcal{M}}\cdot F$ and
$u_\mathcal{M}= \tau_{\mathcal{M}}\cdot u$ are the restrictions of
$F$ and $u$ to $\mathcal{M}$, that is, $F_\mathcal{K}$ and $u_\mathcal{K}$
are the restrictions of the changes of Hamiltonian vector field $X_H$
under the actions of $F_\mathcal{M}$ and $u_\mathcal{M}$ to $\mathcal{K}$.
Then the 5-tuple $(\mathcal{K},\omega_{\mathcal{K}},
H_\mathcal{K}, F_\mathcal{K}, u_\mathcal{K})$
is a distributional RCH system of the nonholonomic
RCH system $(T^*Q,\omega,\mathcal{D},H,F,W)$ with a control law $u\in W$.
Thus, the geometric formulation of the distributional RCH
system may be summarized as follows.\\

\begin{defi} (Distributional RCH System)
Assume that the 6-tuple $(T^*Q,\omega,\mathcal{D},H,F,W)$ is a nonholonomic
RCH system, where $\omega$ is the canonical
symplectic form on $T^* Q$, and $\mathcal{D}\subset TQ$ is a
$\mathcal{D}$-completely and $\mathcal{D}$-regularly nonholonomic
constraint of the system, and the external force
$F: T^*Q\rightarrow T^*Q$ is the fiber-preserving map,
and the control subset $W\subset T^*Q$ is a fiber submanifold of $T^*Q$.
For a control law $u\in W,$ if there exist a distribution
$\mathcal{K}$, an associated non-degenerate distributional two-form
$\omega_\mathcal{K}$ induced by the canonical symplectic form
and a vector field $X_\mathcal {K}$ on the
constraint submanifold $\mathcal{M}=\mathcal{F}L(\mathcal{D})\subset
T^*Q$, such that the distributional Hamiltonian equation
$\mathbf{i}_{X_\mathcal{K}}\omega_\mathcal{K}=\mathbf{d}H_\mathcal
{K}$ holds, where $\mathbf{d}H_\mathcal{K}$ is the restriction of
$\mathbf{d}H_\mathcal{M}$ to $\mathcal{K}$, and
the function $H_{\mathcal{K}}$ satisfies
$\mathbf{d}H_{\mathcal{K}}= \tau_{\mathcal{K}}\cdot \mathbf{d}H_{\mathcal {M}},$
and $F_{\mathcal{K}}=\tau_{\mathcal{K}}\cdot \textnormal{vlift}(F_{\mathcal{M}})X_H$,
and $u_{\mathcal{K}}=\tau_{\mathcal{K}}\cdot \textnormal{vlift}(u_{\mathcal{M}})X_H$
as defined above,
then the 5-tuple $(\mathcal{K},\omega_{\mathcal{K}},
H_{\mathcal{K}}, F_{\mathcal{K}}, u_{\mathcal{K}})$
is called a distributional RCH system of the nonholonomic
RCH system $(T^*Q,\omega,\mathcal{D},H,F,u)$, and $X_\mathcal
{K}$ is called a nonholonomic dynamical vector field. Denote by
\begin{align}
\tilde{X}=X_{(\mathcal{K},\omega_{\mathcal{K}},
H_{\mathcal{K}}, F_{\mathcal{K}}, u_{\mathcal{K}})}
=X_\mathcal {K}+ F_{\mathcal{K}}+u_{\mathcal{K}}
\label{2.4} \end{align}
is the dynamical vector field of the distributional RCH system
$(\mathcal{K},\omega_{\mathcal {K}}, H_{\mathcal{K}}, F_{\mathcal{K}}, u_{\mathcal{K}})$,
which is the synthetic
of the nonholonomic dynamical vector field $X_{\mathcal{K}}$ and
the vector fields $F_{\mathcal{K}}$ and $u_{\mathcal{K}}$.
Under the above circumstances, we refer to
$(T^*Q,\omega,\mathcal{D}, H, F, u)$ as a nonholonomic RCH system
with an associated distributional RCH system
$(\mathcal{K},\omega_{\mathcal {K}}, H_{\mathcal{K}}, F_{\mathcal{K}}, u_{\mathcal{K}})$.
\end{defi}

It is worthy of noting that,
if the external force and control of a distributional RCH
system $(\mathcal{K},\omega_{\mathcal {K}},H_{\mathcal {K}},
F_{\mathcal {K}}, u_{\mathcal {K}})$ are both zeros,
that is, $F_{\mathcal {K}}=0 $ and $u_{\mathcal {K}}=0$, in this case,
the distributional RCH system is just a distributional Hamiltonian system
$(\mathcal{K},\omega_{\mathcal {K}},H_{\mathcal {K}})$,
which is given in Le\'{o}n and Wang \cite{lewa15}.
Thus, the distributional RCH system
can be regarded as an extension of the distributional Hamiltonian system to
the system with external force and control.
Moreover, in section 4 and section 5, we consider the nonholonomic RCH system
with symmetry, as well as momentum map.
By using the similar method for nonholonomic reduction
given in Bates and $\acute{S}$niatycki \cite{basn93},
Le\'{o}n and Wang \cite{lewa15} and Wang \cite{wa21c},
and analyzing carefully the structures for the nonholonomic reduced
dynamical vector fields, we also give the geometric formulations of
the nonholonomic reduced distributional RCH systems.

\section{Hamilton-Jacobi Equations for a Distributional RCH System }

In this section, for a nonholonomic RCH system
$(T^*Q,\omega,\mathcal{D},H, F, W)$
with an associated distributional RCH system
$(\mathcal{K},\omega_{\mathcal {K}}, H_{\mathcal{K}}, F_{\mathcal{K}}, u_{\mathcal{K}})$,
we shall derive precisely the geometric constraint conditions
of the distributional two-form $\omega_\mathcal{K}$
for the dynamical vector field of the distributional RCH system,
that is, the two types of Hamilton-Jacobi equation
for the distributional RCH system.
In order to do this, in the following we first give
some important notions and prove a key lemma, which is an important
tool for the proofs of two types of
Hamilton-Jacobi theorem for the distributional RCH system.\\

Let $Q$ be an $n$-dimensional smooth manifold, and
denote by $\Omega^i(Q)$ the set of all i-forms on $Q$, $i=1,2.$
For any $\gamma \in \Omega^1(Q),\; q\in Q, $ then $\gamma(q)\in T_q^*Q, $
and we can define a map $\gamma: Q \rightarrow T^*Q, \; q \rightarrow (q, \gamma(q)).$
Hence we say often that the map $\gamma: Q
\rightarrow T^*Q$ is an one-form on $Q$.
Assume that $\omega$ is the canonical symplectic form on $T^*Q$,
and $\mathcal{D}\subset TQ$ is a $\mathcal{D}$-regularly nonholonomic
constraint, and the projection $\pi_Q: T^* Q \rightarrow Q $ induces
the map $T\pi_{Q}: TT^* Q \rightarrow TQ. $
If the one-form $\gamma$ is closed,
then $\mathbf{d}\gamma(x,y)=0, \; \forall\; x, y \in TQ$.
Note that for any $v, w \in TT^* Q, $ we have that
$\mathbf{d}\gamma(T\pi_{Q}(v),T\pi_{Q}(w))=\pi^*(\mathbf{d}\gamma )(v, w)$
is a two-form on the cotangent bundle $T^*Q$, where
$\pi^*: T^*Q \rightarrow T^*T^*Q.$ Thus,
in the following we can introduce two weaker notions.

\begin{defi}
\noindent $(\mathrm{i})$ The one-form $\gamma$ is called to be closed with respect to $T\pi_{Q}:
TT^* Q \rightarrow TQ, $ if for any $v, w \in TT^* Q, $ we have that
$\mathbf{d}\gamma(T\pi_{Q}(v),T\pi_{Q}(w))=0; $\\

\noindent $(\mathrm{ii})$ The one-form $\gamma$ is called to be closed
on $\mathcal{D}$ with respect to $T\pi_{Q}:
TT^* Q \rightarrow TQ, $ if for any $v, w \in TT^* Q, $
and $T\pi_{Q}(v), \; T\pi_{Q}(w) \in \mathcal{D},$  we have that
$\mathbf{d}\gamma(T\pi_{Q}(v),T\pi_{Q}(w))=0. $
\end{defi}

From the above Definition 3.1 we know that, the notion that
$\gamma$ is closed on $\mathcal{D}$ with respect to $T\pi_{Q}:
TT^* Q \rightarrow TQ, $ is weaker than the notion that $\gamma$ is closed
with respect to $T\pi_{Q}: TT^* Q \rightarrow TQ. $
From Wang \cite{wa17} we also know that the latter, that is, $\gamma$ is closed
with respect to $T\pi_{Q}: TT^* Q \rightarrow TQ, $
is weaker than the notion that $\gamma$ is closed. Thus,
the notion that $\gamma$ is closed on $\mathcal{D}$ with respect to $T\pi_{Q}:
TT^* Q \rightarrow TQ, $ is weaker than that $\gamma$ is closed on $\mathcal{D}$,
that is, $\mathbf{d}\gamma(x,y)=0, \; \forall\; x, y \in \mathcal{D}$.
In fact, if $\gamma$ is a closed one-form on $\mathcal{D}$,
then it must be closed on $\mathcal{D}$ with respect to
$T\pi_{Q}: TT^* Q \rightarrow TQ. $
Conversely, if $\gamma$ is closed on $\mathcal{D}$ with respect to
$T\pi_{Q}: TT^* Q \rightarrow TQ, $ then it may not be closed on $\mathcal{D}$.
We can prove a general result as follows, which states that
the notion that $\gamma$ is closed on $\mathcal{D}$
with respect to $T\pi_{Q}: TT^* Q \rightarrow TQ, $
is not equivalent to the notion that $\gamma$ is closed on $\mathcal{D}$,
also see Le\'{o}n and Wang \cite{lewa15}.

\begin{prop}
Assume that $\gamma: Q \rightarrow T^*Q$ is an one-form on $Q$ and
it is not closed on $\mathcal{D}$. We define the set $N$, which is a subset of $TQ$,
such that the one-form $\gamma$ on $N$ satisfies the condition that
for any $x,y \in N, \; \mathbf{d}\gamma(x,y)\neq 0. $ Denote
$Ker(T\pi_Q)= \{u \in TT^*Q| \; T\pi_Q(u)=0 \}, $ and $T\gamma: TQ
\rightarrow TT^* Q .$ If $T\gamma(N)\subset Ker(T\pi_Q), $ then
$\gamma$ is closed with respect to $T\pi_{Q}: TT^* Q \rightarrow TQ.$
and hence $\gamma$ is closed on $\mathcal{D}$ with respect to
$T\pi_{Q}: TT^* Q \rightarrow TQ.$
\end{prop}

\noindent{\bf Proof: } If the $\gamma: Q \rightarrow T^*Q$
is not closed on $\mathcal{D}$, then it is not yet closed itself.
For any $v, w \in TT^* Q, $ if
$T\pi_{Q}(v) \notin N, $ or $T\pi_{Q}(w))\notin N, $ then by the
definition of $N$, we know that
$\mathbf{d}\gamma(T\pi_{Q}(v),T\pi_{Q}(w))=0; $ If $T\pi_{Q}(v)\in
N, $ and $T\pi_{Q}(w))\in N, $ from the condition $T\gamma(N)\subset
Ker(T\pi_Q), $ we know that $T\pi_{Q}\cdot T\gamma \cdot
T\pi_{Q}(v)= T\pi_{Q}(v)=0, $ and $T\pi_{Q}\cdot T\gamma \cdot
T\pi_{Q}(w)= T\pi_{Q}(w)=0, $ where we have used the
relation $\pi_Q\cdot \gamma\cdot \pi_Q= \pi_Q, $ and hence
$\mathbf{d}\gamma(T\pi_{Q}(v),T\pi_{Q}(w))=0. $ Thus, for any $v, w
\in TT^* Q, $ we have always that
$\mathbf{d}\gamma(T\pi_{Q}(v),T\pi_{Q}(w))=0. $
In particular, for any $v, w \in TT^* Q, $
and $T\pi_{Q}(v), \; T\pi_{Q}(w) \in \mathcal{D},$  we have
$\mathbf{d}\gamma(T\pi_{Q}(v),T\pi_{Q}(w))=0. $
that is, $\gamma$ is closed on $\mathcal{D}$
with respect to $T\pi_{Q}: TT^* Q \rightarrow TQ. $
\hskip 0.3cm $\blacksquare$\\

Now, we prove the following Lemma 3.3. It is worthy of noting that
this lemma is an extension of Lemma 2.4 given in Wang \cite{wa17} to the
nonholonomic context. The proofs of the following $(\mathrm{i})$ and $(\mathrm{ii})$
are given in Wang \cite{wa17}, and $(\mathrm{iii})$ given in Le\'{o}n and Wang \cite{lewa15}.
This lemma offers also an important tool for the proofs of the two types of Hamilton-Jacobi
theorems for the distributional RCH system and the nonholonomic
reduced distributional RCH system, and hence we also give its proof here.

\begin{lemm}
Assume that $\gamma: Q \rightarrow T^*Q$ is an one-form on $Q$, and
$\lambda=\gamma \cdot \pi_{Q}: T^* Q \rightarrow T^* Q .$ Then
we have that\\
\noindent $\bf (\mathrm{i})$ for any $x, y \in TQ, \;
\gamma^*\omega (x,y)= -\mathbf{d}\gamma (x,y),$ and for any $v, w \in
TT^* Q, $ \\ $ \lambda^*\omega(v,w)=
-\mathbf{d}\gamma(T\pi_{Q}(v), \; T\pi_{Q}(w)),$
since $\omega$ is the canonical symplectic form on $T^*Q$; \\

\noindent $\bf (\mathrm{ii})$ for any $v, w \in TT^* Q, \;
\omega(T\lambda \cdot v,w)= \omega(v, w-T\lambda \cdot
w)-\mathbf{d}\gamma(T\pi_{Q}(v), \; T\pi_{Q}(w))$ ;\\

\noindent $\bf (\mathrm{iii})$ If the Lagrangian $L$ is $\mathcal{D}$-regular, and
$\textmd{Im}(\gamma)\subset \mathcal{M}=\mathcal{F}L(\mathcal{D}), $
then we have that $ X_{H}\cdot \gamma \in \mathcal{F}$ along
$\gamma$, and $ X_{H}\cdot \lambda \in \mathcal{F}$ along
$\lambda$, that is, $T\pi_{Q}(X_H\cdot\gamma(q))\in
\mathcal{D}_{q}, \; \forall q \in Q $, and $T\pi_{Q}(X_H\cdot\lambda(q,p))\in
\mathcal{D}_{q}, \; \forall q \in Q, \; (q,p) \in T^* Q. $
Moreover, if a symplectic map $\varepsilon: T^* Q \rightarrow T^* Q $
satisfies the condition $\varepsilon(\mathcal{M})\subset \mathcal{M},$ then
we have that $ X_{H}\cdot \varepsilon \in \mathcal{F}$ along
$\varepsilon. $
\end{lemm}

\noindent{\bf Proof:}
We first prove the assertion $(\mathrm{i})$.
Since $\omega$ is the canonical symplectic form on $T^*Q$, we know
that there is an unique canonical one-form $\theta$, such that
$\omega= -\mathbf{d} \theta. $ From the Proposition 3.2.11 in
Abraham and Marsden \cite{abma78}, we have that for the one-form
$\gamma: Q \rightarrow T^*Q, \; \gamma^* \theta= \gamma. $ Then we
can obtain that
\begin{align*}
\gamma^*\omega(x,y) = \gamma^* (-\mathbf{d} \theta) (x, y) =
-\mathbf{d}(\gamma^* \theta)(x, y)= -\mathbf{d}\gamma (x, y).
\end{align*}
Note that $\lambda=\gamma \cdot \pi_{Q}: T^* Q \rightarrow T^* Q, $
and $\lambda^*= \pi_{Q}^* \cdot \gamma^*: T^*T^* Q \rightarrow
T^*T^* Q, $ then we have that
\begin{align*}
\lambda^*\omega(v,w) &= \lambda^* (-\mathbf{d} \theta) (v, w)
=-\mathbf{d}(\lambda^* \theta)(v, w)= -\mathbf{d}(\pi_{Q}^* \cdot
\gamma^* \theta)(v, w)\\ &= -\mathbf{d}(\pi_{Q}^* \cdot\gamma )(v,
w)= -\mathbf{d}\gamma(T\pi_{Q}(v), \; T\pi_{Q}(w)).
\end{align*}
It follows that the assertion $(\mathrm{i})$ holds.\\

Next, we prove the assertion $(\mathrm{ii})$. For any $v, w \in TT^*
Q,$ note that $v- T(\gamma \cdot \pi_Q)\cdot v$ is vertical, because
$$
T\pi_Q(v- T(\gamma \cdot \pi_Q)\cdot v)=T\pi_Q(v)-T(\pi_Q\cdot
\gamma\cdot \pi_Q)\cdot v= T\pi_Q(v)-T\pi_Q(v)=0,
$$
where we have used the relation $\pi_Q\cdot \gamma\cdot \pi_Q= \pi_Q. $
Thus, $\omega(v- T(\gamma \cdot \pi_Q)\cdot v,w- T(\gamma \cdot
\pi_Q)\cdot w)= 0, $ and hence,
$$\omega(T(\gamma \cdot \pi_Q)\cdot v, \; w)=
\omega(v, \; w-T(\gamma \cdot \pi_Q)\cdot w)+ \omega(T(\gamma \cdot
\pi_Q)\cdot v, \; T(\gamma \cdot \pi_Q)\cdot w). $$ However, the
second term on the right-hand side is given by
$$
\omega(T(\gamma \cdot \pi_Q)\cdot v, \; T(\gamma \cdot \pi_Q)\cdot
w)= \gamma^*\omega(T\pi_Q(v), \; T\pi_Q(w))=
-\mathbf{d}\gamma(T\pi_{Q}(v), \; T\pi_{Q}(w)),
$$
where we have used the assertion $(\mathrm{i})$. It follows that
\begin{align*}
\omega(T\lambda \cdot v,w) &=\omega(T(\gamma \cdot \pi_Q)\cdot v, \;
w)\\ &= \omega(v, \; w-T(\gamma \cdot \pi_Q)\cdot w)-\mathbf{d}\gamma(T\pi_{Q}(v), \; T\pi_{Q}(w))
\\ &= \omega(v,
w-T\lambda \cdot w)-\mathbf{d}\gamma(T\pi_{Q}(v), \; T\pi_{Q}(w)).
\end{align*}
Thus, the assertion $(\mathrm{ii})$ holds.\\

Now, we prove $(\mathrm{iii})$. Under the canonical cotangent bundle coordinates,
for any $q \in Q, \; (q,p)\in T^* Q, $ we have that
$$
X_H\cdot \gamma(q)=\sum^n_{i=1} (\frac{\partial H}{\partial
p_i}\frac{\partial}{\partial q^i}-\frac{\partial H}{\partial
q^i}\frac{\partial}{\partial p_i})\gamma(q).
$$
and
$$
X_H\cdot \lambda(q,p)=\sum^n_{i=1} (\frac{\partial H}{\partial
p_i}\frac{\partial}{\partial q^i}-\frac{\partial H}{\partial
q^i}\frac{\partial}{\partial p_i})\gamma\cdot \pi_Q(q,p).
$$
Then,
$$
T\pi_Q(X_H\cdot \gamma(q))=T\pi_Q(X_H\cdot \lambda(q,p))
= \sum^n_{i=1}(\frac{\partial H}{\partial p_i}\frac{\partial}{\partial q^i})
\gamma(q) \in T_q Q.
$$
Since $\textmd{Im}(\gamma)\subset \mathcal{M}, $ and
$\gamma(q)\in \mathcal{M}_{(q,p)}=\mathcal{F}L(\mathcal{D}_q), $ from the Lagrangian $L$ is
$\mathcal{D}$-regular, and $\mathcal{F}L$ is a diffeomorphism, then
there exists a point $(q, \; v_q)\in \mathcal{D}_q, $ such that
$\mathcal{F}L(q,\; v_q)=\gamma(q). $ Thus,
$$
T\pi_Q(X_H\cdot \gamma(q))=T\pi_Q(X_H\cdot \lambda(q,p))
=\mathcal{F}L(q,\; v_q)  \sum^n_{i=1}(\frac{\partial H}{\partial p_i}\frac{\partial}{\partial q^i}) \in
\mathcal{D}_q,
$$
it follows that $ X_{H}\cdot \gamma \in \mathcal{F}$ along
$\gamma$, and $ X_{H}\cdot \lambda \in \mathcal{F}$ along $\lambda$.
Moreover, for the symplectic map $\varepsilon: T^* Q \rightarrow T^* Q $,
we have that
$$
X_H\cdot \varepsilon (q,p)= (\sum^n_{i=1} (\frac{\partial H}{\partial
p_i}\frac{\partial}{\partial q^i} - \frac{\partial H}{\partial
q^i}\frac{\partial}{\partial p_i}))\varepsilon (q,p).
$$
If $\varepsilon$ satisfies the
condition $\varepsilon(\mathcal{M})\subset \mathcal{M},$
then for any $(q,p)\in \mathcal{M}_{(q,p)}$, we have that
$\varepsilon(q,p)\in \mathcal{M}_{(q,p)},$
and there exists a point $(q,\; v_q)\in \mathcal{D}_q, $ such that
$\mathcal{F}L(q,\; v_q)=\varepsilon (q,p). $ Thus,
$$
T\pi_Q(X_H\cdot \varepsilon(q, p))
= \sum^n_{i=1}(\frac{\partial H}{\partial p_i}\frac{\partial}{\partial q^i})
\varepsilon(q, p)
=\mathcal{F}L(q, \; v_q)\sum^n_{i=1}(\frac{\partial H}{\partial
p_i}\frac{\partial}{\partial q^i}) \in
\mathcal{D}_q,
$$
it follows that $ X_{H}\cdot \varepsilon \in \mathcal{F}$ along
$\varepsilon$.
\hskip 0.3cm $\blacksquare$\\

We note that the distributional RCH system
is determined by a non-degenerate distributional two-form
induced from the canonical symplectic form, but, the distributional two-form
is not a "true two-form" on a manifold, and hence the leading
distributional RCH system can not be Hamiltonian.
Thus, we can not describe the Hamilton-Jacobi equations for the distributional
RCH system from the viewpoint of generating function
as in the classical Hamiltonian case, that is,
we cannot prove the Hamilton-Jacobi
theorem for the distributional RCH system,
just like same as the above Theorem 1.1.
Since the distributional RCH system is a
dynamical system closely related to a Hamiltonian system,
in the following by using Lemma 3.3,
for a given nonholonomic RCH system
$(T^*Q,\omega,\mathcal{D},H,F,W)$ with an associated distributional RCH
system $(\mathcal{K},\omega_{\mathcal {K}},H_{\mathcal{K}}, F_{\mathcal{K}}, u_{\mathcal{K}})$,
we can derive precisely the geometric constraint conditions of
the non-degenerate distributional two-form $\omega_\mathcal{K}$
for the dynamical vector field $X_{(\mathcal{K},\omega_{\mathcal{K}},
H_{\mathcal{K}}, F_{\mathcal{K}}, u_{\mathcal{K}})}$,
that is, the two types of Hamilton-Jacobi equation for the distributional RCH
system $(\mathcal{K},\omega_{\mathcal {K}},H_{\mathcal{K}},
F_{\mathcal{K}}, u_{\mathcal{K}})$.
At first, by using the fact that the one-form $\gamma: Q
\rightarrow T^*Q $ is closed on $\mathcal{D}$ with respect to
$T\pi_Q: TT^* Q \rightarrow TQ, $
$\textmd{Im}(\gamma)\subset \mathcal{M}, $ and
$\textmd{Im}(T\gamma)\subset \mathcal{K}, $
we can prove the Type I of
Hamilton-Jacobi theorem for the distributional RCH system.
For convenience, the
maps involved in the following theorem and its proof are shown in
Diagram-1.

\begin{center}
\hskip 0cm \xymatrix{& \mathcal{M} \ar[d]_{X_{\mathcal{K}}}
\ar[r]^{i_{\mathcal{M}}} & T^* Q \ar[d]_{X_H}
 \ar[r]^{\pi_Q}
& Q \ar[d]_{\tilde{X}^\gamma} \ar[r]^{\gamma} & T^*Q \ar[d]^{\tilde{X}} \\
& \mathcal{K}  & T(T^*Q) \ar[l]_{\tau_{\mathcal{K}}} & TQ
\ar[l]_{T\gamma} & T(T^* Q)\ar[l]_{T\pi_Q}}
\end{center}
$$\mbox{Diagram-1}$$

\begin{theo} (Type I of Hamilton-Jacobi Theorem for a Distributional RCH System)
For the nonholonomic RCH system $(T^*Q,\omega,\mathcal{D},H,F,u)$
with an associated distributional RCH system \\
$(\mathcal{K},\omega_{\mathcal {K}},H_{\mathcal{K}}, F_{\mathcal{K}}, u_{\mathcal{K}})$,
assume that $\gamma: Q \rightarrow T^*Q$ is an one-form on $Q$, and
$\tilde{X}^\gamma = T\pi_{Q}\cdot \tilde{X} \cdot \gamma$,
where $\tilde{X}=X_{(\mathcal{K},\omega_{\mathcal{K}},
H_{\mathcal{K}}, F_{\mathcal{K}}, u_{\mathcal{K}})}
=X_\mathcal {K}+ F_{\mathcal{K}}+u_{\mathcal{K}}$
is the dynamical vector field of the distributional RCH system.
Moreover, assume that $\textmd{Im}(\gamma)\subset
\mathcal{M}=\mathcal{F}L(\mathcal{D}), $ and $
\textmd{Im}(T\gamma)\subset \mathcal{K}. $ If the
one-form $\gamma: Q \rightarrow T^*Q $ is closed on $\mathcal{D}$ with respect to
$T\pi_Q: TT^* Q \rightarrow TQ, $ then $\gamma$ is a
solution of the equation $T\gamma \cdot
\tilde{X}^\gamma= X_{\mathcal{K}} \cdot \gamma, $ where $X_{\mathcal{K}}$
is the nonholonomic dynamical vector field of the distributional RCH system.
The equation $T\gamma \cdot \tilde{X}^\gamma= X_{\mathcal{K}} \cdot \gamma $
is called the Type I of Hamilton-Jacobi equation for the distributional RCH system
$(\mathcal{K},\omega_{\mathcal {K}},H_{\mathcal{K}}, F_{\mathcal{K}}, u_{\mathcal{K}})$.
\end{theo}

\noindent{\bf Proof: } From Definition 2.3 we have that
the dynamical vector field of the distributional RCH system
$(\mathcal{K},\omega_{\mathcal {K}}, H_{\mathcal{K}}, F_{\mathcal{K}}, u_{\mathcal{K}})$
is the synthetic of the nonholonomic dynamical vector field $X_{\mathcal{K}}$ and
the vector fields $F_{\mathcal{K}}$ and $u_{\mathcal{K}}$, that is,
$\tilde{X}=X_{(\mathcal{K},\omega_{\mathcal{K}},
H_{\mathcal{K}}, F_{\mathcal{K}}, u_{\mathcal{K}})}
=X_\mathcal {K}+ F_{\mathcal{K}}+u_{\mathcal{K}}$,
and $F_{\mathcal{K}}=\tau_{\mathcal{K}}\cdot \textnormal{vlift}(F_{\mathcal{M}})X_H$,
and $u_{\mathcal{K}}=\tau_{\mathcal{K}}\cdot \textnormal{vlift}(u_{\mathcal{M}})X_H$.
Note that $T\pi_{Q}\cdot \textnormal{vlift}(F_{\mathcal{M}})X_H=T\pi_{Q}\cdot \textnormal{vlift}(u_{\mathcal{M}})X_H=0, $
then we have that $T\pi_{Q}\cdot F_{\mathcal{K}}=T\pi_{Q}\cdot u_{\mathcal{K}}=0,$
and hence $T\pi_{Q}\cdot \tilde{X}\cdot \gamma=T\pi_{Q}\cdot X_{\mathcal{K}}\cdot \gamma. $
On the other hand, we note that
$\textmd{Im}(\gamma)\subset \mathcal{M}, $ and
$\textmd{Im}(T\gamma)\subset \mathcal{K}, $ in this case,
$\omega_{\mathcal{K}}\cdot
\tau_{\mathcal{K}}=\tau_{\mathcal{K}}\cdot \omega_{\mathcal{M}}=
\tau_{\mathcal{K}}\cdot i_{\mathcal{M}}^* \cdot \omega, $ along
$\textmd{Im}(T\gamma)$. Moreover, from the distributional Hamiltonian equation (2.3),
we have that $X_{\mathcal{K}}= \tau_{\mathcal{K}}\cdot X_H,$ and
$\tau_{\mathcal{K}}\cdot X_{H}\cdot \gamma =
X_{\mathcal{K}}\cdot \gamma \in \mathcal{K}$.
Thus, using the non-degenerate
distributional two-form $\omega_{\mathcal{K}}$, from Lemma 3.3(ii) and (iii),
if we take that $v= X_{\mathcal{K}}\cdot \gamma \in \mathcal{K} (\subset \mathcal{F}),$
and for any $w \in \mathcal{F}, \; T\lambda(w)\neq 0, $ and
$\tau_{\mathcal{K}}\cdot w \neq 0, $ then we have that
\begin{align*}
& \omega_{\mathcal{K}}(T\gamma \cdot \tilde{X}^\gamma, \;
\tau_{\mathcal{K}}\cdot w)=
\omega_{\mathcal{K}}(\tau_{\mathcal{K}}\cdot T\gamma \cdot
\tilde{X}^\gamma, \; \tau_{\mathcal{K}}\cdot w)\\ & =
\tau_{\mathcal{K}}\cdot i_{\mathcal{M}}^* \cdot\omega(T\gamma \cdot
T\pi_Q \cdot \tilde{X}\cdot \gamma, \; w ) = \tau_{\mathcal{K}}\cdot
i_{\mathcal{M}}^* \cdot \omega(T\gamma \cdot T\pi_Q \cdot X_{\mathcal{K}}\cdot \gamma, \; w)\\
& =\tau_{\mathcal{K}}\cdot i_{\mathcal{M}}^* \cdot \omega
(T\gamma \cdot T\pi_{Q}\cdot X_{\mathcal{K}}\cdot \gamma, \; w)
= \tau_{\mathcal{K}}\cdot i_{\mathcal{M}}^* \cdot \omega
(T(\gamma \cdot \pi_{Q})\cdot X_{\mathcal{K}}\cdot \gamma, \; w)\\
& =\tau_{\mathcal{K}}\cdot i_{\mathcal{M}}^* \cdot (\omega
(X_{\mathcal{K}}\cdot \gamma, \; w-T(\gamma
\cdot \pi_Q)\cdot w)-\mathbf{d}\gamma(T\pi_{Q}(X_{\mathcal{K}}\cdot \gamma), \; T\pi_{Q}(w)))\\
& = \tau_{\mathcal{K}}\cdot i_{\mathcal{M}}^* \cdot\omega
(X_{\mathcal{K}}\cdot \gamma, \; w) - \tau_{\mathcal{K}}\cdot i_{\mathcal{M}}^* \cdot
\omega (X_{\mathcal{K}}\cdot \gamma, \; T(\gamma \cdot \pi_Q) \cdot w)\\
& \;\;\;\;\;\; -\tau_{\mathcal{K}}\cdot i_{\mathcal{M}}^* \cdot \mathbf{d}\gamma
(T\pi_{Q}(X_{\mathcal{K}}\cdot \gamma), \; T\pi_{Q}(w))\\
& = \omega_{\mathcal{K}}( \tau_{\mathcal{K}}\cdot X_{\mathcal{K}}\cdot \gamma,
\; \tau_{\mathcal{K}}\cdot w) -
\omega_{\mathcal{K}}(\tau_{\mathcal{K}}\cdot X_{\mathcal{K}}\cdot \gamma, \;
\tau_{\mathcal{K}}\cdot T(\gamma \cdot \pi_Q) \cdot w)\\
& \;\;\;\;\;\;
-\tau_{\mathcal{K}}\cdot i_{\mathcal{M}}^* \cdot\mathbf{d}\gamma
(T\pi_{Q}(X_{\mathcal{K}}\cdot \gamma), \; T\pi_{Q}(w))\\
& = \omega_{\mathcal{K}}(X_{\mathcal{K}}\cdot \gamma, \;
\tau_{\mathcal{K}} \cdot w) -
\omega_{\mathcal{K}}(X_{\mathcal{K}} \cdot \gamma, \; T\gamma \cdot T\pi_{Q}(w))\\
& \;\;\;\;\;\; - \tau_{\mathcal{K}}\cdot
i_{\mathcal{M}}^* \cdot\mathbf{d}\gamma
(T\pi_{Q}(X_{\mathcal{K}}\cdot \gamma), \; T\pi_{Q}(w)),
\end{align*}
where we have used that $ \tau_{\mathcal{K}}\cdot T\gamma= T\gamma, $ and
$\tau_{\mathcal{K}}\cdot X_{\mathcal{K}}\cdot \gamma = X_{\mathcal{K}}\cdot
\gamma, $ since $\textmd{Im}(T\gamma)\subset \mathcal{K}. $
If the one-form $\gamma: Q \rightarrow T^*Q $ is closed on $\mathcal{D}$ with respect to
$T\pi_Q: TT^* Q \rightarrow TQ, $ then we have that
$\mathbf{d}\gamma(T\pi_{Q}(X_{\mathcal{K}}\cdot \gamma), \; T\pi_{Q}(w))=0, $
since $X_{\mathcal{K}}\cdot \gamma, \; w \in \mathcal{F},$
and $T\pi_{Q}(X_{\mathcal{K}}\cdot \gamma), \; T\pi_{Q}(w) \in \mathcal{D}, $ and hence
$$
\tau_{\mathcal{K}}\cdot i_{\mathcal{M}}^* \cdot\mathbf{d}\gamma
(T\pi_{Q}(X_{\mathcal{K}}\cdot \gamma), \; T\pi_{Q}(w))=0,
$$
and
\begin{equation}
\omega_{\mathcal{K}}(T\gamma \cdot \tilde{X}^\gamma, \;
\tau_{\mathcal{K}}\cdot w)- \omega_{\mathcal{K}}(X_{\mathcal{K}}\cdot \gamma, \;
\tau_{\mathcal{K}} \cdot w)
= -\omega_{\mathcal{K}}( X_{\mathcal{K}} \cdot \gamma, \; T\gamma \cdot T\pi_{Q}(w)).
\label{3.1} \end{equation}
If $\gamma$ satisfies the equation $T\gamma\cdot \tilde{X}^\gamma= X_{\mathcal{K}}\cdot \gamma ,$
from Lemma 3.3(i) we know that the right side of (3.1) becomes that
\begin{align*}
 -\omega_{\mathcal{K}}(X_{\mathcal{K}} \cdot \gamma, \; T\gamma \cdot T\pi_{Q}(w))
& = -\omega_{\mathcal{K}}(T\gamma\cdot \tilde{X}^\gamma, \; T\gamma \cdot T\pi_{Q}(w))\\
& = -\omega_{\mathcal{K}}(\tau_{\mathcal{K}}\cdot T\gamma \cdot T\pi_Q \cdot
\tilde{X}\cdot \gamma, \; \tau_{\mathcal{K}}\cdot T\gamma \cdot T\pi_{Q}(w))\\
& = -\tau_{\mathcal{K}}\cdot
i_{\mathcal{M}}^* \cdot \omega(T\gamma \cdot T\pi_{Q}(X_{\mathcal{K}}\cdot\gamma),
\; T\gamma \cdot T\pi_{Q}(w))\\
& = -\tau_{\mathcal{K}}\cdot
i_{\mathcal{M}}^* \cdot\gamma^*\omega( T\pi_{Q}(X_{\mathcal{K}}\cdot\gamma), \; T\pi_{Q}(w))\\
& = \tau_{\mathcal{K}}\cdot i_{\mathcal{M}}^* \cdot
\mathbf{d}\gamma(T\pi_{Q}( X_{\mathcal{K}}\cdot\gamma ), \; T\pi_{Q}(w))=0.
\end{align*}
Because the distributional two-form $\omega_{\mathcal{K}}$ is non-degenerate,
the left side of (3.1) equals zero, only when
$\gamma$ satisfies the equation
$T\gamma\cdot \tilde{X}^\gamma= X_{\mathcal{K}}\cdot \gamma .$ Thus,
if the one-form $\gamma: Q \rightarrow T^*Q $ is closed on $\mathcal{D}$ with respect to
$T\pi_Q: TT^* Q \rightarrow TQ, $ then $\gamma$ must be a solution of the Type I of Hamilton-Jacobi equation
$T\gamma\cdot \tilde{X}^\gamma= X_{\mathcal{K}}\cdot \gamma .$
\hskip 0.3cm $\blacksquare$\\

Next, for any symplectic map $\varepsilon: T^* Q \rightarrow T^* Q $,
we can prove the following Type II of
Hamilton-Jacobi theorem for the distributional RCH system.
For convenience, the
maps involved in the following theorem and its proof are shown in
Diagram-2.

\begin{center}
\hskip 0cm \xymatrix{& \mathcal{M} \ar[d]_{X_{\mathcal{K}}}
\ar[r]^{i_{\mathcal{M}}} & T^* Q \ar[d]_{X_{H\cdot \varepsilon}}
\ar[dr]^{\tilde{X}^\varepsilon} \ar[r]^{\pi_Q}
& Q \ar[r]^{\gamma} & T^*Q \ar[d]^{\tilde{X}} \\
& \mathcal{K}  & T(T^*Q) \ar[l]_{\tau_{\mathcal{K}}} & TQ
\ar[l]_{T\gamma} & T(T^* Q)\ar[l]_{T\pi_Q}}
\end{center}
$$\mbox{Diagram-2}$$

\begin{theo} (Type II of Hamilton-Jacobi Theorem for a Distributional RCH System)
For the nonholonomic RCH system $(T^*Q,\omega,\mathcal{D},H,F,u)$
with an associated distributional RCH system \\
$(\mathcal{K},\omega_{\mathcal {K}},H_{\mathcal{K}}, F_{\mathcal{K}}, u_{\mathcal{K}})$,
assume that $\gamma: Q \rightarrow T^*Q$ is an one-form on $Q$,
and $\lambda=\gamma \cdot \pi_{Q}: T^* Q \rightarrow T^* Q, $ and for any
symplectic map $\varepsilon: T^* Q \rightarrow T^* Q $, denote by
$\tilde{X}^\varepsilon = T\pi_{Q}\cdot \tilde{X}\cdot \varepsilon$, where
$\tilde{X}=X_{(\mathcal{K},\omega_{\mathcal{K}},
H_{\mathcal{K}}, F_{\mathcal{K}}, u_{\mathcal{K}})}
=X_\mathcal {K}+ F_{\mathcal{K}}+u_{\mathcal{K}}$
is the dynamical vector field of the distributional RCH system.
Moreover, assume that $\textmd{Im}(\gamma)\subset
\mathcal{M}=\mathcal{F}L(\mathcal{D}), $ and $\varepsilon(\mathcal{M})\subset \mathcal{M}$,
and $\textmd{Im}(T\gamma)\subset \mathcal{K}. $
If $\varepsilon$ is a solution of the equation
$\tau_{\mathcal{K}}\cdot T\varepsilon(X_{H\cdot\varepsilon})
= T\lambda \cdot \tilde{X}\cdot \varepsilon,$
if and only if it is a solution of the equation
$T\gamma \cdot
\tilde{X}^\varepsilon= X_{\mathcal{K}} \cdot \varepsilon $. Here $
X_{H \cdot \varepsilon}$ is the Hamiltonian vector field of the function
$H \cdot \varepsilon: T^* Q\rightarrow \mathbb{R}, $ and $X_{\mathcal{K}}$
is the nonholonomic dynamical vector field of the distributional RCH system.
The equation $T\gamma \cdot \tilde{X}^\varepsilon= X_{\mathcal{K}} \cdot \varepsilon,$
is called the Type II of Hamilton-Jacobi equation for the distributional RCH system
$(\mathcal{K},\omega_{\mathcal {K}},H_{\mathcal{K}}, F_{\mathcal{K}}, u_{\mathcal{K}})$.
\end{theo}

\noindent{\bf Proof: } In the same way, from Definition 2.3 we have that
the dynamical vector field of the distributional RCH system
$(\mathcal{K},\omega_{\mathcal {K}}, H_{\mathcal{K}}, F_{\mathcal{K}}, u_{\mathcal{K}})$
is the synthetic of the nonholonomic dynamical vector field $X_{\mathcal{K}}$ and
the vector fields $F_{\mathcal{K}}$ and $u_{\mathcal{K}}$, that is,
$\tilde{X}=X_{(\mathcal{K},\omega_{\mathcal{K}},
H_{\mathcal{K}}, F_{\mathcal{K}}, u_{\mathcal{K}})}
=X_\mathcal {K}+ F_{\mathcal{K}}+u_{\mathcal{K}}$,
and $F_{\mathcal{K}}=\tau_{\mathcal{K}}\cdot \textnormal{vlift}(F_{\mathcal{M}})X_H$,
and $u_{\mathcal{K}}=\tau_{\mathcal{K}}\cdot \textnormal{vlift}(u_{\mathcal{M}})X_H$.
Note that $T\pi_{Q}\cdot \textnormal{vlift}(F_{\mathcal{M}})X_H=T\pi_{Q}\cdot \textnormal{vlift}(u_{\mathcal{M}})X_H=0, $
then we have that $T\pi_{Q}\cdot F_{\mathcal{K}}=T\pi_{Q}\cdot u_{\mathcal{K}}=0,$
and hence $T\pi_{Q}\cdot \tilde{X}\cdot \varepsilon =T\pi_{Q}\cdot X_{\mathcal{K}}\cdot \varepsilon. $
On the other hand, we note that
$\textmd{Im}(\gamma)\subset \mathcal{M}, $ and
$\textmd{Im}(T\gamma)\subset \mathcal{K}, $ in this case,
$\omega_{\mathcal{K}}\cdot
\tau_{\mathcal{K}}=\tau_{\mathcal{K}}\cdot \omega_{\mathcal{M}}=
\tau_{\mathcal{K}}\cdot i_{\mathcal{M}}^* \cdot \omega, $ along
$\textmd{Im}(T\gamma)$. Moreover, from the distributional Hamiltonian equation (2.3),
we have that $X_{\mathcal{K}}= \tau_{\mathcal{K}}\cdot X_H.$
Note that $\varepsilon(\mathcal{M})\subset \mathcal{M},$ and
$T\pi_{Q}(X_H\cdot \varepsilon(q,p))\in
\mathcal{D}_{q}, \; \forall q \in Q, \; (q,p) \in \mathcal{M}(\subset T^* Q), $
and hence $X_H\cdot \varepsilon \in \mathcal{F}$ along $\varepsilon$,
and $\tau_{\mathcal{K}}\cdot X_{H}\cdot \varepsilon
= X_{\mathcal{K}}\cdot \varepsilon \in \mathcal{K}$.
Thus, using the non-degenerate
distributional two-form $\omega_{\mathcal{K}}$, from Lemma 3.3, if
we take that $v= X_{\mathcal{K}}\cdot
\varepsilon \in \mathcal{K} (\subset \mathcal{F}), $ and for any $w
\in \mathcal{F}, \; T\lambda(w)\neq 0, $ and
$\tau_{\mathcal{K}}\cdot w \neq 0, $ then we have that
\begin{align*}
& \omega_{\mathcal{K}}(T\gamma \cdot \tilde{X}^\varepsilon, \;
\tau_{\mathcal{K}}\cdot w)=
\omega_{\mathcal{K}}(\tau_{\mathcal{K}}\cdot T\gamma \cdot
X^\varepsilon, \; \tau_{\mathcal{K}}\cdot w)\\ & =
\tau_{\mathcal{K}}\cdot i_{\mathcal{M}}^* \cdot\omega(T\gamma \cdot
T\pi_{Q}\cdot \tilde{X}\cdot \varepsilon, \; w ) = \tau_{\mathcal{K}}\cdot
i_{\mathcal{M}}^* \cdot\omega (T(\gamma \cdot \pi_Q)\cdot X_{\mathcal{K}}\cdot \varepsilon, \; w)\\
& =\tau_{\mathcal{K}}\cdot i_{\mathcal{M}}^* \cdot(\omega(X_{\mathcal{K}}\cdot
\varepsilon, \; w-T(\gamma \cdot \pi_Q)\cdot w)
-\mathbf{d}\gamma(T\pi_{Q}(X_{\mathcal{K}}\cdot \varepsilon), \; T\pi_{Q}(w)))\\
& = \tau_{\mathcal{K}}\cdot i_{\mathcal{M}}^* \cdot\omega(X_{\mathcal{K}}\cdot
\varepsilon, \; w) - \tau_{\mathcal{K}}\cdot i_{\mathcal{M}}^* \cdot
\omega(X_{\mathcal{K}}\cdot \varepsilon, \; T\lambda \cdot w)\\
& \;\;\;\;\;\;
-\tau_{\mathcal{K}}\cdot i_{\mathcal{M}}^* \cdot\mathbf{d}\gamma
(T\pi_{Q}(X_{\mathcal{K}}\cdot \varepsilon), \; T\pi_{Q}(w))\\
& = \omega_{\mathcal{K}}( \tau_{\mathcal{K}}\cdot X_{\mathcal{K}} \cdot \varepsilon,
\; \tau_{\mathcal{K}}\cdot w) -
\omega_{\mathcal{K}}(\tau_{\mathcal{K}}\cdot X_{\mathcal{K}}\cdot \varepsilon, \;
\tau_{\mathcal{K}}\cdot T\lambda \cdot w)\\
& \;\;\;\;\;\;
+\tau_{\mathcal{K}}\cdot i_{\mathcal{M}}^* \cdot \lambda^* \omega(X_{\mathcal{K}}\cdot \varepsilon, \; w)\\
& = \omega_{\mathcal{K}}(X_{\mathcal{K}}\cdot \varepsilon, \;
\tau_{\mathcal{K}} \cdot w) -
\omega_{\mathcal{K}}(X_{\mathcal{K}}\cdot \varepsilon,
\; T\lambda \cdot w)+ \omega_{\mathcal{K}}(T\lambda\cdot X_{\mathcal{K}}\cdot \varepsilon,
\; T\lambda \cdot w),
\end{align*}
where we have used that
$ \tau_{\mathcal{K}}\cdot T\gamma= T\gamma, \; \tau_{\mathcal{K}}\cdot T\lambda= T\lambda, $ and
$\tau_{\mathcal{K}}\cdot X_{\mathcal{K}}\cdot \varepsilon = X_{\mathcal{K}}\cdot
\varepsilon, $ since $\textmd{Im}(T\gamma)\subset \mathcal{K}. $
From the distributional Hamiltonian equation
$\mathbf{i}_{X_\mathcal{K}}\omega_{\mathcal{K}}=\mathbf{d}H_\mathcal
{K}$, we have that $X_{\mathcal{K}}=\tau_{\mathcal{K}}\cdot X_H. $
On the other hand,
$\varepsilon: T^* Q \rightarrow T^* Q $ is symplectic, and $
X_H\cdot \varepsilon = T\varepsilon \cdot X_{H\cdot\varepsilon}, $ along
$\varepsilon$, and hence $ X_{\mathcal{K}}\cdot \varepsilon
=\tau_{\mathcal{K}}\cdot X_H \cdot \varepsilon=
\tau_{\mathcal{K}}\cdot T\varepsilon \cdot X_{H \cdot \varepsilon}, $ along $\varepsilon$.
Note that
$T\lambda \cdot X_\mathcal{K}\cdot \varepsilon=T\gamma \cdot
T\pi_Q\cdot X_\mathcal{K}\cdot \varepsilon=T\gamma \cdot
T\pi_Q\cdot \tilde{X}\cdot \varepsilon=T\lambda\cdot \tilde{X}\cdot \varepsilon.$
Then we have that
\begin{align*}
& \omega_{\mathcal{K}}(T\gamma \cdot \tilde{X}^\varepsilon, \;
\tau_{\mathcal{K}}\cdot w)-
\omega_{\mathcal{K}}(X_{\mathcal{K}}\cdot \varepsilon, \;
\tau_{\mathcal{K}} \cdot w) \nonumber \\
& = - \omega_{\mathcal{K}}(X_{\mathcal{K}}\cdot \varepsilon,
\; T\lambda \cdot w)+ \omega_{\mathcal{K}}(T\lambda\cdot X_{\mathcal{K}}\cdot \varepsilon,
\; T\lambda \cdot w) \\
& = -\omega_{\mathcal{K}}(\tau_{\mathcal{K}}\cdot X_H \cdot \varepsilon, \;
 T\lambda \cdot w)+ \omega_{\mathcal{K}}(T\lambda\cdot \tilde{X}\cdot \varepsilon,
\; T\lambda \cdot w)\\
&= \omega_{\mathcal{K}}(T\lambda\cdot \tilde{X}\cdot \varepsilon
-\tau_{\mathcal{K}}\cdot T\varepsilon \cdot X_{H \cdot \varepsilon},
\; T\lambda \cdot w).
\end{align*}
Because the distributional two-form
$\omega_{\mathcal{K}}$ is non-degenerate, it follows that the equation
$T\gamma\cdot \tilde{X}^\varepsilon= X_{\mathcal{K}}\cdot
\varepsilon ,$ is equivalent to the equation
$\tau_{\mathcal{K}}\cdot T\varepsilon \cdot X_{H\cdot\varepsilon}
= T\lambda\cdot \tilde{X}\cdot \varepsilon $.
Thus, $\varepsilon$ is a solution of the equation
$\tau_{\mathcal{K}}\cdot T\varepsilon\cdot X_{H\cdot\varepsilon}
= T\lambda \cdot \tilde{X}\cdot\varepsilon,$
if and only if it is a solution of
the Type II of Hamilton-Jacobi equation $T\gamma\cdot \tilde{X}^\varepsilon
= X_{\mathcal{K}}\cdot \varepsilon .$
\hskip 0.3cm $\blacksquare$

\begin{rema}
It is worthy of noting that, the Type I of Hamilton-Jacobi equation
$T\gamma\cdot \tilde{X}^\gamma= X_{\mathcal{K}}\cdot \gamma .$
is the equation of the differential one-form $\gamma$; and
the Type II of Hamilton-Jacobi equation $T\gamma\cdot \tilde{X}^\varepsilon
= X_{\mathcal{K}}\cdot \varepsilon .$ is the equation of
the symplectic diffeomorphism map $\varepsilon$.
If the nonholonomic RCH system we considered has not any constrains, in this case,
the distributional RCH system is just the RCH system itself.
From the above Type I and Type II of Hamilton-Jacobi theorems, that is,
Theorem 3.4 and Theorem 3.5, we can get the Theorem 2.6
and Theorem 2.7 in Wang \cite{wa13d}.
It shows that Theorem 3.4 and Theorem 3.5 can be regarded as an extension of two types of
Hamilton-Jacobi theorem for the RCH system given in \cite{wa13d} to that for the system with
nonholonomic context. In particular, in this case,
if both the external force and control of a nonholonomic RCH
system $(T^*Q,\omega,\mathcal{D},H, F,W)$ are also zeroes,
that is, $F=0 $ and $W=\emptyset$,
in this case the nonholonomic RCH system
is just a Hamiltonian system $(T^*Q,\omega,H)$
with the canonical symplectic form $\omega$ on $T^*Q$,
we can obtain two types of Hamilton-Jacobi
equation for the associated Hamiltonian system, which is given in Wang \cite{wa17}.
Thus, Theorem 3.4 and Theorem 3.5 can be regarded as an extension of two types of Hamilton-Jacobi
theorem for a Hamiltonian system to that for the system with external force,
control and nonholonomic constrain.
\end{rema}

\section{Hamilton-Jacobi Equations for a Nonholonomic Reduced Distributional RCH System }

It is well-known that the reduction of nonholonomically constrained mechanical systems
is a very important subject in geometric mechanics, and it is
regarded as a useful tool for simplifying and studying
concrete nonholonomic systems, see Koiller \cite{ko92},
Bates and $\acute{S}$niatycki \cite{basn93},  Cantrijn et al.
\cite{calemama99, calemama98}, Cendra et al. \cite{cemara01},
Cushman et al. \cite{cudusn10} and \cite{cukesnba95},
and de Le\'{o}n and Rodrigues \cite{lero89} and so on,
for more details and development.\\

In this section, we shall consider the nonholonomic reduction and
Hamilton-Jacobi theory of a nonholonomic RCH
system with symmetry. We first give the definition of
a nonholonomic RCH system with symmetry.
Then, by using the similar method in Le\'{o}n and Wang \cite{lewa15}
and Bates and $\acute{S}$niatycki \cite{basn93}.
and by analyzing carefully the structure of dynamical vector field
of the nonholonomic RCH system with symmetry,
we give a geometric formulation of
the nonholonomic reduced distributional RCH system,
Moreover, we derive precisely the geometric constraint conditions of
the non-degenerate, and nonholonomic reduced distributional two-form
for the nonholonomic reducible dynamical vector field,
that is, the two types of Hamilton-Jacobi equation for the
nonholonomic reduced distributional RCH system,
which are an extension of the above two types of Hamilton-Jacobi equation
for the distributional RCH system
given in section 3 under nonholonomic reduction.\\

Assume that the Lie group $G$ acts smoothly on the manifold $Q$ by the left,
and we also consider the natural lifted actions on $TQ$ and $T^* Q$,
and assume that the cotangent lifted left action $\Phi^{T^\ast}:
G\times T^\ast Q\rightarrow T^\ast Q$ is free, proper and
symplectic with respect to the canonical symplectic form
$\omega$ on $T^* Q$.
The orbit space $T^* Q/ G$ is a smooth manifold and the
canonical projection $\pi_{/G}: T^* Q \rightarrow T^* Q /G $ is
a surjective submersion. For the cotangent lifted left action
$\Phi^{T^\ast}: G\times T^\ast Q\rightarrow T^\ast Q$,
assume that $H: T^*Q \rightarrow \mathbb{R}$ is a
$G$-invariant Hamiltonian, and the fiber-preserving map
$F:T^\ast Q\rightarrow T^\ast Q$ and the control subset
$W$ of\; $T^\ast Q$ are both $G$-invariant,
and the $\mathcal{D}$-completely and
$\mathcal{D}$-regularly nonholonomic constraint $\mathcal{D}\subset
TQ$ is a $G$-invariant distribution for the tangent lifted left action $\Phi^{T}:
G\times TQ\rightarrow TQ$, that is, the tangent of group action maps
$\mathcal{D}_q$ to $\mathcal{D}_{gq}$ for any
$q\in Q $. A nonholonomic RCH system with symmetry
is 7-tuple $(T^*Q,G,\omega,\mathcal{D},H, F,W)$, which is an
RCH system with symmetry and $G$-invariant
nonholonomic constraint $\mathcal{D}$.\\

In the following we first consider the nonholonomic reduction
of a nonholonomic RCH system with symmetry
$(T^*Q,G,\omega,\mathcal{D},H, F,W)$.
Note that the Legendre transformation $\mathcal{F}L: TQ
\rightarrow T^*Q$ is a fiber-preserving map,
and $\mathcal{D}\subset TQ$ is $G$-invariant
for the tangent lifted left action $\Phi^{T}: G\times TQ\rightarrow TQ, $
then the constraint submanifold
$\mathcal{M}=\mathcal{F}L(\mathcal{D})\subset T^*Q$ is
$G$-invariant for the cotangent lifted left action $\Phi^{T^\ast}:
G\times T^\ast Q\rightarrow T^\ast Q$,
For the nonholonomic RCH system with symmetry
$(T^*Q,G, \omega,\mathcal{D},H,F,W)$,
in the same way, we define the distribution $\mathcal{F}$, which is the pre-image of the
nonholonomic constraints $\mathcal{D}$ for the map $T\pi_Q: TT^* Q
\rightarrow TQ$, that is, $\mathcal{F}=(T\pi_Q)^{-1}(\mathcal{D})$,
and the distribution $\mathcal{K}=\mathcal{F} \cap T\mathcal{M}$.
Moreover, we can also define the distributional two-form $\omega_\mathcal{K}$,
which is induced from the canonical symplectic form $\omega$ on $T^* Q$, that is,
$\omega_\mathcal{K}= \tau_{\mathcal{K}}\cdot \omega_{\mathcal{M}},$ and
$\omega_{\mathcal{M}}= i_{\mathcal{M}}^* \omega $.
If the admissibility condition $\mathrm{dim}\mathcal{M}=
\mathrm{rank}\mathcal{F}$ and the compatibility condition
$T\mathcal{M}\cap \mathcal{F}^\bot= \{0\}$ hold, then
$\omega_\mathcal{K}$ is non-degenerate as a
bilinear form on each fibre of $\mathcal{K}$, there exists a vector
field $X_\mathcal{K}$ on $\mathcal{M}$ which takes values in the
constraint distribution $\mathcal{K}$, such that for the function $H_\mathcal{K}$,
the following distributional Hamiltonian equation holds, that is,
\begin{align}
\mathbf{i}_{X_\mathcal{K}}\omega_\mathcal{K}
=\mathbf{d}H_\mathcal{K},
\label{4.1} \end{align}
where the function $H_{\mathcal{K}}$ satisfies
$\mathbf{d}H_{\mathcal{K}}= \tau_{\mathcal{K}}\cdot \mathbf{d}H_{\mathcal {M}}$,
and $H_\mathcal{M}= \tau_{\mathcal{M}}\cdot H$
is the restriction of $H$ to $\mathcal{M}$, and
from the equation (4.1), we have that
$X_{\mathcal{K}}=\tau_{\mathcal{K}}\cdot X_H $.\\

In the following we define that the quotient space
$\bar{\mathcal{M}}=\mathcal{M}/G$ of the $G$-orbit in $\mathcal{M}$
is a smooth manifold with projection $\pi_{/G}:
\mathcal{M}\rightarrow \bar{\mathcal{M}}( \subset T^* Q /G),$ which
is a surjective submersion. The reduced symplectic form
$\omega_{\bar{\mathcal{M}}}= \pi^*_{/G} \cdot \omega_{\mathcal{M}}$
on $\bar{\mathcal{M}}$ is induced from the symplectic form $\omega_{\mathcal{M}}
= i_{\mathcal{M}}^* \omega $ on $\mathcal{M}$.
Since $G$ is the symmetry group of the system
$(T^*Q,G,\omega,\mathcal{D},H, F,W)$, all intrinsically
defined vector fields and distributions are pushed down to
$\bar{\mathcal{M}}$. In particular, the vector field $X_\mathcal{M}$
on $\mathcal{M}$ is pushed down to a vector field
$X_{\bar{\mathcal{M}}}=T\pi_{/G}\cdot X_\mathcal{M}$, and the
distribution $\mathcal{K}$ is pushed down to a distribution
$T\pi_{/G}\cdot \mathcal{K}$ on $\bar{\mathcal{M}}$, and the
Hamiltonian $H$ is pushed down to $h_{\bar{\mathcal{M}}}$, such that
$h_{\bar{\mathcal{M}}}\cdot \pi_{/G}=
\tau_{\mathcal{M}}\cdot H$. However, $\omega_\mathcal{K}$ need not
to be pushed down to a distributional two-form defined on $T\pi_{/G}\cdot
\mathcal{K}$, despite of the fact that $\omega_\mathcal{K}$ is
$G$-invariant. This is because there may be infinitesimal symmetry
$\eta_{\mathcal{K}}$ that lies in $\mathcal{M}$, such that
$\mathbf{i}_{\eta_\mathcal{K}} \omega_\mathcal{K}\neq 0$. From Bates
and $\acute{S}$niatycki \cite{basn93}, we know that in order to eliminate
this difficulty, $\omega_\mathcal{K}$ is restricted to a
sub-distribution $\mathcal{U}$ of $\mathcal{K}$ defined by
$$\mathcal{U}=\{u\in\mathcal{K} \; | \; \omega_\mathcal{K}(u,v)
=0,\quad \forall \; v \in \mathcal{V}\cap \mathcal{K}\},$$ where
$\mathcal{V}$ is the distribution on $\mathcal{M}$ tangent to the
orbits of $G$ in $\mathcal{M}$ and it is spanned by the infinitesimal
symmetries. Clearly, $\mathcal{U}$ and $\mathcal{V}$ are both
$G$-invariant, project down to $\bar{\mathcal{M}}$ and
$T\pi_{/G}\cdot \mathcal{V}=0$, and define the distribution $\bar{\mathcal{K}}$ by
$\bar{\mathcal{K}}= T\pi_{/G}\cdot \mathcal{U}$. Moreover, we take
that $\omega_\mathcal{U}= \tau_{\mathcal{U}}\cdot
\omega_{\mathcal{M}}$ is the restriction of the induced symplectic form
$\omega_{\mathcal{M}}$ on $T^*\mathcal{M}$ fibrewise to the
distribution $\mathcal{U}$, where $\tau_{\mathcal{U}}$ is the
restriction map to distribution $\mathcal{U}$, and the
$\omega_{\mathcal{U}}$ is pushed down to a
distributional two-form $\omega_{\bar{\mathcal{K}}}$ on
$\bar{\mathcal{K}}$, such that $\pi_{/G}^*
\omega_{\bar{\mathcal{K}}}= \omega_{\mathcal{U}}$.
we know that distributional two-form
$\omega_{\bar{\mathcal{K}}}$ is not a "true two-form"
on a manifold, which is called the nonholonomic reduced
distributional two-form to avoid any confusion.\\

From the above construction we know that,
if the admissibility condition $\mathrm{dim}\bar{\mathcal{M}}=
\mathrm{rank}\bar{\mathcal{F}}$ and the compatibility condition
$T\bar{\mathcal{M}} \cap \bar{\mathcal{F}}^\bot= \{0\}$ hold, where
$\bar{\mathcal{F}}^\bot$ denotes the symplectic orthogonal of
$\bar{\mathcal{F}}$ with respect to the reduced symplectic form
$\omega_{\bar{\mathcal{M}}}$, then the nonholonomic reduced
distributional two-form
$\omega_{\bar{\mathcal{K}}}$ is non-degenerate as a bilinear form on
each fibre of $\bar{\mathcal{K}}$, and hence there exists a vector field
$X_{\bar{\mathcal{K}}}$ on $\bar{\mathcal{M}}$ which takes values in
the constraint distribution $\bar{\mathcal{K}}$, such that the
reduced distributional Hamiltonian equation holds, that is,
\begin{align}
\mathbf{i}_{X_{\bar{\mathcal{K}}}}\omega_{\bar{\mathcal{K}}}
=\mathbf{d}h_{\bar{\mathcal{K}}},
\label{4.2} \end{align}
where $\mathbf{d}h_{\bar{\mathcal{K}}}$ is the restriction of
$\mathbf{d}h_{\bar{\mathcal{M}}}$ to $\bar{\mathcal{K}}$ and
the function $h_{\bar{\mathcal{K}}}:\bar{M}(\subset T^* Q/G)\rightarrow \mathbb{R}$ satisfies
$\mathbf{d}h_{\bar{\mathcal{K}}}= \tau_{\bar{\mathcal{K}}}\cdot \mathbf{d}h_{\bar{\mathcal{M}}}$,
and $h_{\bar{\mathcal{M}}}\cdot \pi_{/G}= H_{\mathcal{M}}$ and
$H_{\mathcal{M}}$ is the restriction of the Hamiltonian function $H$
to $\mathcal{M}$, and the function
$h_{\bar{\mathcal{M}}}:\bar{M}(\subset T^* Q/G)\rightarrow \mathbb{R}$.
In addition, from the distributional Hamiltonian equation (4.1),
$\mathbf{i}_{X_\mathcal{K}}\omega_\mathcal{K}=\mathbf{d}H_\mathcal
{K},$ we have that $X_{\mathcal{K}}=\tau_{\mathcal{K}}\cdot X_H, $
and from the reduced distributional Hamiltonian equation (4.2),
$\mathbf{i}_{X_{\bar{\mathcal{K}}}}\omega_{\bar{\mathcal{K}}}
=\mathbf{d}h_{\bar{\mathcal{K}}}$, we have that
$X_{\bar{\mathcal{K}}}
=\tau_{\bar{\mathcal{K}}}\cdot X_{h_{\bar{\mathcal{K}}}},$
where $ X_{h_{\bar{\mathcal{K}}}}$ is the Hamiltonian vector field of
the function $h_{\bar{\mathcal{K}}}$ with respect to the reduced symplectic
form $\omega_{\bar{\mathcal{M}}}$,
and the vector fields $X_{\mathcal{K}}$
and $X_{\bar{\mathcal{K}}}$ are $\pi_{/G}$-related,
that is, $X_{\bar{\mathcal{K}}}\cdot \pi_{/G}=T\pi_{/G}\cdot X_{\mathcal{K}}.$ \\

Moreover, if considering the external force $F$ and control subset $W$,
and we define the vector fields $F_\mathcal{K}
=\tau_{\mathcal{K}}\cdot \textnormal{vlift}(F_{\mathcal{M}})X_H,$
and for a control law $u\in W$,
$u_\mathcal{K}= \tau_{\mathcal{K}}\cdot  \textnormal{vlift}(u_{\mathcal{M}})X_H,$
where $F_\mathcal{M}= \tau_{\mathcal{M}}\cdot F$ and
$u_\mathcal{M}= \tau_{\mathcal{M}}\cdot u$ are the restrictions of
$F$ and $u$ to $\mathcal{M}$, that is, $F_\mathcal{K}$ and $u_\mathcal{K}$
are the restrictions of the changes of Hamiltonian vector field $X_H$
under the actions of $F_\mathcal{M}$ and $u_\mathcal{M}$ to $\mathcal{K}$,
then the 5-tuple $(\mathcal{K},\omega_{\mathcal{K}},
H_\mathcal{K}, F_\mathcal{K}, u_\mathcal{K})$
is a distributional RCH system corresponding to the nonholonomic RCH system with symmetry
$(T^*Q,G,\omega,\mathcal{D},H,F,u)$,
and the dynamical vector field of the distributional RCH system
can be expressed by
\begin{align}
\tilde{X}=X_{(\mathcal{K},\omega_{\mathcal{K}},
H_{\mathcal{K}}, F_{\mathcal{K}}, u_{\mathcal{K}})}
=X_\mathcal {K}+ F_{\mathcal{K}}+u_{\mathcal{K}},
\label{4.3} \end{align}
which is the synthetic
of the nonholonomic dynamical vector field $X_{\mathcal{K}}$ and
the vector fields $F_{\mathcal{K}}$ and $u_{\mathcal{K}}$.
Assume that the vector fields $F_\mathcal{K}$ and $u_\mathcal{K}$
on $\mathcal{M}$ are pushed down to the vector fields
$f_{\bar{\mathcal{M}}}= T\pi_{/G}\cdot F_\mathcal{K}$ and
$u_{\bar{\mathcal{M}}}=T\pi_{/G}\cdot u_\mathcal{K}$ on $\bar{\mathcal{M}}$.
Then we define that $f_{\bar{\mathcal{K}}}=T\tau_{\bar{\mathcal{K}}}\cdot f_{\bar{\mathcal{M}}}$ and
$u_{\bar{\mathcal{K}}}=T\tau_{\bar{\mathcal{K}}}\cdot u_{\bar{\mathcal{M}}},$
that is, $f_{\bar{\mathcal{K}}}$ and
$u_{\bar{\mathcal{K}}}$ are the restrictions of
$f_{\bar{\mathcal{M}}}$ and $u_{\bar{\mathcal{M}}}$ to $\bar{\mathcal{K}}$,
where $\tau_{\bar{\mathcal{K}}}$
is the restriction map to distribution $\bar{\mathcal{K}}$,
and $T\tau_{\bar{\mathcal{K}}}$ is the tangent map of $\tau_{\bar{\mathcal{K}}}$.
Then the 5-tuple $(\bar{\mathcal{K}},\omega_{\bar{\mathcal{K}}},
h_{\bar{\mathcal{K}}}, f_{\bar{\mathcal{K}}}, u_{\bar{\mathcal{K}}})$
is a nonholonomic reduced distributional RCH system of the nonholonomic
reducible RCH system with symmetry $(T^*Q,G,\omega,\mathcal{D},H,F,W)$,
as well as with a control law $u \in W$.
Thus, the geometrical formulation of a nonholonomic reduced distributional
RCH system may be summarized as follows.

\begin{defi} (Nonholonomic Reduced Distributional RCH System)
Assume that the 7-tuple \\ $(T^*Q,G,\omega,\mathcal{D},H,F,W)$ is a nonholonomic
reducible RCH system with symmetry, where $\omega$ is the canonical
symplectic form on $T^* Q$, and $\mathcal{D}\subset TQ$ is a
$\mathcal{D}$-completely and $\mathcal{D}$-regularly nonholonomic
constraint of the system, and $\mathcal{D}$, $H, F$ and $W$ are all
$G$-invariant. If there exists a nonholonomic reduced distribution $\bar{\mathcal{K}}$,
an associated non-degenerate  and nonholonomic reduced
distributional two-form $\omega_{\bar{\mathcal{K}}}$
and a vector field $X_{\bar{\mathcal {K}}}$ on the reduced constraint
submanifold $\bar{\mathcal{M}}=\mathcal{M}/G, $ where
$\mathcal{M}=\mathcal{F}L(\mathcal{D})\subset T^*Q$, such that the
nonholonomic reduced distributional Hamiltonian equation
$ \mathbf{i}_{X_{\bar{\mathcal{K}}}}\omega_{\bar{\mathcal{K}}} =
\mathbf{d}h_{\bar{\mathcal{K}}}, $ holds,
where $\mathbf{d}h_{\bar{\mathcal{K}}}$ is the restriction of
$\mathbf{d}h_{\bar{\mathcal{M}}}$ to $\bar{\mathcal{K}}$ and
the function $h_{\bar{\mathcal{K}}}$ satisfies
$\mathbf{d}h_{\bar{\mathcal{K}}}= \tau_{\bar{\mathcal{K}}}\cdot \mathbf{d}h_{\bar{\mathcal{M}}}$
and $h_{\bar{\mathcal{M}}}\cdot \pi_{/G}= H_{\mathcal{M}}$,
and the vector fields $f_{\bar{\mathcal{K}}}=T\tau_{\bar{\mathcal{K}}}\cdot f_{\bar{\mathcal{M}}}$ and
$u_{\bar{\mathcal{K}}}=T\tau_{\bar{\mathcal{K}}}\cdot u_{\bar{\mathcal{M}}}$ as defined above.
Then the 5-tuple $(\bar{\mathcal{K}},\omega_{\bar{\mathcal {K}}},h_{\bar{\mathcal{K}}},
f_{\bar{\mathcal{K}}}, u_{\bar{\mathcal{K}}})$
is called a nonholonomic reduced distributional RCH system
of the nonholonomic reducible RCH system $(T^*Q,G,\omega,\mathcal{D},H,F,W)$
with a control law $u \in W$, and $X_{\bar{\mathcal {K}}}$ is
called a nonholonomic reduced dynamical vector field.
Denote by
\begin{align}
\hat{X}=X_{(\bar{\mathcal{K}},\omega_{\bar{\mathcal{K}}},
h_{\bar{\mathcal{K}}}, f_{\bar{\mathcal{K}}}, u_{\bar{\mathcal{K}}})}
=X_{\bar{\mathcal{K}}}+ f_{\bar{\mathcal{K}}}+u_{\bar{\mathcal{K}}}
\label{4.4} \end{align}
is the dynamical vector field of the
nonholonomic reduced distributional RCH system
$(\bar{\mathcal{K}},\omega_{\bar{\mathcal{K}}},h_{\bar{\mathcal{K}}}, f_{\bar{\mathcal{K}}}, u_{\bar{\mathcal{K}}})$, which is the synthetic
of the nonholonomic reduced dynamical vector field $X_{\bar{\mathcal{K}}}$ and
the vector fields $F_{\bar{\mathcal{K}}}$ and $u_{\bar{\mathcal{K}}}$.
Under the above
circumstances, we refer to $(T^*Q,G,\omega,\mathcal{D},H,F,u)$ as a
nonholonomic reducible RCH system with the associated
distributional RCH system
$(\mathcal{K},\omega_{\mathcal {K}},H_{\mathcal{K}}, F_{\mathcal{K}}, u_{\mathcal{K}})$
and the nonholonomic reduced distributional RCH system
$(\bar{\mathcal{K}},\omega_{\bar{\mathcal{K}}},h_{\bar{\mathcal{K}}},
f_{\bar{\mathcal{K}}}, u_{\bar{\mathcal{K}}})$.
The dynamical vector fields
$\tilde{X}=X_{(\mathcal{K},\omega_{\mathcal{K}},
H_{\mathcal{K}}, F_{\mathcal{K}}, u_{\mathcal{K}})}$
and $\hat{X}= X_{(\bar{\mathcal{K}},\omega_{\bar{\mathcal{K}}},
h_{\bar{\mathcal{K}}}, f_{\bar{\mathcal{K}}}, u_{\bar{\mathcal{K}}})}$
are $\pi_{/G}$-related, that is,
$\hat{X}\cdot \pi_{/G}=T\pi_{/G}\cdot \tilde{X}.$
\end{defi}

Since the non-degenerate and nonholonomic reduced distributional two-form
$\omega_{\bar{\mathcal{K}}}$ is not a "true two-form"
on a manifold, and it is not symplectic, and hence
the nonholonomic reduced distributional RCH system
$(\bar{\mathcal{K}},\omega_{\bar{\mathcal{K}}},h_{\bar{\mathcal{K}}}, f_{\bar{\mathcal{K}}}, u_{\bar{\mathcal{K}}})$ is not a Hamiltonian system,
and has no yet generating function,
and hence we can not describe the Hamilton-Jacobi equation for the nonholonomic reduced
distributional RCH system just like as in Theorem 1.1.
But, for a given nonholonomic reducible RCH system
$(T^*Q,G,\omega,\mathcal{D},H,F,u)$ with the associated
distributional RCH system
$(\mathcal{K},\omega_{\mathcal {K}},H_{\mathcal{K}}, F_{\mathcal{K}}, u_{\mathcal{K}})$
and the nonholonomic reduced distributional RCH system
$(\bar{\mathcal{K}},\omega_{\bar{\mathcal {K}}},h_{\bar{\mathcal{K}}}, f_{\bar{\mathcal{K}}}, u_{\bar{\mathcal{K}}})$, by using Lemma 3.3,
we can derive precisely
the geometric constraint conditions of the nonholonomic reduced distributional two-form
$\omega_{\bar{\mathcal{K}}}$ for the nonholonomic reducible dynamical vector field
$\tilde{X}=X_{(\mathcal{K},\omega_{\mathcal{K}},
H_{\mathcal{K}}, F_{\mathcal{K}}, u_{\mathcal{K}})}$,
that is, the two types of Hamilton-Jacobi equation for the
nonholonomic reduced distributional RCH system
$(\bar{\mathcal{K}},\omega_{\bar{\mathcal {K}}},h_{\bar{\mathcal{K}}}, f_{\bar{\mathcal{K}}}, u_{\bar{\mathcal{K}}})$.
At first, using the fact that the one-form $\gamma: Q
\rightarrow T^*Q $ is closed on $\mathcal{D}$ with respect to
$T\pi_Q: TT^* Q \rightarrow TQ, $
$\textmd{Im}(\gamma)\subset \mathcal{M}, $ and it is $G$-invariant,
as well as $ \textmd{Im}(T\gamma)\subset \mathcal{K}, $
we can prove the Type I of
Hamilton-Jacobi theorem for the nonholonomic reduced distributional
RCH system. For convenience, the maps involved in the
following theorem and its proof are shown in Diagram-3.
\begin{center}
\hskip 0cm \xymatrix{ & \mathcal{M} \ar[d]_{X_{\mathcal{K}}}
\ar[r]^{i_{\mathcal{M}}} & T^* Q \ar[d]_{X_{H}}
 \ar[r]^{\pi_Q}
  & Q \ar[d]_{\tilde{X}^\gamma} \ar[r]^{\gamma}
  & T^* Q \ar[d]_{\tilde{X}} \ar[r]^{\pi_{/G}} & T^* Q/G \ar[d]_{X_{h}}
  & \mathcal{\bar{M}} \ar[l]_{i_{\mathcal{\bar{M}}}} \ar[d]_{X_{\mathcal{\bar{K}}}}\\
  & \mathcal{K}
  & T(T^*Q) \ar[l]_{\tau_{\mathcal{K}}}
  & TQ \ar[l]_{T\gamma}
  & T(T^* Q) \ar[l]_{T\pi_Q} \ar[r]^{T\pi_{/G}}
  & T(T^* Q/G) \ar[r]^{\tau_{\mathcal{\bar{K}}}} & \mathcal{\bar{K}} }
\end{center}
$$\mbox{Diagram-3}$$

\begin{theo} (Type I of Hamilton-Jacobi Theorem for a Nonholonomic
Reduced Distributional RCH System)
For a given nonholonomic reducible RCH system
$(T^*Q,G,\omega,\mathcal{D},H,F,u)$ with the associated
distributional RCH system
$(\mathcal{K},\omega_{\mathcal {K}},H_{\mathcal{K}}, F_{\mathcal{K}}, u_{\mathcal{K}})$
and the nonholonomic reduced distributional RCH system
$(\bar{\mathcal{K}},\omega_{\bar{\mathcal{K}}},h_{\bar{\mathcal{K}}}, f_{\bar{\mathcal{K}}}, u_{\bar{\mathcal{K}}})$, assume that
$\gamma: Q \rightarrow T^*Q$ is an one-form on $Q$, and
$\tilde{X}^\gamma = T\pi_{Q}\cdot \tilde{X}\cdot \gamma$,
where $\tilde{X}=X_{(\mathcal{K},\omega_{\mathcal{K}},
H_{\mathcal{K}}, F_{\mathcal{K}}, u_{\mathcal{K}})}
=X_\mathcal {K}+ F_{\mathcal{K}}+u_{\mathcal{K}}$
is the dynamical vector field of the distributional RCH system
corresponding to the nonholonomic reducible RCH system with symmetry
$(T^*Q,G,\omega,\mathcal{D},H,F,u)$. Moreover,
assume that $\textmd{Im}(\gamma)\subset \mathcal{M}, $ and it is
$G$-invariant, $ \textmd{Im}(T\gamma)\subset \mathcal{K}, $ and
$\bar{\gamma}=\pi_{/G}(\gamma): Q \rightarrow T^* Q/G .$ If the
one-form $\gamma: Q \rightarrow T^*Q $ is closed on $\mathcal{D}$ with respect to
$T\pi_Q: TT^* Q \rightarrow TQ, $ then $\bar{\gamma}$ is a solution
of the equation $T\bar{\gamma}\cdot \tilde{X}^ \gamma =
X_{\bar{\mathcal{K}}}\cdot \bar{\gamma}. $ Here
$X_{\bar{\mathcal{K}}}$ is the nonholonomic reduced dynamical vector
field. The equation $T\bar{\gamma}\cdot \tilde{X}^ \gamma = X_{\bar{\mathcal{K}}}\cdot
\bar{\gamma},$ is called the Type I of Hamilton-Jacobi equation for
the nonholonomic reduced distributional RCH system
$(\bar{\mathcal{K}},\omega_{\bar{\mathcal{K}}},h_{\bar{\mathcal{K}}}, f_{\bar{\mathcal{K}}}, u_{\bar{\mathcal{K}}})$.
\end{theo}

\noindent{\bf Proof: } At first, for the dynamical vector field of the distributional RCH system
$(\mathcal{K},\omega_{\mathcal {K}}, H_{\mathcal{K}}, F_{\mathcal{K}}, u_{\mathcal{K}})$,
$\tilde{X}=X_{(\mathcal{K},\omega_{\mathcal{K}},
H_{\mathcal{K}}, F_{\mathcal{K}}, u_{\mathcal{K}})}
=X_\mathcal {K}+ F_{\mathcal{K}}+u_{\mathcal{K}}$,
and $F_{\mathcal{K}}=\tau_{\mathcal{K}}\cdot \textnormal{vlift}(F_{\mathcal{M}})X_H$,
and $u_{\mathcal{K}}=\tau_{\mathcal{K}}\cdot \textnormal{vlift}(u_{\mathcal{M}})X_H$,
note that $T\pi_{Q}\cdot \textnormal{vlift}(F_{\mathcal{M}})X_H=T\pi_{Q}\cdot \textnormal{vlift}(u_{\mathcal{M}})X_H=0, $
then we have that $T\pi_{Q}\cdot F_{\mathcal{K}}=T\pi_{Q}\cdot u_{\mathcal{K}}=0,$
and hence $T\pi_{Q}\cdot \tilde{X}\cdot \gamma=T\pi_{Q}\cdot X_{\mathcal{K}}\cdot \gamma. $
Moreover, from Theorem 3.4, we know that
$\gamma$ is a solution of the Type I of Hamilton-Jacobi equation
$T\gamma\cdot \tilde{X}^\gamma= X_{\mathcal{K}}\cdot \gamma .$ Next, we note that
$\textmd{Im}(\gamma)\subset \mathcal{M}, $ and it is $G$-invariant,
$ \textmd{Im}(T\gamma)\subset \mathcal{K}, $ and hence
$\textmd{Im}(T\bar{\gamma})\subset \bar{\mathcal{K}}, $ in this case,
$\pi^*_{/G}\cdot\omega_{\bar{\mathcal{K}}}\cdot\tau_{\bar{\mathcal{K}}}= \tau_{\mathcal{U}}\cdot
\omega_{\mathcal{M}}= \tau_{\mathcal{U}}\cdot i_{\mathcal{M}}^*\cdot
\omega, $ along $\textmd{Im}(T\bar{\gamma})$.
From the distributional Hamiltonian equation (4.1),
we have that $X_{\mathcal{K}}= \tau_{\mathcal{K}}\cdot X_H,$
and $\tau_{\mathcal{K}}\cdot X_{H}\cdot \gamma
= X_{\mathcal{K}}\cdot \gamma \in \mathcal{K}$.
Because the vector fields $X_{\mathcal{K}}$
and $X_{\bar{\mathcal{K}}}$ are $\pi_{/G}$-related,
$T\pi_{/G}(X_{\mathcal{K}})=X_{\bar{\mathcal{K}}}\cdot \pi_{/G}$,
and hence $\tau_{\bar{\mathcal{K}}}\cdot T\pi_{/G}(X_{\mathcal{K}}\cdot \gamma)
=\tau_{\bar{\mathcal{K}}}\cdot (T\pi_{/G}(X_{\mathcal{K}}))\cdot (\gamma)
= \tau_{\bar{\mathcal{K}}}\cdot (X_{\bar{\mathcal{K}}}\cdot \pi_{/G})\cdot (\gamma)
= \tau_{\bar{\mathcal{K}}}\cdot X_{\bar{\mathcal{K}}}\cdot \pi_{/G}(\gamma)
= X_{\bar{\mathcal{K}}}\cdot \bar{\gamma}.$
Thus, using the non-degenerate nonholonomic reduced distributional two-form
$\omega_{\bar{\mathcal{K}}}$, from Lemma 3.3(ii) and (iii), if we take that
$v=X_{\mathcal{K}}\cdot \gamma \in \mathcal{K} (\subset \mathcal{F}),$
and for any $w \in \mathcal{F}, \; T\lambda(w)\neq 0, $ and
$\tau_{\bar{\mathcal{K}}}\cdot T\pi_{/G}\cdot w \neq 0, $ then we have that
\begin{align*}
& \omega_{\bar{\mathcal{K}}}(T\bar{\gamma} \cdot \tilde{X}^\gamma, \;
\tau_{\bar{\mathcal{K}}}\cdot T\pi_{/G} \cdot w)
= \omega_{\bar{\mathcal{K}}}(\tau_{\bar{\mathcal{K}}}\cdot T(\pi_{/G} \cdot
\gamma) \cdot \tilde{X}^\gamma, \; \tau_{\bar{\mathcal{K}}}\cdot T\pi_{/G} \cdot w )\\
& = \pi^*_{/G}
\cdot \omega_{\bar{\mathcal{K}}}\cdot\tau_{\bar{\mathcal{K}}}
(T\gamma \cdot T\pi_{Q}\cdot \tilde{X}\cdot \gamma, \; w) =
\tau_{\mathcal{U}}\cdot i^*_{\mathcal{M}} \cdot\omega(T\gamma \cdot
T\pi_{Q}\cdot X_{\mathcal{K}}\cdot \gamma, \; w)\\
& = \tau_{\mathcal{U}}\cdot i^*_{\mathcal{M}} \cdot
\omega(T(\gamma \cdot \pi_Q)\cdot X_{\mathcal{K}}\cdot \gamma, \; w) \\
& = \tau_{\mathcal{U}}\cdot i^*_{\mathcal{M}} \cdot
(\omega(X_{\mathcal{K}}\cdot \gamma, \; w-T(\gamma \cdot \pi_Q)\cdot w)- \mathbf{d}\gamma(T\pi_{Q}
(X_{\mathcal{K}}\cdot \gamma), \; T\pi_{Q}(w)))\\
& = \tau_{\mathcal{U}}\cdot i^*_{\mathcal{M}} \cdot \omega(X_{\mathcal{K}} \cdot
\gamma, \; w) - \tau_{\mathcal{U}}\cdot i^*_{\mathcal{M}} \cdot
\omega(X_{\mathcal{K}}\cdot \gamma, \; T(\gamma \cdot \pi_Q) \cdot w)\\
& \;\;\;\;\;\;
- \tau_{\mathcal{U}}\cdot i^*_{\mathcal{M}} \cdot\mathbf{d}\gamma(T\pi_{Q}
(X_{\mathcal{K}}\cdot \gamma), \; T\pi_{Q}(w))\\
& =\pi^*_{/G}\cdot \omega_{\bar{\mathcal{K}}}\cdot\tau_{\bar{\mathcal{K}}}(X_{\mathcal{K}}\cdot \gamma, \;
w) - \pi^*_{/G}\cdot \omega_{\bar{\mathcal{K}}}\cdot\tau_{\bar{\mathcal{K}}}(X_{\mathcal{K}}\cdot \gamma,
\; T(\gamma \cdot \pi_Q) \cdot w)\\
& \;\;\;\;\;\; - \tau_{\mathcal{U}}\cdot i^*_{\mathcal{M}}
\cdot\mathbf{d}\gamma(T\pi_{Q}(X_{\mathcal{K}}\cdot \gamma), \; T\pi_{Q}(w))\\
& = \omega_{\bar{\mathcal{K}}}(\tau_{\bar{\mathcal{K}}}\cdot T\pi_{/G}(X_{\mathcal{K}}\cdot \gamma), \;
\tau_{\bar{\mathcal{K}}}\cdot T\pi_{/G} \cdot w) - \omega_{\bar{\mathcal{K}}}(\tau_{\bar{\mathcal{K}}}\cdot T\pi_{/G}(X_{\mathcal{K}}\cdot \gamma), \;
\tau_{\bar{\mathcal{K}}}\cdot T(\pi_{/G} \cdot\gamma) \cdot T\pi_{Q}(w))\\
& \;\;\;\;\;\; -\tau_{\mathcal{U}}\cdot i^*_{\mathcal{M}}
\cdot\mathbf{d}\gamma(T\pi_{Q}(X_{\mathcal{K}}\cdot \gamma), \; T\pi_{Q}(w))\\
& = \omega_{\bar{\mathcal{K}}}(X_{\bar{\mathcal{K}}} \cdot
\bar{\gamma}, \; \tau_{\bar{\mathcal{K}}}\cdot T\pi_{/G} \cdot w)-
\omega_{\bar{\mathcal{K}}}(X_{\bar{\mathcal{K}}} \cdot
\bar{\gamma}, \; T\bar{\gamma} \cdot T\pi_{Q}(w))- \tau_{\mathcal{U}}\cdot
i^*_{\mathcal{M}} \cdot\mathbf{d}\gamma(T\pi_{Q}(X_{\mathcal{K}}\cdot \gamma),
\; T\pi_{Q}(w)),
\end{align*}
where we have used that $\tau_{\bar{\mathcal{K}}}\cdot T\pi_{/G}(X_{\mathcal{K}}\cdot \gamma)
=X_{\bar{\mathcal{K}}}\cdot \bar{\gamma}, $ and
$\tau_{\bar{\mathcal{K}}}\cdot T\bar{\gamma}=T\bar{\gamma}, $ since
$\textmd{Im}(T\bar{\gamma})\subset \bar{\mathcal{K}}. $
If the one-form $\gamma: Q \rightarrow T^*Q $ is closed on $\mathcal{D}$ with respect to
$T\pi_Q: TT^* Q \rightarrow TQ, $ then we have that
$\mathbf{d}\gamma(T\pi_{Q}(X_{\mathcal{K}}\cdot \gamma), \; T\pi_{Q}(w))=0, $
since $X_{\mathcal{K}}\cdot \gamma, \; w \in \mathcal{F},$ and
$T\pi_{Q}(X_{\mathcal{K}}\cdot \gamma), \; T\pi_{Q}(w) \in \mathcal{D}, $ and hence
$$
\tau_{\mathcal{U}}\cdot
i_{\mathcal{M}}^* \cdot\mathbf{d}\gamma(T\pi_{Q}(X_{\mathcal{K}}\cdot \gamma),
\; T\pi_{Q}(w))=0,
$$
and
\begin{equation}
\omega_{\bar{\mathcal{K}}}(T\bar{\gamma} \cdot \tilde{X}^\gamma, \;
\tau_{\bar{\mathcal{K}}}\cdot T\pi_{/G} \cdot w)- \omega_{\bar{\mathcal{K}}}(X_{\bar{\mathcal{K}}} \cdot
\bar{\gamma}, \; \tau_{\bar{\mathcal{K}}}\cdot T\pi_{/G} \cdot w)
= -\omega_{\bar{\mathcal{K}}}(X_{\bar{\mathcal{K}}} \cdot
\bar{\gamma}, \; T\bar{\gamma} \cdot T\pi_{Q}(w)).
\label{4.5} \end{equation}
If $\bar{\gamma}$ satisfies the equation $
T\bar{\gamma}\cdot \tilde{X}^ \gamma
= X_{\bar{\mathcal{K}}}\cdot \bar{\gamma} ,$
from Lemma 3.3(i) we know that the right side of (4.5) becomes that
\begin{align*}
 -\omega_{\bar{\mathcal{K}}}(X_{\bar{\mathcal{K}}} \cdot
\bar{\gamma}, \; T\bar{\gamma} \cdot T\pi_{Q}(w))
& = -\omega_{\bar{\mathcal{K}}}\cdot\tau_{\bar{\mathcal{K}}}(T\bar{\gamma}\cdot \tilde{X}^\gamma, \; T\bar{\gamma} \cdot T\pi_{Q}(w))\\
& = -\bar{\gamma}^*\omega_{\bar{\mathcal{K}}}\cdot\tau_{\bar{\mathcal{K}}}
(T\pi_{Q} \cdot \tilde{X}\cdot \gamma, \; T\pi_{Q}(w))\\
& = - \gamma^* \cdot \pi^*_{/G}\cdot \omega_{\bar{\mathcal{K}}}\cdot\tau_{\bar{\mathcal{K}}}
(T\pi_{Q} \cdot X_{\mathcal{K}} \cdot \gamma, \; T\pi_{Q}(w))\\
& = - \gamma^* \cdot \tau_{\mathcal{U}}\cdot
i_{\mathcal{M}}^* \cdot \omega(T\pi_{Q}(X_{\mathcal{K}}\cdot\gamma), \; T\pi_{Q}(w))\\
& = -\tau_{\mathcal{U}}\cdot
i_{\mathcal{M}}^* \cdot\gamma^*\omega( T\pi_{Q}(X_{\mathcal{K}}\cdot\gamma), \; T\pi_{Q}(w))\\
& = \tau_{\mathcal{U}}\cdot i_{\mathcal{M}}^* \cdot
\mathbf{d}\gamma(T\pi_{Q}( X_{\mathcal{K}}\cdot\gamma ), \; T\pi_{Q}(w))=0,
\end{align*}
where $\gamma^*\cdot \tau_{\mathcal{U}}\cdot i^*_{\mathcal{M}}
\cdot \omega= \tau_{\mathcal{U}}\cdot i^*_{\mathcal{M}}
\cdot\gamma^*\cdot \omega, $ because $\textmd{Im}(\gamma)\subset
\mathcal{M}. $
But, since the nonholonomic reduced distributional two-form
$\omega_{\bar{\mathcal{K}}}$ is non-degenerate,
the left side of (4.5) equals zero, only when
$\bar{\gamma}$ satisfies the equation $
T\bar{\gamma}\cdot \tilde{X}^ \gamma = X_{\bar{\mathcal{K}}}\cdot
\bar{\gamma} .$ Thus,
if the one-form $\gamma: Q \rightarrow T^*Q $ is closed on $\mathcal{D}$ with respect to
$T\pi_Q: TT^* Q \rightarrow TQ, $ then $\bar{\gamma}$ must be a solution of
the Type I of Hamilton-Jacobi equation
$T\bar{\gamma}\cdot \tilde{X}^ \gamma = X_{\bar{\mathcal{K}}}\cdot
\bar{\gamma}. $
\hskip 0.3cm $\blacksquare$\\

Next, for any $G$-invariant symplectic map $\varepsilon: T^* Q \rightarrow T^* Q $, we can prove
the following Type II of
Hamilton-Jacobi theorem for the nonholonomic reduced distributional RCH system.
For convenience, the maps involved in the following theorem and its
proof are shown in Diagram-4.
\begin{center}
\hskip 0cm \xymatrix{ & \mathcal{M} \ar[d]_{X_{\mathcal{K}}}
\ar[r]^{i_{\mathcal{M}}} & T^* Q \ar[d]_{X_{H\cdot \varepsilon}}
\ar[dr]^{\tilde{X}^\varepsilon} \ar[r]^{\pi_Q}
  & Q \ar[r]^{\gamma} & T^* Q \ar[d]_{\tilde{X}} \ar[r]^{\pi_{/G}} & T^* Q/G \ar[d]_{X_{h}}
  & \mathcal{\bar{M}} \ar[l]_{i_{\mathcal{\bar{M}}}} \ar[d]_{X_{\mathcal{\bar{K}}}}\\
  & \mathcal{K}
  & T(T^*Q) \ar[l]_{\tau_{\mathcal{K}}}
  & TQ \ar[l]_{T\gamma}
  & T(T^* Q) \ar[l]_{T\pi_Q} \ar[r]^{T\pi_{/G}}
  & T(T^* Q/G) \ar[r]^{\tau_{\mathcal{\bar{K}}}} & \mathcal{\bar{K}} }
\end{center}
$$\mbox{Diagram-4}$$

\begin{theo} (Type II of Hamilton-Jacobi Theorem for a Nonholonomic
Reduced Distributional RCH System)
For a given nonholonomic reducible RCH system
$(T^*Q,G,\omega,\mathcal{D},H,F,u)$ with the associated
distributional RCH system
$(\mathcal{K},\omega_{\mathcal {K}},H_{\mathcal{K}}, F_{\mathcal{K}}, u_{\mathcal{K}})$
and the nonholonomic reduced distributional RCH system
$(\bar{\mathcal{K}},\omega_{\bar{\mathcal{K}}},h_{\bar{\mathcal{K}}},
f_{\bar{\mathcal{K}}}, u_{\bar{\mathcal{K}}})$, assume that
$\gamma: Q \rightarrow T^*Q$ is an one-form on $Q$, and $\lambda=
\gamma \cdot \pi_{Q}: T^* Q \rightarrow T^* Q, $ and for any $G$-invariant
symplectic map $\varepsilon: T^* Q \rightarrow T^* Q $, denote by
$\tilde{X}^\varepsilon = T\pi_{Q}\cdot \tilde{X}\cdot \varepsilon$,
where $\tilde{X}=X_{(\mathcal{K},\omega_{\mathcal{K}},
H_{\mathcal{K}}, F_{\mathcal{K}}, u_{\mathcal{K}})}
=X_\mathcal {K}+ F_{\mathcal{K}}+u_{\mathcal{K}}$
is the dynamical vector field of the distributional RCH system
corresponding to the nonholonomic reducible RCH system with symmetry
$(T^*Q,G,\omega,\mathcal{D},H,F,u)$. Moreover,
assume that $\textmd{Im}(\gamma)\subset \mathcal{M}, $ and it is
$G$-invariant, $\varepsilon(\mathcal{M})\subset \mathcal{M}$,
$ \textmd{Im}(T\gamma)\subset \mathcal{K}, $ and
$\bar{\gamma}=\pi_{/G}(\gamma): Q \rightarrow T^* Q/G $, and
$\bar{\lambda}=\pi_{/G}(\lambda): T^* Q \rightarrow T^* Q/G, $ and
$\bar{\varepsilon}=\pi_{/G}(\varepsilon): T^* Q \rightarrow T^* Q/G. $ Then
$\varepsilon$ and $\bar{\varepsilon}$ satisfy the equation
$\tau_{\bar{\mathcal{K}}}\cdot T\bar{\varepsilon}\cdot X_{h_{\bar{\mathcal{K}}}\cdot
\bar{\varepsilon}}= T\bar{\lambda} \cdot \tilde{X}\cdot \varepsilon, $ if and only if they satisfy the
equation $T\bar{\gamma}\cdot \tilde{X}^ \varepsilon =
X_{\bar{\mathcal{K}}}\cdot \bar{\varepsilon}. $ Here
$ X_{h_{\bar{\mathcal{K}}} \cdot\bar{\varepsilon}}$ is the Hamiltonian
vector field of the function $h_{\bar{\mathcal{K}}}\cdot \bar{\varepsilon}: T^* Q\rightarrow
\mathbb{R}, $ and $X_{\bar{\mathcal{K}}}$ is the nonholonomic
reduced dynamical vector field. The equation $
T\bar{\gamma}\cdot \tilde{X}^\varepsilon = X_{\bar{\mathcal{K}}}\cdot
\bar{\varepsilon},$ is called the Type II of Hamilton-Jacobi equation for the
nonholonomic reduced distributional RCH system
$(\bar{\mathcal{K}},\omega_{\bar{\mathcal{K}}},h_{\bar{\mathcal{K}}}, f_{\bar{\mathcal{K}}}, u_{\bar{\mathcal{K}}})$.
\end{theo}

\noindent{\bf Proof: } In the same way, for the dynamical vector field of the distributional RCH system
$(\mathcal{K},\omega_{\mathcal {K}}, H_{\mathcal{K}}, F_{\mathcal{K}}, u_{\mathcal{K}})$,
$\tilde{X}=X_{(\mathcal{K},\omega_{\mathcal{K}},
H_{\mathcal{K}}, F_{\mathcal{K}}, u_{\mathcal{K}})}
=X_\mathcal {K}+ F_{\mathcal{K}}+u_{\mathcal{K}}$,
and $F_{\mathcal{K}}=\tau_{\mathcal{K}}\cdot \textnormal{vlift}(F_{\mathcal{M}})X_H$,
and $u_{\mathcal{K}}=\tau_{\mathcal{K}}\cdot \textnormal{vlift}(u_{\mathcal{M}})X_H$,
note that $T\pi_{Q}\cdot \textnormal{vlift}(F_{\mathcal{M}})X_H=T\pi_{Q}\cdot \textnormal{vlift}(u_{\mathcal{M}})X_H=0, $
then we have that $T\pi_{Q}\cdot F_{\mathcal{K}}=T\pi_{Q}\cdot u_{\mathcal{K}}=0,$
and hence $T\pi_{Q}\cdot \tilde{X}\cdot \varepsilon
=T\pi_{Q}\cdot X_{\mathcal{K}}\cdot \varepsilon. $
Next, we note that
$\textmd{Im}(\gamma)\subset \mathcal{M}, $ and it is $G$-invariant,
$ \textmd{Im}(T\gamma)\subset \mathcal{K}, $ and hence
$\textmd{Im}(T\bar{\gamma})\subset \bar{\mathcal{K}}, $ in this case,
$\pi^*_{/G}\cdot\omega_{\bar{\mathcal{K}}}\cdot\tau_{\bar{\mathcal{K}}}= \tau_{\mathcal{U}}\cdot
\omega_{\mathcal{M}}= \tau_{\mathcal{U}}\cdot i_{\mathcal{M}}^*\cdot
\omega, $ along $\textmd{Im}(T\bar{\gamma})$.
Moreover, from the distributional Hamiltonian equation (4.1),
we have that $X_{\mathcal{K}}= \tau_{\mathcal{K}}\cdot X_H.$
Note that $\varepsilon(\mathcal{M})\subset \mathcal{M},$ and
$T\pi_{Q}(X_H\cdot \varepsilon(q,p))\in
\mathcal{D}_{q}, \; \forall q \in Q, \; (q,p) \in \mathcal{M}(\subset T^* Q), $
and hence $X_H\cdot \varepsilon \in \mathcal{F}$ along $\varepsilon$,
and $\tau_{\mathcal{K}}\cdot X_{H}\cdot \varepsilon
= X_{\mathcal{K}}\cdot \varepsilon \in \mathcal{K}$.
Because the vector fields $X_{\mathcal{K}}$
and $X_{\bar{\mathcal{K}}}$ are $\pi_{/G}$-related, then
$T\pi_{/G}(X_{\mathcal{K}})=X_{\bar{\mathcal{K}}}\cdot \pi_{/G}$,
and hence $\tau_{\bar{\mathcal{K}}}\cdot T\pi_{/G}(X_{\mathcal{K}}\cdot \varepsilon)
=\tau_{\bar{\mathcal{K}}}\cdot (T\pi_{/G}(X_{\mathcal{K}}))\cdot (\varepsilon)
= \tau_{\bar{\mathcal{K}}}\cdot (X_{\bar{\mathcal{K}}}\cdot \pi_{/G})\cdot (\varepsilon)
= \tau_{\bar{\mathcal{K}}}\cdot X_{\bar{\mathcal{K}}}\cdot \pi_{/G}(\varepsilon)
= X_{\bar{\mathcal{K}}}\cdot \bar{\varepsilon}.$
Thus, using the non-degenerate and nonholonomic reduced distributional two-form
$\omega_{\bar{\mathcal{K}}}$, from Lemma 3.3, if we take that
$v=X_{\mathcal{K}}\cdot \varepsilon \in \mathcal{K} (\subset \mathcal{F}),$
and for any $w \in \mathcal{F}, \; T\lambda(w)\neq 0 $,
$\tau_{\bar{\mathcal{K}}}\cdot T\pi_{/G}\cdot w \neq 0, $
and $\tau_{\bar{\mathcal{K}}}\cdot T\pi_{/G}\cdot T\lambda(w) \neq 0, $
then we have that
\begin{align*}
& \omega_{\bar{\mathcal{K}}}(T\bar{\gamma} \cdot \tilde{X}^\varepsilon, \;
\tau_{\bar{\mathcal{K}}}\cdot T\pi_{/G} \cdot w)
= \omega_{\bar{\mathcal{K}}}(\tau_{\bar{\mathcal{K}}}\cdot T(\pi_{/G} \cdot
\gamma) \cdot \tilde{X}^\varepsilon, \; \tau_{\bar{\mathcal{K}}}\cdot T\pi_{/G} \cdot w )\\
& = \pi^*_{/G}\cdot \omega_{\bar{\mathcal{K}}}
\cdot\tau_{\bar{\mathcal{K}}}(T\gamma \cdot \tilde{X}^\varepsilon, \; w) =
\tau_{\mathcal{U}}\cdot i^*_{\mathcal{M}} \cdot\omega
(T\gamma \cdot T\pi_Q \cdot \tilde{X} \cdot \varepsilon, \; w)\\
& = \tau_{\mathcal{U}}\cdot i^*_{\mathcal{M}} \cdot
\omega(T(\gamma \cdot \pi_Q)\cdot X_{\mathcal{K}}\cdot \varepsilon, \; w) \\
& = \tau_{\mathcal{U}}\cdot i^*_{\mathcal{M}} \cdot
(\omega(X_{\mathcal{K}}\cdot \varepsilon, \; w-T(\gamma \cdot \pi_Q)\cdot w)
- \mathbf{d}\gamma(T\pi_{Q}(X_{\mathcal{K}}\cdot \varepsilon), \; T\pi_{Q}(w)))\\
& = \tau_{\mathcal{U}}\cdot i^*_{\mathcal{M}} \cdot \omega(X_{\mathcal{K}}\cdot
\varepsilon, \; w) - \tau_{\mathcal{U}}\cdot i^*_{\mathcal{M}} \cdot
\omega(X_{\mathcal{K}}\cdot \varepsilon, \; T\lambda \cdot w)
- \tau_{\mathcal{U}}\cdot i^*_{\mathcal{M}} \cdot\mathbf{d}\gamma
(T\pi_{Q}(X_{\mathcal{K}}\cdot \varepsilon), \; T\pi_{Q}(w))\\
& =\pi^*_{/G}\cdot \omega_{\bar{\mathcal{K}}}\cdot \tau_{\bar{\mathcal{K}}}
(X_{\mathcal{K}}\cdot \varepsilon, \;
w) - \pi^*_{/G}\cdot \omega_{\bar{\mathcal{K}}}\cdot \tau_{\bar{\mathcal{K}}}
(X_{\mathcal{K}}\cdot \varepsilon,
\; T\lambda \cdot w)+ \tau_{\mathcal{U}}\cdot i^*_{\mathcal{M}}
\cdot \lambda^* \omega(X_{\mathcal{K}}\cdot \varepsilon, \; w)\\
& = \omega_{\bar{\mathcal{K}}}(\tau_{\bar{\mathcal{K}}}\cdot T\pi_{/G}
(X_{\mathcal{K}}\cdot \varepsilon), \;
\tau_{\bar{\mathcal{K}}}\cdot T\pi_{/G} \cdot w)
- \omega_{\bar{\mathcal{K}}}(\tau_{\bar{\mathcal{K}}}\cdot T\pi_{/G}(X_{\mathcal{K}}\cdot
\varepsilon), \; \tau_{\bar{\mathcal{K}}}\cdot T(\pi_{/G} \cdot\lambda) \cdot w)\\
& \;\;\;\;\;\; +\pi^*_{/G}\cdot \omega_{\bar{\mathcal{K}}}\cdot
\tau_{\bar{\mathcal{K}}}(T\lambda\cdot X_{\mathcal{K}}\cdot \varepsilon, \;
T\lambda \cdot w)\\
& = \omega_{\bar{\mathcal{K}}}(X_{\bar{\mathcal{K}}} \cdot
\bar{\varepsilon}, \; \tau_{\bar{\mathcal{K}}}\cdot T\pi_{/G} \cdot w)-
\omega_{\bar{\mathcal{K}}}(X_{\bar{\mathcal{K}}}\cdot
\bar{\varepsilon}, \;  \tau_{\bar{\mathcal{K}}}\cdot T\bar{\lambda} \cdot w)+ \omega_{\bar{\mathcal{K}}}
(T\bar{\lambda}\cdot X_{\mathcal{K}}\cdot \varepsilon,
\;  \tau_{\bar{\mathcal{K}}}\cdot T\bar{\lambda} \cdot w),
\end{align*}
where we have used that $\tau_{\bar{\mathcal{K}}}\cdot T\pi_{/G}(X_{\mathcal{K}}\cdot \varepsilon)
=X_{\bar{\mathcal{K}}}\cdot \bar{\varepsilon}, $
and $\tau_{\bar{\mathcal{K}}}\cdot T\bar{\gamma}=T\bar{\gamma},$
and $\tau_{\bar{\mathcal{K}}}\cdot T\pi_{/G}\cdot T\lambda
=\tau_{\bar{\mathcal{K}}}\cdot T\bar{\lambda}=T\bar{\lambda},$ since
$\textmd{Im}(T\bar{\gamma})\subset \bar{\mathcal{K}}.$
From the nonholonomic reduced distributional Hamiltonian equation (4.2),
$\mathbf{i}_{X_{\bar{\mathcal{K}}}}\omega_{\bar{\mathcal{K}}} =
\mathbf{d}h_{\bar{\mathcal{K}}},$ we have that $X_{\bar{\mathcal{K}}}
=\tau_{\bar{\mathcal{K}}}\cdot X_{h_{\bar{\mathcal{K}}}},$
where $ X_{h_{\bar{\mathcal{K}}}}$ is the Hamiltonian vector field of
the function $h_{\bar{\mathcal{K}}}: \bar{M}(\subset T^* Q/G)\rightarrow \mathbb{R}.$
Note that
$\varepsilon: T^* Q \rightarrow T^* Q $ is symplectic, and
$\bar{\varepsilon}=\pi_{/G}(\varepsilon): T^* Q \rightarrow T^* Q/G$
is also symplectic along $\bar{\varepsilon}$, and
hence $X_{h_{\bar{\mathcal{K}}}}\cdot \bar{\varepsilon}
= T\bar{\varepsilon} \cdot X_{h_{\bar{\mathcal{K}}} \cdot
\bar{\varepsilon}}, $ along $\bar{\varepsilon}$, and hence
$X_{\bar{\mathcal{K}}}\cdot \bar{\varepsilon}
=\tau_{\bar{\mathcal{K}}}\cdot X_{h_{\bar{\mathcal{K}}}} \cdot\bar{\varepsilon}
= \tau_{\bar{\mathcal{K}}}\cdot T\bar{\varepsilon} \cdot X_{h_{\bar{\mathcal{K}}}
\cdot \bar{\varepsilon}}, $ along $\bar{\varepsilon}$.
Note that
$T\bar{\lambda} \cdot X_\mathcal{K}\cdot \varepsilon
=T\pi_{/G}\cdot T\lambda \cdot X_\mathcal{K}\cdot \varepsilon
=T\pi_{/G}\cdot T\gamma \cdot T\pi_Q\cdot X_\mathcal{K}\cdot \varepsilon
=T\pi_{/G}\cdot T\gamma \cdot T\pi_Q\cdot \tilde{X}\cdot \varepsilon
=T\pi_{/G}\cdot T\lambda\cdot \tilde{X}\cdot \varepsilon
=T\bar{\lambda} \cdot \tilde{X}\cdot \varepsilon.$
Then we have that
\begin{align*}
& \omega_{\bar{\mathcal{K}}}(T\bar{\gamma} \cdot \tilde{X}^\varepsilon, \;
\tau_{\bar{\mathcal{K}}}\cdot T\pi_{/G}\cdot w)-
\omega_{\bar{\mathcal{K}}}(X_{\bar{\mathcal{K}}}\cdot \bar{\varepsilon},
\; \tau_{\bar{\mathcal{K}}}\cdot T\pi_{/G} \cdot w) \nonumber \\
&= - \omega_{\bar{\mathcal{K}}}(X_{\bar{\mathcal{K}}}\cdot
\bar{\varepsilon}, \;  \tau_{\bar{\mathcal{K}}}\cdot T\bar{\lambda} \cdot w)+ \omega_{\bar{\mathcal{K}}}
(T\bar{\lambda}\cdot X_{\mathcal{K}}\cdot \varepsilon,
\;  \tau_{\bar{\mathcal{K}}}\cdot T\bar{\lambda} \cdot w) \\
& =-\omega_{\bar{\mathcal{K}}}(\tau_{\bar{\mathcal{K}}} \cdot X_{h_{\bar{\mathcal{K}}}}\cdot
\bar{\varepsilon}, \; \tau_{\bar{\mathcal{K}}}\cdot T\bar{\lambda} \cdot w)+ \omega_{\bar{\mathcal{K}}}
(T\bar{\lambda}\cdot \tilde{X} \cdot \varepsilon,
\; \tau_{\bar{\mathcal{K}}}\cdot T\bar{\lambda} \cdot w)\\
& = -\omega_{\bar{\mathcal{K}}}(\tau_{\bar{\mathcal{K}}} \cdot T\bar{\varepsilon}
\cdot X_{h_{\bar{\mathcal{K}}} \cdot \bar{\varepsilon}}, \;
 \tau_{\bar{\mathcal{K}}}\cdot T\bar{\lambda} \cdot w)+ \omega_{\bar{\mathcal{K}}}
(T\bar{\lambda}\cdot \tilde{X} \cdot \varepsilon,
\;  \tau_{\bar{\mathcal{K}}}\cdot T\bar{\lambda} \cdot w)\\
& = \omega_{\bar{\mathcal{K}}}(T\bar{\lambda}\cdot \tilde{X} \cdot \varepsilon- \tau_{\bar{\mathcal{K}}}\cdot T\bar{\varepsilon} \cdot X_{h_{\bar{\mathcal{K}}}\cdot \bar{\varepsilon}},
\;  \tau_{\bar{\mathcal{K}}}\cdot T\bar{\lambda} \cdot w).
\end{align*}
Because the nonholonomic reduced distributional two-form
$\omega_{\bar{\mathcal{K}}}$ is non-degenerate, it follows that the equation
$T\bar{\gamma}\cdot \tilde{X}^\varepsilon = X_{\bar{\mathcal{K}}}\cdot
\bar{\varepsilon},$ is equivalent to the equation $T\bar{\lambda}\cdot \tilde{X} \cdot \varepsilon
= \tau_{\bar{\mathcal{K}}}\cdot T\bar{\varepsilon} \cdot X_{h_{\bar{\mathcal{K}}}
\cdot \bar{\varepsilon}}. $
Thus, $\varepsilon$ and $\bar{\varepsilon}$ satisfy the equation
$T\bar{\lambda}\cdot \tilde{X} \cdot \varepsilon
= \tau_{\bar{\mathcal{K}}}\cdot T\bar{\varepsilon} \cdot X_{h_{\bar{\mathcal{K}}}
\cdot \bar{\varepsilon}},$ if and only if they satisfy
the Type II of Hamilton-Jacobi equation
$T\bar{\gamma}\cdot \tilde{X}^\varepsilon
= X_{\bar{\mathcal{K}}}\cdot \bar{\varepsilon}.$
\hskip 0.3cm $\blacksquare$\\

For a given nonholonomic reducible RCH system
$(T^*Q,G,\omega,\mathcal{D},H,F,u)$ with the associated
the distributional RCH system $(\mathcal{K},\omega_{\mathcal{K}},
H_{\mathcal{K}}, F_{\mathcal{K}}, u_{\mathcal{K}})$
and the nonholonomic reduced distributional RCH system
$(\bar{\mathcal{K}},\omega_{\bar{\mathcal{K}}},h_{\bar{\mathcal{K}}},
f_{\bar{\mathcal{K}}}, u_{\bar{\mathcal{K}}})$,
we know that the nonholonomic dynamical vector field
$X_{\mathcal{K}}$ and the nonholonomic reduced dynamical vector field
$X_{\bar{\mathcal{K}}}$ are $\pi_{/G}$-related, that is,
$X_{\bar{\mathcal{K}}}\cdot \pi_{/G}=T\pi_{/G}\cdot X_{\mathcal{K}}.$
Then we can prove the following Theorem 4.4,
which states the relationship between the solutions of Type II of
Hamilton-Jacobi equations and nonholonomic reduction.

\begin{theo}
For a given nonholonomic reducible RCH system
$(T^*Q,G,\omega,\mathcal{D},H,F,u)$ with the associated
the distributional RCH system $(\mathcal{K},\omega_{\mathcal{K}},
H_{\mathcal{K}}, F_{\mathcal{K}}, u_{\mathcal{K}})$
and the nonholonomic reduced distributional RCH system
$(\bar{\mathcal{K}},\omega_{\bar{\mathcal{K}}},h_{\bar{\mathcal{K}}},
f_{\bar{\mathcal{K}}}, u_{\bar{\mathcal{K}}})$, assume that
$\gamma: Q \rightarrow T^*Q$ is an one-form on $Q$,
and $\varepsilon: T^* Q \rightarrow T^* Q $ is a $G$-invariant symplectic map,
and $\bar{\gamma}=\pi_{/G}(\gamma): Q \rightarrow T^* Q/G $, and
$\bar{\varepsilon}=\pi_{/G}(\varepsilon): T^* Q \rightarrow T^* Q/G. $
Under the hypotheses and notations of Theorem 4.3, then we have that
$\varepsilon$ is a solution of the Type II of Hamilton-Jacobi equation, $T\gamma\cdot
\tilde{X}^\varepsilon= X_{\mathcal{K}}\cdot \varepsilon, $ for the distributional
RCH system $(\mathcal{K},\omega_{\mathcal{K}},H_{\mathcal{K}},
F_{\mathcal{K}}, u_{\mathcal{K}})$, if and only if
$\varepsilon$ and $\bar{\varepsilon}$ satisfy the Type II of
Hamilton-Jacobi equation $T\bar{\gamma}\cdot \tilde{X}^\varepsilon =
X_{\bar{\mathcal{K}}}\cdot \bar{\varepsilon}, $ for the nonholonomic reduced
distributional RCH system $ (\bar{\mathcal{K}},
\omega_{\bar{\mathcal{K}}}, h_{\bar{\mathcal{K}}}, f_{\bar{\mathcal{K}}}, u_{\bar{\mathcal{K}}} ). $
\end{theo}

\noindent{\bf Proof: }At first, for the dynamical vector field of the distributional RCH system
$(\mathcal{K},\omega_{\mathcal {K}}, H_{\mathcal{K}}, F_{\mathcal{K}}, u_{\mathcal{K}})$,
$\tilde{X}=X_{(\mathcal{K},\omega_{\mathcal{K}},
H_{\mathcal{K}}, F_{\mathcal{K}}, u_{\mathcal{K}})}
=X_\mathcal {K}+ F_{\mathcal{K}}+u_{\mathcal{K}}$,
and $F_{\mathcal{K}}=\tau_{\mathcal{K}}\cdot \textnormal{vlift}(F_{\mathcal{M}})X_H$,
and $u_{\mathcal{K}}=\tau_{\mathcal{K}}\cdot \textnormal{vlift}(u_{\mathcal{M}})X_H$,
note that $T\pi_{Q}\cdot \textnormal{vlift}(F_{\mathcal{M}})X_H=T\pi_{Q}\cdot \textnormal{vlift}(u_{\mathcal{M}})X_H=0, $
then we have that $T\pi_{Q}\cdot F_{\mathcal{K}}=T\pi_{Q}\cdot u_{\mathcal{K}}=0,$
and hence $T\pi_{Q}\cdot \tilde{X}\cdot \varepsilon
=T\pi_{Q}\cdot X_{\mathcal{K}}\cdot \varepsilon. $
Next, under the hypotheses and notations of Theorem 4.3,
$\textmd{Im}(\gamma)\subset \mathcal{M},$ and
it is $G$-invariant, $\textmd{Im}(T\gamma)\subset \mathcal{K}, $
and hence $\textmd{Im}(T\bar{\gamma})\subset \bar{\mathcal{K}}, $ in
this case, $\pi^*_{/G}\cdot\omega_{\bar{\mathcal{K}}}\cdot \tau_{\bar{\mathcal{K}}}= \tau_{\mathcal{U}}\cdot
\omega_{\mathcal{M}}= \tau_{\mathcal{U}}\cdot i_{\mathcal{M}}^*\cdot
\omega, $ along $\textmd{Im}(T\bar{\gamma})$.
In addition, from the distributional Hamiltonian equation (4.1),
we have that $X_{\mathcal{K}}=\tau_{\mathcal{K}}\cdot X_H, $
and from the reduced distributional Hamiltonian equation (4.2),
we have that $X_{\bar{\mathcal{K}}}
=\tau_{\bar{\mathcal{K}}}\cdot X_{h_{\bar{\mathcal{K}}}},$
and the nonholonomic dynamical vector field $X_{\mathcal{K}}$ and the
nonholonomic reduced dynamical vector field $X_{\bar{\mathcal{K}}}$ are $\pi_{/G}$-related,
that is, $X_{\bar{\mathcal{K}}}\cdot \pi_{/G}=T\pi_{/G}\cdot
X_{\mathcal{K}}. $ Note that $\varepsilon(\mathcal{M})\subset \mathcal{M},$
and hence $X_H\cdot \varepsilon \in \mathcal{F}$ along $\varepsilon$,
and $\tau_{\mathcal{K}}\cdot X_{H}\cdot \varepsilon
= X_{\mathcal{K}}\cdot \varepsilon \in \mathcal{K}$.
Then $\tau_{\bar{\mathcal{K}}}\cdot T\pi_{/G}(X_{\mathcal{K}}\cdot \varepsilon)
=\tau_{\bar{\mathcal{K}}}\cdot (T\pi_{/G}(X_{\mathcal{K}}))\cdot (\varepsilon)
= \tau_{\bar{\mathcal{K}}}\cdot (X_{\bar{\mathcal{K}}}\cdot \pi_{/G})\cdot (\varepsilon)
= \tau_{\bar{\mathcal{K}}}\cdot X_{\bar{\mathcal{K}}}\cdot \pi_{/G}(\varepsilon)
= X_{\bar{\mathcal{K}}}\cdot \bar{\varepsilon}.$
Thus, using the non-degenerate and nonholonomic reduced distributional two-form
$\omega_{\bar{\mathcal{K}}}$, note that $\tau_{\bar{\mathcal{K}}}\cdot T\bar{\gamma} =T\bar{\gamma},$
for any $w \in \mathcal{F}, \; \tau_{\mathcal{K}} \cdot w\neq 0, $ and
$\tau_{\bar{\mathcal{K}}}\cdot T\pi_{/G}\cdot w \neq 0, $ then
we have that
\begin{align*}
& \omega_{\bar{\mathcal{K}}}(T\bar{\gamma} \cdot \tilde{X}^\varepsilon
- X_{\bar{\mathcal{K}}}\cdot \bar{\varepsilon}, \; \tau_{\bar{\mathcal{K}}}\cdot T\pi_{/G}\cdot w)\\
& = \omega_{\bar{\mathcal{K}}}(T\bar{\gamma} \cdot \tilde{X}^\varepsilon, \;
\tau_{\bar{\mathcal{K}}}\cdot T\pi_{/G}\cdot w)-
\omega_{\bar{\mathcal{K}}}(X_{\bar{\mathcal{K}}}\cdot \bar{\varepsilon},
\; \tau_{\bar{\mathcal{K}}}\cdot T\pi_{/G} \cdot w) \\
& = \omega_{\bar{\mathcal{K}}}(\tau_{\bar{\mathcal{K}}}\cdot T\bar{\gamma}\cdot \tilde{X}^
\varepsilon, \; \tau_{\bar{\mathcal{K}}}\cdot T\pi_{/G}\cdot w)
-\omega_{\bar{\mathcal{K}}}(\tau_{\bar{\mathcal{K}}}\cdot X_{\bar{\mathcal{K}}}\cdot \pi_{/G}\cdot \varepsilon,
\; \tau_{\bar{\mathcal{K}}}\cdot T\pi_{/G}\cdot w)\\
& = \omega_{\bar{\mathcal{K}}}\cdot \tau_{\bar{\mathcal{K}}}(T\pi_{/G}\cdot T\gamma \cdot \tilde{X}^
\varepsilon, \; T\pi_{/G} \cdot w)
- \omega_{\bar{\mathcal{K}}}\cdot \tau_{\bar{\mathcal{K}}}(T\pi_{/G}\cdot
X_{\mathcal{K}}\cdot
\varepsilon, \; T\pi_{/G}\cdot w)\\
& = \pi^*_{/G}\cdot\omega_{\bar{\mathcal{K}}}\cdot \tau_{\bar{\mathcal{K}}}(T\gamma \cdot \tilde{X}^
\varepsilon, \; w)
- \pi^*_{/G}\cdot\omega_{\bar{\mathcal{K}}}\cdot \tau_{\bar{\mathcal{K}}}(X_{\mathcal{K}} \cdot
\varepsilon, \; w)\\
& = \tau_{\mathcal{U}}\cdot i_{\mathcal{M}}^* \cdot
\omega(T\gamma \cdot \tilde{X}^\varepsilon, \; w)
- \tau_{\mathcal{U}}\cdot i_{\mathcal{M}}^* \cdot
\omega(X_{\mathcal{K}} \cdot
\varepsilon, \; w).
\end{align*}
In the case we note that $\tau_{\mathcal{U}}\cdot i_{\mathcal{M}}^* \cdot
\omega=\tau_{\mathcal{K}}\cdot i_{\mathcal{M}}^* \cdot
\omega= \omega_{\mathcal{K}}\cdot \tau_{\mathcal{K}}, $
and
$\tau_{\mathcal{K}}\cdot T\gamma =T\gamma, \; \tau_{\mathcal{K}} \cdot X_{\mathcal{K}}
= X_{\mathcal{K}}$, since $\textmd{Im}(\gamma)\subset
\mathcal{M}, $ and $\textmd{Im}(T\gamma)\subset \mathcal{K}. $
Thus, we have that
\begin{align*}
& \omega_{\bar{\mathcal{K}}}(T\bar{\gamma} \cdot \tilde{X}^\varepsilon
- X_{\bar{\mathcal{K}}}\cdot \bar{\varepsilon}, \; \tau_{\bar{\mathcal{K}}}\cdot T\pi_{/G}\cdot w)\\
& = \omega_{\mathcal{K}}\cdot \tau_{\mathcal{K}}(T\gamma \cdot \tilde{X}^
\varepsilon, \; w)- \omega_{\mathcal{K}}\cdot \tau_{\mathcal{K}}(X_{\mathcal{K}} \cdot \varepsilon, \; w)\\
& = \omega_{\mathcal{K}}(\tau_{\mathcal{K}} \cdot T\gamma \cdot \tilde{X}^
\varepsilon, \; \tau_{\mathcal{K}} \cdot w)
- \omega_{\mathcal{K}}(\tau_{\mathcal{K}} \cdot X_{\mathcal{K}}\cdot \varepsilon, \; \tau_{\mathcal{K}} \cdot w)\\
& = \omega_{\mathcal{K}}(T\gamma \cdot \tilde{X}^
\varepsilon- X_{\mathcal{K}}\cdot \varepsilon, \; \tau_{\mathcal{K}} \cdot w).
\end{align*}
Because the distributional two-form $\omega_{\mathcal{K}}$
and the nonholonomic reduced distributional
two-form $\omega_{\bar{\mathcal{K}}}$ are both non-degenerate,
it follows that the equation
$T\bar{\gamma}\cdot \tilde{X}^\varepsilon=
X_{\bar{\mathcal{K}}}\cdot \bar{\varepsilon}, $
is equivalent to the equation $T\gamma\cdot \tilde{X}^\varepsilon= X_{\mathcal{K}}\cdot \varepsilon. $
Thus, $\varepsilon$ is a solution of the Type II of Hamilton-Jacobi equation
$T\gamma\cdot \tilde{X}^\varepsilon= X_{\mathcal{K}}\cdot \varepsilon, $ for the distributional
RCH system $(\mathcal{K},\omega_{\mathcal {K}},H_{\mathcal{K}},
F_{\mathcal{K}}, u_{\mathcal{K}})$, if and only if
$\varepsilon$ and $\bar{\varepsilon} $ satisfy the Type II of Hamilton-Jacobi
equation $T\bar{\gamma}\cdot \tilde{X}^\varepsilon=
X_{\bar{\mathcal{K}}}\cdot \bar{\varepsilon}, $ for the
nonholonomic reduced distributional RCH system
$(\bar{\mathcal{K}},\omega_{\bar{\mathcal{K}}},h_{\bar{\mathcal{K}}}, f_{\bar{\mathcal{K}}}, u_{\bar{\mathcal{K}}})$.
\hskip 0.3cm $\blacksquare$ \\

It is worthy of noting that,
the Type I of Hamilton-Jacobi equation
$T\bar{\gamma}\cdot \tilde{X}^ \gamma = X_{\bar{\mathcal{K}}}\cdot
\bar{\gamma}. $ is the equation of
the reduced differential one-form $\bar{\gamma}$; and
the Type II of Hamilton-Jacobi equation $T\bar{\gamma}\cdot \tilde{X}^\varepsilon
= X_{\bar{\mathcal{K}}}\cdot \bar{\varepsilon}.$ is the equation of the symplectic
diffeomorphism map $\varepsilon$ and the reduced symplectic
diffeomorphism map $\bar{\varepsilon}$.
If the nonholonomic RCH system with symmetry we considered
has not any the external force and control, that is, $F=0 $ and $W=\emptyset$,
in this case, the nonholonomic RCH system with symmetry
$(T^*Q,G,\omega,\mathcal{D},H,F,W)$ is just the nonholonomic Hamiltonian system
with symmetry $(T^*Q,G,\omega,\mathcal{D},H)$.
and with the canonical symplectic form $\omega$ on $T^*Q$.
From the above Type I and Type II of Hamilton-Jacobi theorems, that is,
Theorem 4.2 and Theorem 4.3, we can get the Theorem 4.2
and Theorem 4.3 in Le\'{o}n and Wang \cite{lewa15}.
It shows that Theorem 4.2 and Theorem 4.3 can be regarded as an extension of two types of
Hamilton-Jacobi theorem for the nonholonomic Hamiltonian system with symmetry given in
Le\'{o}n and Wang \cite{lewa15} to that for the system with the external force and control.

\section{Nonholonomic RCH System with Symmetry and Momentum Map}

As it is well known that momentum map is a very important notion in modern study
of geometric mechanics, and it is a geometric generalization of the classical
linear momentum and angular momentum. A fundamental fact about momentum map is that
if the Hamiltonian $H$ is invariant under the action of a Lie group $G$,
then the momentum map, that is, the vector valued function
$\mathbf{J}: T^\ast Q\rightarrow \mathfrak{g}^\ast ,$ is a constant of the motion
for the Hamiltonian dynamical vector field $X_H$ associated to $H$,
and hence all momentum maps are conserved quantities. Moreover, because momentum map
has infinitesimal equivariance, such that it plays an important role in the
study of reduction theory of Hamiltonian systems with symmetries, see Marsden
\cite{ma92} and Marsden et al.\cite{mamiorpera07, mamora90}.
Now, it is a natural problem what and how we could do,
when the RCH system with symmetry we considered,
has nonholonomic constrains and momentum map,
and the action Lie group $G$ is not Abelian,
and $G_\mu\neq G ,$ where $G_\mu$ is the isotropy subgroup
of coadjoint $G$-action at the point $\mu\in \mathfrak{g}^*$,
and hence the above procedure of nonholonomic reduction
given in $\S 4$ does not work or is not efficient enough.
In this section, we shall consider the nonholonomic RCH
system with symmetry and momentum map, by analyzing carefully the
structures of the nonholonomic dynamical vector fields
and by combining with
the regular symplectic reduction with momentum map,
we give the geometric formulations of
the nonholonomic $R_p$-reduced distributional RCH system
and the nonholonomic $R_o$-reduced distributional RCH system,
where the "regular point reduced" is simply written as $R_p$-reduced,
and the "regular point reduced" is simply written as $R_o$-reduced.
Moreover, we derive precisely the two types of Hamilton-Jacobi equations
for the nonholonomic $R_p$-reduced distributional RCH system
and the nonholonomic $R_o$-reduced distributional RCH
system with respect to momentum map.\\

\subsection{Hamilton-Jacobi equations in the case compatible with regular point reduction}

It is worthy of noting that the regular point
symplectic reduction for the Hamiltonian system with symmetry and
coadjoint equivariant momentum map was set up by famous professors
Jerrold E. Marsden and Alan Weinstein, which is called
Marsden-Weinstein reduction, and it is a very important subject
in the research of geometric mechanics, see
Abraham and Marsden \cite{abma78}, Arnold \cite{ar89},
Libermann and Marle \cite{lima87}, Marsden \cite{ma92}, Marsden et al.
\cite{mamiorpera07, mamora90}, Marsden and Perlmutter \cite{mape00},
Marsden and Ratiu \cite{mara99}, Marsden and
Weinstein \cite{mawe74},  Meyer \cite{me73},
Nijmeijer and Van der Schaft \cite {nivds90} and Ortega and Ratiu \cite{orra04}
and so on, for more details and development.\\

In this subsection, at first, we consider the regular point reducible RCH system,
which is given by Marsden et al \cite{mawazh10}.
Let $Q$ be a smooth manifold and $T^\ast Q$ its cotangent bundle
with the symplectic form $\omega$. Let $\Phi: G \times Q\rightarrow Q$
be a smooth left action of the Lie group $G$ on $Q$, which is free
and proper. Assume that the cotangent lifted left action
$\Phi^{T^\ast}:G\times T^\ast Q\rightarrow T^\ast Q$ is symplectic,
free and proper, and admits an $\operatorname{Ad}^\ast$-equivariant
momentum map $\mathbf{J}:T^\ast Q\rightarrow \mathfrak{g}^\ast$,
where $\mathfrak{g}$ is a Lie algebra of $G$ and $\mathfrak{g}^\ast$
is the dual of $\mathfrak{g}$. Let $\mu\in\mathfrak{g}^\ast$ be a
regular value of $\mathbf{J}$ and denote by $G_\mu$ the isotropy
subgroup of the coadjoint $G$-action at the point
$\mu\in\mathfrak{g}^\ast$, which is defined by $G_\mu=\{g\in
G|\operatorname{Ad}_g^\ast \mu=\mu \}$. Since $G_\mu (\subset G)$
acts freely and properly on $Q$ and on $T^\ast Q$, then
$Q_\mu=Q/G_\mu$ is a smooth manifold and that the canonical
projection $\rho_\mu:Q\rightarrow Q_\mu$ is a surjective submersion.
It follows that $G_\mu$ acts also freely and properly on
$\mathbf{J}^{-1}(\mu)$, so that the space $(T^\ast
Q)_\mu=\mathbf{J}^{-1}(\mu)/G_\mu$ is a symplectic manifold with
the $R_p$-reduced symplectic form $\omega_\mu$ uniquely characterized
by the following relation
\begin{equation}
\pi_\mu^\ast \omega_\mu=i_\mu^\ast \omega.
\label{5.1} \end{equation}
The map
$i_\mu:\mathbf{J}^{-1}(\mu)\rightarrow T^\ast Q$ is the inclusion
and $\pi_\mu:\mathbf{J}^{-1}(\mu)\rightarrow (T^\ast Q)_\mu$ is the
projection. The pair $((T^\ast Q)_\mu,\omega_\mu)$ is called
the $R_p$-reduced space of $(T^\ast Q,\omega)$ at $\mu$,
which is symplectic diffeomorphic to a symplectic fiber bundle. Let
$H: T^\ast Q\rightarrow \mathbb{R}$ be a $G$-invariant Hamiltonian,
the flow $F_t$ of the Hamiltonian vector field $X_H$ leaves the
connected components of $\mathbf{J}^{-1}(\mu)$ invariant and
commutes with the $G$-action, so it induces a flow $f_t^\mu$ on
$(T^\ast Q)_\mu$, defined by $f_t^\mu\cdot \pi_\mu=\pi_\mu \cdot
F_t\cdot i_\mu$, and the vector field $X_{h_\mu}$ generated by the
flow $f_t^\mu$ on $((T^\ast Q)_\mu,\omega_\mu)$ is Hamiltonian with
the associated $R_p$-reduced Hamiltonian function
$h_\mu:(T^\ast Q)_\mu\rightarrow \mathbb{R}$ defined by
$h_\mu\cdot\pi_\mu=H\cdot i_\mu$, and the Hamiltonian vector fields
$X_H$ and $X_{h_\mu}$ are $\pi_\mu$-related.
Moreover, assume that the fiber-preserving map $F:T^\ast Q\rightarrow T^\ast
Q$ and the control subset $W$ of\; $T^\ast Q$ are both $G$-invariant.
In order to get the $R_p$-reduced RCH system, we also assume that
$F(\mathbf{J}^{-1}(\mu))\subset \mathbf{J}^{-1}(\mu)£¬$ and $W \cap
\mathbf{J}^{-1}(\mu)\neq \emptyset $.
Thus, we can introduce a regular point
reducible RCH system as follows, see Marsden et al \cite{mawazh10}
and Wang \cite{wa18}.

\begin{defi}
(Regular Point Reducible RCH System) A 6-tuple $(T^\ast Q, G, \omega, H, F, W)$ with
the canonical symplectic form $\omega$ on $T^*Q$, where the Hamiltonian $H:T^\ast Q\rightarrow
\mathbb{R}$, the fiber-preserving map $F:T^\ast Q\rightarrow T^\ast
Q$ and the fiber submanifold $W$ of\; $T^\ast Q$ are all
$G$-invariant, is called a regular point reducible RCH system, if
there exists a point $\mu\in\mathfrak{g}^\ast$, which is a regular
value of the momentum map $\mathbf{J}$, such that the regular point
reduced system, that is, the 5-tuple $((T^\ast Q)_\mu,
\omega_\mu,h_\mu,f_\mu,W_\mu)$, where $(T^\ast
Q)_\mu=\mathbf{J}^{-1}(\mu)/G_\mu$, $\pi_\mu^\ast
\omega_\mu=i_\mu^\ast\omega$, $h_\mu\cdot \pi_\mu=H\cdot i_\mu$,
$F(\mathbf{J}^{-1}(\mu))\subset \mathbf{J}^{-1}(\mu) $, $f_\mu\cdot
\pi_\mu=\pi_\mu \cdot F\cdot i_\mu$, $W \cap
\mathbf{J}^{-1}(\mu)\neq \emptyset $, $W_\mu=\pi_\mu(W\cap
\mathbf{J}^{-1}(\mu))$, is an RCH system, which is
simply written as $R_p$-reduced RCH system. Where $((T^\ast
Q)_\mu,\omega_\mu)$ is the $R_p$-reduced space, the function $h_\mu:(T^\ast
Q)_\mu\rightarrow \mathbb{R}$ is called the $R_p$-reduced Hamiltonian, the
fiber-preserving map $f_\mu:(T^\ast Q)_\mu\rightarrow (T^\ast
Q)_\mu$ is called the $R_p$-reduced (external) force map, $W_\mu$ is a
fiber submanifold of \;$(T^\ast Q)_\mu$ and is called the $R_p$-reduced
control subset.
\end{defi}

In the following we consider that a nonholonomic RCH system
with symmetry and momentum map
is 8-tuple $(T^*Q,G, \omega,\mathbf{J},\mathcal{D},H,F,W)$,
where $\omega$ is the canonical symplectic form on $T^* Q$,
and the Lie group $G$, which may not be Abelian, acts smoothly by the left on $Q$,
its tangent lifted action on $TQ$ and its cotangent lifted action on $T^\ast Q$,
and $\mathcal{D}\subset TQ$ is a
$\mathcal{D}$-completely and $\mathcal{D}$-regularly nonholonomic
constraint of the system, and $\mathcal{D}$, $H, F$ and $W$ are all
$G$-invariant. Thus, the nonholonomic RCH system with symmetry and momentum map
is a regular point reducible RCH system with $G$-invariant
nonholonomic constraint $\mathcal{D}$.
Moreover, in the following we shall give
carefully a geometric formulation of the $\mathbf{J}$-nonholonomic
$R_p$-reduced distributional RCH system, by using momentum map and the
nonholonomic reduction compatible with regular point reduction.\\

Note that the Legendre transformation $\mathcal{F}L: TQ \rightarrow T^*Q $
is a fiber-preserving map, and $\mathcal{D}\subset TQ$ is $G$-invariant
for the tangent lifted left action $\Phi^{T}: G\times TQ\rightarrow TQ, $
then the constraint submanifold
$\mathcal{M}=\mathcal{F}L(\mathcal{D})\subset T^*Q$ is
$G$-invariant for the cotangent lifted left action $\Phi^{T^\ast}:
G\times T^\ast Q\rightarrow T^\ast Q$.
For the nonholonomic RCH system with symmetry
and momentum map  $(T^*Q,G, \omega,\mathbf{J},\mathcal{D},H,F,W)$,
in the same way, we define the distribution $\mathcal{F}$, which is the pre-image of the
nonholonomic constraints $\mathcal{D}$ for the map $T\pi_Q: TT^* Q
\rightarrow TQ$, that is, $\mathcal{F}=(T\pi_Q)^{-1}(\mathcal{D})$,
and the distribution $\mathcal{K}=\mathcal{F} \cap T\mathcal{M}$.
Moreover, we can also define the distributional two-form $\omega_\mathcal{K}$,
which is induced from the canonical symplectic form $\omega$ on $T^* Q$, that is,
$\omega_\mathcal{K}= \tau_{\mathcal{K}}\cdot \omega_{\mathcal{M}},$ and
$\omega_{\mathcal{M}}= i_{\mathcal{M}}^* \omega $.
If the admissibility condition $\mathrm{dim}\mathcal{M}=
\mathrm{rank}\mathcal{F}$ and the compatibility condition
$T\mathcal{M}\cap \mathcal{F}^\bot= \{0\}$ hold, then
$\omega_\mathcal{K}$ is non-degenerate as a
bilinear form on each fibre of $\mathcal{K}$, there exists a vector
field $X_\mathcal{K}$ on $\mathcal{M}$ which takes values in the
constraint distribution $\mathcal{K}$, such that for the function $H_\mathcal{K}$,
the following distributional Hamiltonian equation holds, that is,
\begin{align}
\mathbf{i}_{X_\mathcal{K}}\omega_\mathcal{K}
=\mathbf{d}H_\mathcal{K},
\label{5.2} \end{align}
where the function $H_{\mathcal{K}}$ satisfies
$\mathbf{d}H_{\mathcal{K}}= \tau_{\mathcal{K}}\cdot \mathbf{d}H_{\mathcal {M}}$,
and $H_\mathcal{M}= \tau_{\mathcal{M}}\cdot H$
is the restriction of $H$ to $\mathcal{M}$, and
from the equation (5.2), we have that
$X_{\mathcal{K}}=\tau_{\mathcal{K}}\cdot X_H $.\\

Since the nonholonomic RCH system with symmetry and momentum map
is a regular point reducible RCH system with $G$-invariant
nonholonomic constraint $\mathcal{D}$,
for a regular value $\mu\in\mathfrak{g}^\ast$ of the
momentum map $\mathbf{J}:T^\ast Q\rightarrow \mathfrak{g}^\ast$,
we shall assume that the constraint submanifold $\mathcal{M}$
is clean intersection with $\mathbf{J}^{-1}(\mu)$, that is,
$\mathcal{M} \cap \mathbf{J}^{-1}(\mu)\neq \emptyset$.
Note that $\mathcal{M}$ is also $G_\mu (\subset G)$ action
invariant, and so is $\mathbf{J}^{-1}(\mu)$, because $\mathbf{J}$ is
$\operatorname{Ad}^\ast$-equivariant. It follows that the quotient
space $\mathcal{M}_\mu =(\mathcal{M}\cap \mathbf{J}^{-1}(\mu))
/G_\mu \subset (T^\ast Q)_\mu$ of the $G_\mu$-orbit in
$\mathcal{M}\cap \mathbf{J}^{-1}(\mu)$, is a smooth manifold with
projection $\pi_\mu: \mathcal{M}\cap \mathbf{J}^{-1}(\mu)
\rightarrow \mathcal{M}_\mu$ which is a surjective submersion.
Denote $i_{\mathcal{M}_\mu}: \mathcal{M}_\mu\rightarrow
(T^*Q)_\mu, $ and $\omega_{\mathcal{M}_\mu}= i_{\mathcal{M}_\mu}^*
\omega_\mu $, that is, the symplectic form
$\omega_{\mathcal{M}_\mu}$ is induced from the $R_p$-reduced symplectic
form $\omega_\mu$ on $(T^* Q)_\mu$ given in (5.1), where $i_{\mathcal{M}_\mu}^*:
T^*(T^*Q)_\mu \rightarrow T^*\mathcal{M}_\mu. $ Moreover, the
distribution $\mathcal{F}$ is pushed down to a distribution
$\mathcal{F}_\mu= T\pi_\mu\cdot \mathcal{F}$ on $(T^\ast Q)_\mu$,
and we define $\mathcal{K}_\mu=\mathcal{F}_\mu \cap
T\mathcal{M}_\mu$. Assume that $\omega_{\mathcal{K}_\mu}=
\tau_{\mathcal{K}_\mu}\cdot \omega_{\mathcal{M}_\mu}$ is the
restriction of the symplectic form $\omega_{\mathcal{M}_\mu}$ on
$T^*\mathcal{M}_\mu$ fibrewise to the distribution $\mathcal{K}_\mu$,
where $\tau_{\mathcal{K}_\mu}$ is the restriction map to distribution
$\mathcal{K}_\mu$. The $\omega_{\mathcal{K}_\mu}$ is not
a true two-form on a manifold, which is called
as a $\mathbf{J}$-nonholonomic $R_p$-reduced
distributional two-form to avoid any confusion.\\

From the above construction we know that,
if the admissibility condition $\mathrm{dim}\mathcal{M}_\mu=
\mathrm{rank}\mathcal{F}_\mu$ and the compatibility condition
$T\mathcal{M}_\mu\cap \mathcal{F}_\mu^\bot= \{0\}$ hold, where
$\mathcal{F}_\mu^\bot$ denotes the symplectic orthogonal of
$\mathcal{F}_\mu$ with respect to the $R_p$-reduced symplectic form
$\omega_\mu$, then $\omega_{\mathcal{K}_\mu}$ is non-degenerate
as a bilinear form on each fibre of $\mathcal{K}_\mu$,
and hence there exists a vector field $X_{\mathcal{K}_\mu}$
on $\mathcal{M}_\mu$, which takes values in the constraint
distribution $\mathcal{K}_\mu$, such that for the function
$h_{\mathcal {K}_\mu}$, the $\mathbf{J}$-nonholonomic
$R_p$-reduced distributional Hamiltonian equation holds, that is,
\begin{align}
\mathbf{i}_{X_{\mathcal{K}_\mu}}\omega_{\mathcal{K}_\mu}
=\mathbf{d}h_{\mathcal{K}_\mu},
\label{5.3} \end{align}
where $\mathbf{d}h_{\mathcal{K}_\mu}$ is the restriction
of $\mathbf{d}h_{\mathcal{M}_\mu}$ to $\mathcal{K}_\mu$,
and the function $h_{\mathcal {K}_\mu}$ satisfies
$\mathbf{d}h_{\mathcal{K}_\mu}= \tau_{\mathcal{K}_\mu}\cdot \mathbf{d}h_{\mathcal{M}_\mu} $,
and $h_{\mathcal{M}_\mu}= \tau_{\mathcal{M}_\mu}\cdot h_\mu$ is the
restriction of $h_\mu$ to $\mathcal{M}_\mu$,
and $h_\mu$ is the $R_p$-reduced
Hamiltonian function $h_\mu: (T^* Q)_\mu \rightarrow \mathbb{R}$ defined
by $h_\mu\cdot \pi_\mu= H\cdot i_\mu$.
In addition, from the distributional Hamiltonian equation (5.2),
$\mathbf{i}_{X_\mathcal{K}}\omega_\mathcal{K}=\mathbf{d}H_\mathcal{K},$
we have that $X_{\mathcal{K}}=\tau_{\mathcal{K}}\cdot X_H, $
and from the $\mathbf{J}$-nonholonomic $R_p$-reduced
distributional Hamiltonian equation (5.3),
$\mathbf{i}_{X_{\mathcal{K_\mu}}}\omega_{\mathcal{K_\mu}}
=\mathbf{d}h_{\mathcal{K_\mu}}$, we have that
$X_{\mathcal{K_\mu}}
=\tau_{\mathcal{K_\mu}}\cdot X_{h_{\mathcal{K_\mu}}},$
where $ X_{h_{\mathcal{K_\mu}}}$ is the Hamiltonian vector field of
the function $h_{\mathcal{K_\mu}},$
and the vector fields $X_{\mathcal{K}}$
and $X_{\mathcal{K_\mu}}$ are $\pi_{\mu}$-related,
that is, $X_{\mathcal{K_\mu}}\cdot \pi_{\mu}=T\pi_{\mu}\cdot X_{\mathcal{K}}.$ \\

Moreover, if considering the external force $F$ and control subset $W$,
and we define the vector fields $F_\mathcal{K}
=\tau_{\mathcal{K}}\cdot \textnormal{vlift}(F_{\mathcal{M}})X_H,$
and for a control law $u\in W$,
$u_\mathcal{K}= \tau_{\mathcal{K}}\cdot  \textnormal{vlift}(u_{\mathcal{M}})X_H,$
where $F_\mathcal{M}= \tau_{\mathcal{M}}\cdot F$ and
$u_\mathcal{M}= \tau_{\mathcal{M}}\cdot u$ are the restrictions of
$F$ and $u$ to $\mathcal{M}$, that is, $F_\mathcal{K}$ and $u_\mathcal{K}$
are the restrictions of the changes of Hamiltonian vector field $X_H$
under the actions of $F_\mathcal{M}$ and $u_\mathcal{M}$ to $\mathcal{K}$,
then the 5-tuple $(\mathcal{K},\omega_{\mathcal{K}},
H_\mathcal{K}, F_\mathcal{K}, u_\mathcal{K})$
is a distributional RCH system corresponding to the nonholonomic RCH system with symmetry
and momentum map $(T^*Q,G, \omega,\mathbf{J},\mathcal{D},H,F,u)$,
and the dynamical vector field of the distributional RCH system
can be expressed by
\begin{align}
\tilde{X}=X_{(\mathcal{K},\omega_{\mathcal{K}},
H_{\mathcal{K}}, F_{\mathcal{K}}, u_{\mathcal{K}})}
=X_\mathcal {K}+ F_{\mathcal{K}}+u_{\mathcal{K}},
\label{5.4} \end{align}
which is the synthetic
of the nonholonomic dynamical vector field $X_{\mathcal{K}}$ and
the vector fields $F_{\mathcal{K}}$ and $u_{\mathcal{K}}$.
Assume that the vector fields $F_\mathcal{K}$ and $u_\mathcal{K}$
on $\mathcal{M}$ are pushed down to the vector fields
$f_{\mathcal{M}_\mu}=T\pi_\mu \cdot F_\mathcal{K}$ and
$u_{\mathcal{M}_\mu}=T\pi_\mu \cdot u_\mathcal{K}$ on $\mathcal{M}_\mu$.
Then we define that $f_{\mathcal{K}_\mu}=T\pi_{\mathcal{K}_\mu}\cdot f_{\mathcal{M}_\mu}$ and
$u_{\mathcal{K}_\mu}=T\pi_{\mathcal{K}_\mu}\cdot u_{\mathcal{M}_\mu},$
that is, $f_{\mathcal{K}_\mu}$ and
$u_{\mathcal{K}_\mu}$ are the restrictions of
$f_{\mathcal{M}_\mu}$ and $u_{\mathcal{M}_\mu}$ to $\mathcal{K}_\mu$.
In consequence, the 5-tuple $(\mathcal{K}_\mu,\omega_{\mathcal{K}_\mu},
h_{\mathcal{K}_\mu}, f_{\mathcal{K}_\mu}, u_{\mathcal{K}_\mu})$
is a $\mathbf{J}$-nonholonomic $R_p$-reduced distributional RCH system of the nonholonomic
RCH system with symmetry and momentum map
$(T^*Q,G,\omega,\mathbf{J},\mathcal{D},H,F,W)$,
as well as with a control law $u \in W$.
Thus, the geometrical formulation
of the $\mathbf{J}$-nonholonomic $R_p$-reduced distributional RCH
system may be summarized as follows.

\begin{defi} ($\mathbf{J}$-Nonholonomic $R_p$-reduced Distributional RCH System)
Assume that the 8-tuple $(T^*Q,G,\omega,\mathbf{J},\mathcal{D},H,F,W)$
is a nonholonomic RCH system with symmetry and momentum map,
where $\omega$ is the canonical symplectic form on $T^* Q$,
and $\mathcal{D}\subset TQ$ is a $\mathcal{D}$-completely and
$\mathcal{D}$-regularly nonholonomic constraint of the system, and
$\mathcal{D}$, and $H, F$ and $W$ are all $G$-invariant. For a regular value
$\mu\in\mathfrak{g}^\ast$ of the momentum map $\mathbf{J}:T^\ast
Q\rightarrow \mathfrak{g}^\ast$, assume that there exists a
$\mathbf{J}$-nonholonomic $R_p$-reduced distribution
$\mathcal{K}_\mu$, an associated non-degenerate and $\mathbf{J}$-nonholonomic
$R_p$-reduced distributional two-form
$\omega_{\mathcal{K}_\mu}$ and a vector field $X_{\mathcal {K}_\mu}$
on the $\mathbf{J}$-nonholonomic $R_p$-reduced constraint submanifold
$\mathcal{M}_\mu=(\mathcal{M}\cap \mathbf{J}^{-1}(\mu)) /G_\mu, $
where $\mathcal{M}=\mathcal{F}L(\mathcal{D}),$ and $\mathcal{M}\cap
\mathbf{J}^{-1}(\mu)\neq \emptyset, $ and $G_\mu=\{g\in G \;| \;
\operatorname{Ad}_g^\ast \mu=\mu \}$, such that the
$\mathbf{J}$-nonholonomic $R_p$-reduced distributional Hamiltonian
equation (5.6) holds, that is,
$ \mathbf{i}_{X_{\mathcal{K}_\mu}}\omega_{\mathcal{K}_\mu} =
\mathbf{d}h_{\mathcal{K}_\mu}, $
where $\mathbf{d}h_{\mathcal{K}_\mu}$ is the restriction of
$\mathbf{d}h_{\mathcal{M}_\mu}$ to $\mathcal{K}_\mu$,
and the function $h_{\mathcal {K}_\mu}$, and the vector fields $f_{\mathcal {K}_\mu}$
and $u_{\mathcal {K}_\mu}$ are defined above. Then the 5-tuple
$(\mathcal{K}_\mu,\omega_{\mathcal {K}_\mu},h_{\mathcal {K}_\mu},
f_{\mathcal {K}_\mu}, u_{\mathcal {K}_\mu})$ is called a
$\mathbf{J}$-nonholonomic $R_p$-reduced distributional RCH system
of the nonholonomic RCH system with symmetry and momentum map
$(T^*Q,G,\omega,\mathbf{J},\mathcal{D},H, F,u)$
with a control law $u \in W$, and $X_{\mathcal{K}_\mu}$
is called the $\mathbf{J}$-nonholonomic $R_p$-reduced
dynamical vector field. Denote by
\begin{align}
\hat{X}_\mu=X_{(\mathcal{K}_\mu,\omega_{\mathcal{K}_\mu},h_{\mathcal {K}_\mu},
f_{\mathcal {K}_\mu}, u_{\mathcal {K}_\mu})}
=X_{\mathcal{K}_\mu}+ f_{\mathcal{K}_\mu}+u_{\mathcal{K}_\mu}
\label{5.5} \end{align}
is the dynamical vector field of the
$\mathbf{J}$-nonholonomic $R_p$-reduced distributional RCH system
$(\mathcal{K}_\mu,\omega_{\mathcal{K}_\mu}, \\ h_{\mathcal {K}_\mu},
f_{\mathcal {K}_\mu}, u_{\mathcal {K}_\mu})$,
which is the synthetic
of the $\mathbf{J}$-nonholonomic $R_p$-reduced
dynamical vector field $X_{\mathcal{K}_\mu}$ and
the vector fields $F_{\mathcal{K}_\mu}$ and $u_{\mathcal{K}_\mu}$.
Under the above circumstances, we refer to
$(T^*Q,G,\omega,\mathbf{J},\mathcal{D},H, F,u)$ as
a $\mathbf{J}$-nonholonomic point reducible RCH system
with the associated distributional RCH system
$(\mathcal{K},\omega_{\mathcal {K}},H_{\mathcal{K}}, F_{\mathcal{K}}, u_{\mathcal{K}})$
and the $\mathbf{J}$-nonholonomic $R_p$-reduced
distributional RCH system
$(\mathcal{K}_\mu,\omega_{\mathcal{K}_\mu},h_{\mathcal {K}_\mu},\\
f_{\mathcal {K}_\mu}, u_{\mathcal {K}_\mu})$.
\end{defi}

Since the non-degenerate and $\mathbf{J}$-nonholonomic $R_p$-reduced
distributional two-form $\omega_{\mathcal{K}_\mu}$ is not a "true two-form"
on a manifold, and it is not symplectic, and hence
the $\mathbf{J}$-nonholonomic $R_p$-reduced distributional RCH system
$(\mathcal{K}_\mu,\omega_{\mathcal {K}_\mu},h_{\mathcal {K}_\mu},
f_{\mathcal {K}_\mu}, u_{\mathcal {K}_\mu})$
is not a Hamiltonian system, and has no yet generating function,
and hence we can not describe the Hamilton-Jacobi equation for a
$\mathbf{J}$-nonholonomic $R_p$-reduced
distributional RCH system just like as in Theorem 1.1.
But, for a given $\mathbf{J}$-nonholonomic regular point reducible RCH system $(T^*Q,G,\omega,\mathbf{J},\mathcal{D},H, F, u)$ with the associated distributional RCH system
$(\mathcal{K},\omega_{\mathcal {K}},H_{\mathcal{K}}, F_{\mathcal{K}}, u_{\mathcal{K}})$
and the $\mathbf{J}$-nonholonomic $R_p$-reduced distributional RCH system
$(\mathcal{K}_\mu,\omega_{\mathcal {K}_\mu},h_{\mathcal {K}_\mu},
f_{\mathcal {K}_\mu}, u_{\mathcal {K}_\mu})$, by using Lemma 3.3,
we can derive precisely
the geometric constraint conditions of the $\mathbf{J}$-nonholonomic $R_p$-reduced distributional two-form
$\omega_{\mathcal{K}_\mu}$ for the $\mathbf{J}$-nonholonomic regular point reducible dynamical vector field,
that is, the two types of Hamilton-Jacobi equation for the
$\mathbf{J}$-nonholonomic $R_p$-reduced distributional RCH system
$(\mathcal{K}_\mu,\omega_{\mathcal {K}_\mu},h_{\mathcal {K}_\mu},
f_{\mathcal {K}_\mu}, u_{\mathcal {K}_\mu})$.
At first, by using the fact that the one-form $\gamma: Q
\rightarrow T^*Q $ is closed on $\mathcal{D}$ with respect to
$T\pi_Q: TT^* Q \rightarrow TQ, $
and $\textmd{Im}(\gamma)\subset \mathcal{M} \cap
\mathbf{J}^{-1}(\mu), $ and it is $G_\mu$-invariant,
as well as $ \textmd{Im}(T\bar{\gamma}_\mu)\subset \mathcal{K}_\mu, $
we can prove the Type I of
Hamilton-Jacobi theorem for the $\mathbf{J}$-nonholonomic $R_p$-reduced distributional
RCH system. For convenience, the maps involved in the
following theorem and its proof are shown in Diagram-5.
\begin{center}
\hskip 0cm \xymatrix{ \mathbf{J}^{-1}(\mu) \ar[r]^{i_\mu}
& T^* Q  \ar[r]^{\pi_Q}
& Q \ar[d]_{\tilde{X}^\gamma} \ar[r]^{\gamma}
& T^*Q \ar[d]_{\tilde{X}} \ar[r]^{\pi_\mu} & (T^* Q)_\mu \ar[d]_{\hat{X}_\mu}
& \mathcal{M}_\mu  \ar[l]_{i_{\mathcal{M}_\mu}} \ar[d]_{X_{\mathcal{K}_\mu}}\\
& T(T^*Q)  & TQ \ar[l]_{T\gamma} & T(T^*Q) \ar[l]^{T\pi_Q} \ar[r]_{T\pi_\mu}
& T(T^* Q)_\mu \ar[r]^{\tau_{\mathcal{K}_\mu}} & \mathcal{K}_\mu }
\end{center}
$$\mbox{Diagram-5}$$

\begin{theo} (Type I of Hamilton-Jacobi Theorem for a $\mathbf{J}$-Nonholonomic
$R_p$-reduced Distributional RCH System) For a given
$\mathbf{J}$-nonholonomic regular point reducible RCH system
$(T^*Q,G,\omega,\mathbf{J},\\ \mathcal{D},H, F, u)$
with the associated distributional RCH system
$(\mathcal{K},\omega_{\mathcal {K}},H_{\mathcal{K}}, F_{\mathcal{K}}, u_{\mathcal{K}})$
and the $\mathbf{J}$-nonholonomic $R_p$-reduced distributional RCH system
$(\mathcal{K}_\mu,\omega_{\mathcal{K}_\mu}, h_{\mathcal {K}_\mu},
f_{\mathcal {K}_\mu}, u_{\mathcal {K}_\mu})$, assume that $\gamma:
Q \rightarrow T^*Q$ is an one-form on $Q$, and
$\tilde{X}^\gamma = T\pi_{Q}\cdot \tilde{X}\cdot \gamma$,
where $\tilde{X}=X_{(\mathcal{K},\omega_{\mathcal{K}},
H_{\mathcal{K}}, F_{\mathcal{K}}, u_{\mathcal{K}})}
=X_\mathcal {K}+ F_{\mathcal{K}}+u_{\mathcal{K}}$
is the dynamical vector field of the distributional RCH system
$(\mathcal{K},\omega_{\mathcal{K}},
H_{\mathcal{K}}, F_{\mathcal{K}}, u_{\mathcal{K}})$
corresponding to the $\mathbf{J}$-nonholonomic regular point reducible RCH
system with symmetry and momentum map $(T^*Q,G,\omega,\mathbf{J},\mathcal{D},H, F,u)$.
Moreover, assume that $\mu\in\mathfrak{g}^\ast$ is a regular value of the momentum
map $\mathbf{J}$, and $\textmd{Im}(\gamma)\subset \mathcal{M} \cap
\mathbf{J}^{-1}(\mu), $ and it is $G_\mu$-invariant, and
$\bar{\gamma}_\mu=\pi_\mu(\gamma): Q \rightarrow \mathcal{M}_\mu $,
and $ \textmd{Im}(T\bar{\gamma}_\mu)\subset \mathcal{K}_\mu. $
If the one-form $\gamma: Q \rightarrow T^*Q $ is closed on $\mathcal{D}$ with respect to
$T\pi_Q: TT^* Q \rightarrow TQ, $ then
$\bar{\gamma}_\mu$ is a solution of the
equation $T\bar{\gamma}_\mu\cdot
\tilde{X}^\gamma= X_{\mathcal{K}_\mu}\cdot \bar{\gamma}_\mu. $
Here $X_{\mathcal{K}_\mu}$ is the $\mathbf{J}$-nonholonomic $R_p$-reduced
dynamical vector field. The equation
$T\bar{\gamma}_\mu\cdot \tilde{X}^\gamma= X_{\mathcal{K}_\mu}\cdot
\bar{\gamma}_\mu,$ is called the Type I of Hamilton-Jacobi equation for the
$\mathbf{J}$-nonholonomic $R_p$-reduced distributional RCH system
$(\mathcal{K}_\mu,\omega_{\mathcal{K}_\mu},h_{\mathcal {K}_\mu},
f_{\mathcal {K}_\mu}, u_{\mathcal {K}_\mu})$.
\end{theo}

\noindent{\bf Proof: } At first, for the dynamical vector field of the distributional RCH system
$(\mathcal{K},\omega_{\mathcal {K}}, H_{\mathcal{K}}, F_{\mathcal{K}}, u_{\mathcal{K}})$,
$\tilde{X}=X_{(\mathcal{K},\omega_{\mathcal{K}},
H_{\mathcal{K}}, F_{\mathcal{K}}, u_{\mathcal{K}})}
=X_\mathcal {K}+ F_{\mathcal{K}}+u_{\mathcal{K}}$,
and $F_{\mathcal{K}}=\tau_{\mathcal{K}}\cdot \textnormal{vlift}(F_{\mathcal{M}})X_H$,
and $u_{\mathcal{K}}=\tau_{\mathcal{K}}\cdot \textnormal{vlift}(u_{\mathcal{M}})X_H$,
note that $T\pi_{Q}\cdot \textnormal{vlift}(F_{\mathcal{M}})X_H=T\pi_{Q}\cdot \textnormal{vlift}(u_{\mathcal{M}})X_H=0, $
then we have that $T\pi_{Q}\cdot F_{\mathcal{K}}=T\pi_{Q}\cdot u_{\mathcal{K}}=0,$
and hence $T\pi_{Q}\cdot \tilde{X}\cdot \gamma=T\pi_{Q}\cdot X_{\mathcal{K}}\cdot \gamma. $
Moreover, from Theorem 3.4, we know that
$\gamma$ is a solution of the Hamilton-Jacobi equation
$T\gamma\cdot \tilde{X}^\gamma= X_{\mathcal{K}}\cdot \gamma .$ Next, we note that
$\textmd{Im}(\gamma)\subset \mathcal{M} \cap \mathbf{J}^{-1}(\mu), $
and it is $G_\mu$-invariant, in this case,
$\pi_\mu^*\omega_\mu=
i_\mu^*\omega= \omega, $ along $\textmd{Im}(\gamma)$.
On the other hand, because
$\textmd{Im}(T\bar{\gamma}_\mu)\subset \mathcal{K}_\mu, $
then $\omega_{\mathcal{K}_\mu}\cdot
\tau_{\mathcal{K}_\mu}=\tau_{\mathcal{K}_\mu}\cdot
\omega_{\mathcal{M}_\mu}= \tau_{\mathcal{K}_\mu}\cdot
i_{\mathcal{M}_\mu}^* \cdot \omega_\mu, $ along
$\textmd{Im}(T\bar{\gamma}_\mu)$.
From the distributional Hamiltonian equation (5.2),
we have that $X_{\mathcal{K}}= \tau_{\mathcal{K}}\cdot X_H,$
and $\tau_{\mathcal{K}}\cdot X_{H}\cdot \gamma
= X_{\mathcal{K}}\cdot \gamma \in \mathcal{K}$.
Because the vector fields $X_{\mathcal{K}}$
and $X_{\mathcal{K}_\mu}$ are $\pi_{\mu}$-related,
$T\pi_{\mu}(X_{\mathcal{K}})=X_{\mathcal{K}_\mu}\cdot \pi_{\mu}$,
and hence $\tau_{\mathcal{K}_\mu}\cdot T\pi_{\mu}(X_{\mathcal{K}}\cdot \gamma)
=\tau_{\mathcal{K}_\mu}\cdot (T\pi_{\mu}(X_{\mathcal{K}}))\cdot (\gamma)
= \tau_{\mathcal{K}_\mu}\cdot (X_{\mathcal{K}_\mu}\cdot \pi_{\mu})\cdot (\gamma)
= \tau_{\mathcal{K}_\mu}\cdot X_{\mathcal{K}_\mu}\cdot \pi_{\mu}(\gamma)
= X_{\mathcal{K}_\mu}\cdot \bar{\gamma}.$
Thus, using the $\mathbf{J}$-nonholonomic $R_p$-reduced distributional two-form
$\omega_{\mathcal{K}_\mu}$, from Lemma 3.3(ii) and (iii), if we take that $v=
X_{\mathcal{K}}\cdot \gamma \in \mathcal{K} (\subset \mathcal{F}), $
and for any $w \in \mathcal{F}, \; T\lambda(w)\neq 0, $ and
$\tau_{\mathcal{K}_\mu}\cdot T\pi_\mu \cdot w \neq 0, $ then we have that
\begin{align*}
& \omega_{\mathcal{K}_\mu}(T\bar{\gamma}_\mu \cdot \tilde{X}^\gamma, \;
\tau_{\mathcal{K}_\mu}\cdot T\pi_\mu \cdot w)=
\omega_{\mathcal{K}_\mu}(\tau_{\mathcal{K}_\mu}\cdot
T\bar{\gamma}_\mu \cdot \tilde{X}^\gamma, \; \tau_{\mathcal{K}_\mu}\cdot T\pi_\mu \cdot w)\\
& = \tau_{\mathcal{K}_\mu}\cdot \omega_{\mathcal{M}_\mu}(T(\pi_\mu
\cdot\gamma )\cdot \tilde{X}^\gamma, \; T\pi_\mu \cdot w ) =
\tau_{\mathcal{K}_\mu}\cdot i_{\mathcal{M}_\mu}^* \cdot
\omega_\mu(T\pi_\mu \cdot T\gamma \cdot \tilde{X}^\gamma, \; T\pi_\mu
\cdot w )\\
& = \tau_{\mathcal{K}_\mu}\cdot i_{\mathcal{M}_\mu}^*
\cdot \pi_\mu^*\omega_\mu(T\gamma \cdot T\pi_Q \cdot \tilde{X} \cdot
\gamma, \; w) = \tau_{\mathcal{K}_\mu}\cdot i_{\mathcal{M}_\mu}^*
\cdot \omega(T(\gamma \cdot \pi_Q) \cdot X_{\mathcal{K}}\cdot \gamma, \; w)\\
& =\tau_{\mathcal{K}_\mu}\cdot i_{\mathcal{M}_\mu}^* \cdot
(\omega(X_{\mathcal{K}}\cdot \gamma, \; w-T(\gamma
\cdot \pi_Q)\cdot w) -\mathbf{d}\gamma(T\pi_{Q}(X_{\mathcal{K}}\cdot \gamma), \; T\pi_{Q}(w)))\\
& = \tau_{\mathcal{K}_\mu}\cdot i_{\mathcal{M}_\mu}^*
\cdot\pi_\mu^*\omega_\mu(X_{\mathcal{K}}\cdot \gamma, \; w) -
\tau_{\mathcal{K}_\mu}\cdot i_{\mathcal{M}_\mu}^*
\cdot\pi_\mu^*\omega_\mu(X_{\mathcal{K}}\cdot
\gamma, \; T(\gamma \cdot \pi_Q) \cdot w)\\
& \;\;\;\;\;\; -\tau_{\mathcal{K}_\mu}\cdot i_{\mathcal{M}_\mu}^* \cdot\mathbf{d}\gamma
(T\pi_{Q}(X_{\mathcal{K}}\cdot \gamma), \; T\pi_{Q}(w))\\
& = \tau_{\mathcal{K}_\mu}\cdot i_{\mathcal{M}_\mu}^*
\cdot\omega_\mu(T\pi_\mu \cdot (X_{\mathcal{K}}\cdot \gamma), \; T\pi_\mu \cdot w) -
\tau_{\mathcal{K}_\mu}\cdot i_{\mathcal{M}_\mu}^* \cdot
\omega_\mu(T\pi_\mu \cdot (X_{\mathcal{K}}\cdot \gamma), \; T(\pi_\mu
\cdot\gamma) \cdot T\pi_{Q}(w))\\
& \;\;\;\;\;\; -\tau_{\mathcal{K}_\mu}\cdot
i_{\mathcal{M}_\mu}^* \cdot\mathbf{d}\gamma(T\pi_{Q}(X_{\mathcal{K}}\cdot
\gamma), \; T\pi_{Q}(w))\\
& = \omega_{\mathcal{K}_\mu}\cdot
\tau_{\mathcal{K}_\mu}(T\pi_\mu \cdot (X_{\mathcal{K}}\cdot \gamma), \; T\pi_\mu \cdot w)-
\omega_{\mathcal{K}_\mu}\cdot \tau_{\mathcal{K}_\mu}
(T\pi_\mu \cdot (X_{\mathcal{K}}\cdot \gamma), \; T\bar{\gamma}_\mu \cdot T\pi_{Q}(w))\\
& \;\;\;\;\;\; -\tau_{\mathcal{K}_\mu}\cdot i_{\mathcal{M}_\mu}^* \cdot\mathbf{d}\gamma
(T\pi_{Q}(X_{\mathcal{K}}\cdot \gamma), \; T\pi_{Q}(w))\\
& = \omega_{\mathcal{K}_\mu}( \tau_{\mathcal{K}_\mu}\cdot T\pi_\mu \cdot (X_{\mathcal{K}}\cdot \gamma), \; \tau_{\mathcal{K}_\mu}\cdot T\pi_\mu
\cdot w) - \omega_{\mathcal{K}_\mu}(\tau_{\mathcal{K}_\mu}\cdot
T\pi_\mu \cdot (X_{\mathcal{K}}\cdot \gamma), \;
\tau_{\mathcal{K}_\mu}\cdot T\bar{\gamma}_\mu \cdot T\pi_{Q}(w))\\
& \;\;\;\;\;\; -\tau_{\mathcal{K}_\mu}\cdot i_{\mathcal{M}_\mu}^* \cdot\mathbf{d}\gamma
(T\pi_{Q}(X_{\mathcal{K}}\cdot \gamma), \; T\pi_{Q}(w))\\
& = \omega_{\mathcal{K}_\mu}(X_{\mathcal{K}_\mu}\cdot
\bar{\gamma}_\mu, \; \tau_{\mathcal{K}_\mu}\cdot T\pi_\mu \cdot w)
- \omega_{\mathcal{K}_\mu}(X_{\mathcal{K}_\mu} \cdot
\bar{\gamma}_\mu, \; \tau_{\mathcal{K}_\mu}\cdot T\bar{\gamma}_\mu
\cdot T\pi_{Q}(w))\\
& \;\;\;\;\;\; -\tau_{\mathcal{K}_\mu}\cdot i_{\mathcal{M}_\mu}^*
\cdot\mathbf{d}\gamma(T\pi_{Q}(X_{\mathcal{K}}\cdot \gamma), \; T\pi_{Q}(w)),
\end{align*}
where we have used that $ \tau_{\mathcal{K}_\mu}\cdot T\bar{\gamma}_\mu=
T\bar{\gamma}_\mu, $ since $\textmd{Im}(T\bar{\gamma}_\mu)\subset \mathcal{K}_\mu, $ and $\tau_{\mathcal{K}_\mu}\cdot T\pi_\mu \cdot (X_{\mathcal{K}}\cdot \gamma)
= X_{\mathcal{K}_\mu}\cdot \bar{\gamma}_\mu. $
If the one-form $\gamma: Q \rightarrow T^*Q $ is closed on $\mathcal{D}$ with respect to
$T\pi_Q: TT^* Q \rightarrow TQ, $ then we have that
$\mathbf{d}\gamma(T\pi_{Q}(X_{\mathcal{K}}\cdot \gamma), \; T\pi_{Q}(w))=0, $
since $X_{\mathcal{K}}\cdot \gamma, \; w \in \mathcal{F},$ and
$T\pi_{Q}(X_{\mathcal{K}}\cdot \gamma), \; T\pi_{Q}(w) \in \mathcal{D}, $ and hence
$$
\tau_{\mathcal{K}_\mu}\cdot
i_{\mathcal{M}_\mu}^* \cdot\mathbf{d}\gamma(T\pi_{Q}(X_{\mathcal{K}}\cdot \gamma),
\; T\pi_{Q}(w))=0,
$$
and
\begin{align}
& \omega_{\mathcal{K}_\mu}(T\bar{\gamma}_\mu \cdot \tilde{X}^\gamma, \;
\tau_{\mathcal{K}_\mu}\cdot T\pi_\mu \cdot w)- \omega_{\mathcal{K}_\mu}(X_{\mathcal{K}_\mu}\cdot
\bar{\gamma}_\mu, \; \tau_{\mathcal{K}_\mu} \cdot T\pi_\mu \cdot w) \nonumber \\
& = -\omega_{\mathcal{K}_\mu}(X_{\mathcal{K}_\mu} \cdot
\bar{\gamma}_\mu, \; \tau_{\mathcal{K}_\mu}\cdot T\bar{\gamma}_\mu
\cdot T\pi_{Q}(w)).
\label {5.6}\end{align}
If $\bar{\gamma}_\mu$ satisfies the equation
$T\bar{\gamma}_\mu\cdot \tilde{X}^\gamma= X_{\mathcal{K}_\mu}\cdot
\bar{\gamma}_\mu ,$
from Lemma 3.3(i) we know that the right side of (5.6) becomes that
\begin{align*}
& -\omega_{\mathcal{K}_\mu}(X_{\mathcal{K}_\mu} \cdot
\bar{\gamma}_\mu, \; \tau_{\mathcal{K}_\mu}\cdot T\bar{\gamma}_\mu
\cdot T\pi_{Q}(w))\\
& = -\omega_{\mathcal{K}_\mu}(T\bar{\gamma}_\mu\cdot \tilde{X}^\gamma,
\; \tau_{\mathcal{K}_\mu}\cdot T\bar{\gamma}_\mu \cdot T\pi_{Q}(w))\\
& = -\omega_{\mathcal{K}_\mu}(\tau_{\mathcal{K}_\mu}T\bar{\gamma}_\mu\cdot \tilde{X}^\gamma,
\; \tau_{\mathcal{K}_\mu}\cdot T\bar{\gamma}_\mu \cdot T\pi_{Q}(w))\\
& = -\tau_{\mathcal{K}_\mu}\cdot i_{\mathcal{M}_\mu}^* \cdot
\omega_\mu(T\bar{\gamma}_\mu\cdot T\pi_{Q} \cdot \tilde{X}\cdot
\gamma, \; T\bar{\gamma}_\mu \cdot T\pi_{Q}(w))\\
& = -\tau_{\mathcal{K}_\mu}\cdot i_{\mathcal{M}_\mu}^* \cdot \bar{\gamma}_\mu^* \cdot
\omega_\mu( T\pi_{Q} \cdot X_{\mathcal{K}}\cdot
\gamma, \; T\pi_{Q}(w))\\
& = -\tau_{\mathcal{K}_\mu}\cdot i_{\mathcal{M}_\mu}^* \cdot
 \gamma^* \cdot \pi^*_\mu\cdot \omega_\mu(T\pi_{Q} \cdot X_{\mathcal{K}}\cdot
\gamma, \; T\pi_{Q}(w))\\
& = -\tau_{\mathcal{K}_\mu}\cdot
i_{\mathcal{M}_\mu}^* \cdot\gamma^*\omega( T\pi_{Q}(X_{\mathcal{K}}\cdot\gamma), \; T\pi_{Q}(w))\\
& = \tau_{\mathcal{K}_\mu}\cdot i_{\mathcal{M}_\mu}^* \cdot
\mathbf{d}\gamma(T\pi_{Q}(X_{\mathcal{K}}\cdot\gamma ), \; T\pi_{Q}(w))=0.
\end{align*}
But, because the $\mathbf{J}$-nonholonomic $R_p$-reduced distributional two-form
$\omega_{\mathcal{K}_\mu}$ is non-degenerate,
the left side of (5.6) equals zero, only when
$\bar{\gamma}_\mu$ satisfies the equation $T\bar{\gamma}_\mu\cdot \tilde{X}^\gamma
= X_{\mathcal{K}_\mu}\cdot \bar{\gamma}_\mu .$ Thus,
if the one-form $\gamma: Q \rightarrow T^*Q $ is closed on $\mathcal{D}$ with respect to
$T\pi_Q: TT^* Q \rightarrow TQ, $ then $\bar{\gamma}_\mu$ must be a solution
of the Type I of Hamilton-Jacobi equation
$T\bar{\gamma}_\mu\cdot \tilde{X}^\gamma= X_{\mathcal{K}_\mu}\cdot
\bar{\gamma}_\mu. $
\hskip 0.3cm $\blacksquare$\\

Next, for any $G_\mu$-invariant symplectic map $\varepsilon: T^* Q \rightarrow T^* Q $,
we can prove the following Type II of
Hamilton-Jacobi theorem for the $\mathbf{J}$-nonholonomic
$R_p$-reduced distributional RCH system.
For convenience, the maps involved in the following
theorem and its proof are shown in Diagram-6.
\begin{center}
\hskip 0cm \xymatrix{ \mathbf{J}^{-1}(\mu) \ar[r]^{i_\mu} & T^* Q
\ar[d]_{X_{H\cdot \varepsilon}} \ar[dr]^{\tilde{X}^\varepsilon} \ar[r]^{\pi_Q}
& Q \ar[r]^{\gamma} & T^*Q \ar[d]_{\tilde{X}}
\ar[dr]^{X_{h_{\mathcal {K}_\mu}} \cdot\bar{\varepsilon}} \ar[r]^{\pi_\mu}
& (T^* Q)_\mu \ar[d]^{X_{h_{\mathcal {K}_\mu}}}
& \mathcal{M}_\mu  \ar[l]_{i_{\mathcal{M}_\mu}} \ar[d]_{X_{\mathcal{K}_\mu}}\\
& T(T^*Q)  & TQ \ar[l]^{T\gamma} & T(T^*Q) \ar[l]^{T\pi_Q} \ar[r]_{T\pi_\mu}
& T(T^* Q)_\mu \ar[r]^{\tau_{\mathcal{K}_\mu}} & \mathcal{K}_\mu }
\end{center}
$$\mbox{Diagram-6}$$

\begin{theo} (Type II of Hamilton-Jacobi Theorem for a $\mathbf{J}$-Nonholonomic
$R_p$-reduced Distributional RCH System) For a given
$\mathbf{J}$-nonholonomic regular point reducible RCH system
$(T^*Q,G,\omega,\mathbf{J},\\ \mathcal{D},H, F, u)$
with the associated distributional RCH system
$(\mathcal{K},\omega_{\mathcal {K}},H_{\mathcal{K}}, F_{\mathcal{K}}, u_{\mathcal{K}})$
and the $\mathbf{J}$-nonholonomic $R_p$-reduced distributional RCH system
$(\mathcal{K}_\mu,\omega_{\mathcal{K}_\mu},h_{\mathcal {K}_\mu},
f_{\mathcal {K}_\mu}, u_{\mathcal {K}_\mu})$, assume that $\gamma:
Q \rightarrow T^*Q$ is an one-form on $Q$, and $\lambda=\gamma \cdot
\pi_{Q}: T^* Q \rightarrow T^* Q, $ and for any
symplectic map $\varepsilon:T^* Q \rightarrow T^* Q, $
denote $\tilde{X}^\varepsilon = T\pi_{Q}\cdot \tilde{X}\cdot \varepsilon$,
where $\tilde{X}=X_{(\mathcal{K},\omega_{\mathcal{K}},
H_{\mathcal{K}}, F_{\mathcal{K}}, u_{\mathcal{K}})}
=X_\mathcal {K}+ F_{\mathcal{K}}+u_{\mathcal{K}}$
is the dynamical vector field of the distributional RCH system
$(\mathcal{K},\omega_{\mathcal{K}},
H_{\mathcal{K}}, F_{\mathcal{K}}, u_{\mathcal{K}})$
corresponding to the $\mathbf{J}$-nonholonomic regular point reducible RCH
system with symmetry and momentum map $(T^*Q,G,\omega,\mathbf{J},\mathcal{D},H, F,u)$.
Moreover, assume that $\mu\in\mathfrak{g}^\ast$ is a regular value of the
momentum map $\mathbf{J}$, and $\textmd{Im}(\gamma)\subset \mathcal{M} \cap
\mathbf{J}^{-1}(\mu), $ and that it is $G_\mu$-invariant,
and $\varepsilon$ is $G_\mu$-invariant and
$\varepsilon(\mathcal{M}\cap \mathbf{J}^{-1}(\mu))
\subset \mathcal{M}\cap \mathbf{J}^{-1}(\mu). $ Denote
$\bar{\gamma}_\mu=\pi_\mu(\gamma): Q \rightarrow \mathcal{M}_\mu $,
and $ \textmd{Im}(T\bar{\gamma}_\mu)\subset \mathcal{K}_\mu, $ and
$\bar{\lambda}_\mu=\pi_\mu(\lambda): \mathcal{M}\cap \mathbf{J}^{-1}(\mu)
(\subset T^* Q) \rightarrow \mathcal{M}_\mu, $ and
$\bar{\varepsilon}_\mu=\pi_\mu(\varepsilon): \mathcal{M}\cap \mathbf{J}^{-1}(\mu)
(\subset T^* Q) \rightarrow \mathcal{M}_\mu. $
Then $\varepsilon$ and $\bar{\varepsilon}_\mu$ satisfy the
equation $\tau_{\mathcal{K}_\mu} \cdot T\bar{\varepsilon}(X_{h_{\mathcal {K}_\mu}\cdot \bar{\varepsilon}_\mu})
= T\bar{\lambda}_\mu \cdot \tilde{X} \cdot\varepsilon, $ if and only if
they satisfy the equation $T\bar{\gamma}_\mu\cdot
\tilde{X}^\varepsilon= X_{\mathcal{K}_\mu}\cdot \bar{\varepsilon}_\mu. $
Here $X_{h_{\mathcal {K}_\mu} \cdot\bar{\varepsilon}_\mu}$ is the Hamiltonian vector field of the
function $h_{\mathcal {K}_\mu}\cdot \bar{\varepsilon}_\mu: T^* Q\rightarrow \mathbb{R}, $
and $X_{\mathcal{K}_\mu}$ is the $\mathbf{J}$-nonholonomic $R_p$-reduced
dynamical vector field. The equation
$T\bar{\gamma}_\mu\cdot \tilde{X}^\varepsilon= X_{\mathcal{K}_\mu}\cdot
\bar{\varepsilon}_\mu,$ is called the Type II of Hamilton-Jacobi equation for the
$\mathbf{J}$-nonholonomic $R_p$-reduced distributional RCH system
$(\mathcal{K}_\mu,\omega_{\mathcal{K}_\mu},h_{\mathcal {K}_\mu},
f_{\mathcal {K}_\mu}, u_{\mathcal {K}_\mu})$.
\end{theo}

\noindent{\bf Proof: } In the same way, for the dynamical vector field of the distributional RCH system
$(\mathcal{K},\omega_{\mathcal {K}}, H_{\mathcal{K}}, F_{\mathcal{K}}, u_{\mathcal{K}})$,
$\tilde{X}=X_{(\mathcal{K},\omega_{\mathcal{K}},
H_{\mathcal{K}}, F_{\mathcal{K}}, u_{\mathcal{K}})}
=X_\mathcal {K}+ F_{\mathcal{K}}+u_{\mathcal{K}}$,
and $F_{\mathcal{K}}=\tau_{\mathcal{K}}\cdot \textnormal{vlift}(F_{\mathcal{M}})X_H$,
and $u_{\mathcal{K}}=\tau_{\mathcal{K}}\cdot \textnormal{vlift}(u_{\mathcal{M}})X_H$,
note that $T\pi_{Q}\cdot \textnormal{vlift}(F_{\mathcal{M}})X_H=T\pi_{Q}\cdot \textnormal{vlift}(u_{\mathcal{M}})X_H=0, $
then we have that $T\pi_{Q}\cdot F_{\mathcal{K}}=T\pi_{Q}\cdot u_{\mathcal{K}}=0,$
and hence $T\pi_{Q}\cdot \tilde{X}\cdot \varepsilon
=T\pi_{Q}\cdot X_{\mathcal{K}}\cdot \varepsilon. $
Next, we note that
$\textmd{Im}(\gamma)\subset \mathcal{M}\cap \mathbf{J}^{-1}(\mu), $
and it is $G_\mu$-invariant, in this case,
$\pi_\mu^*\omega_\mu=
i_\mu^*\omega= \omega, $ along $\textmd{Im}(\gamma)$.
On the other hand, because
$\textmd{Im}(T\bar{\gamma}_\mu)\subset \mathcal{K}_\mu, $
then $\omega_{\mathcal{K}_\mu}\cdot
\tau_{\mathcal{K}_\mu}=\tau_{\mathcal{K}_\mu}\cdot
\omega_{\mathcal{M}_\mu}= \tau_{\mathcal{K}_\mu}\cdot
i_{\mathcal{M}_\mu}^* \cdot \omega_\mu, $ along
$\textmd{Im}(T\bar{\gamma}_\mu)$.
Moreover, from the distributional Hamiltonian equation (5.2),
we have that $X_{\mathcal{K}}= \tau_{\mathcal{K}}\cdot X_H.$
Note that $\varepsilon(\mathcal{M})\subset \mathcal{M},$ and
$T\pi_{Q}(X_H\cdot \varepsilon(q,p))\in
\mathcal{D}_{q}, \; \forall q \in Q, \; (q,p) \in \mathcal{M}(\subset T^* Q), $
and hence $X_H\cdot \varepsilon \in \mathcal{F}$ along $\varepsilon$,
and $\tau_{\mathcal{K}}\cdot X_{H}\cdot \varepsilon
= X_{\mathcal{K}}\cdot \varepsilon \in \mathcal{K}$.
Because the vector fields $X_{\mathcal{K}}$
and $X_{\mathcal{K}_\mu}$ are $\pi_{\mu}$-related, then
$T\pi_{\mu}(X_{\mathcal{K}})=X_{\mathcal{K}_\mu}\cdot \pi_{\mu}$,
and hence $\tau_{\mathcal{K}_\mu}\cdot T\pi_{\mu}(X_{\mathcal{K}}\cdot \varepsilon)
=\tau_{\mathcal{K}_\mu}\cdot (T\pi_{\mu}(X_{\mathcal{K}}))\cdot (\varepsilon)
= \tau_{\mathcal{K}_\mu}\cdot (X_{\mathcal{K}_\mu}\cdot \pi_{\mu})\cdot (\varepsilon)
= \tau_{\mathcal{K}_\mu}\cdot X_{\mathcal{K}_\mu}\cdot \pi_{\mu}(\varepsilon)
= X_{\mathcal{K}_\mu}\cdot \bar{\varepsilon}.$
Thus, using the $\mathbf{J}$-nonholonomic $R_p$-reduced distributional two-form
$\omega_{\mathcal{K}_\mu}$, from Lemma 3.3(ii) and (iii), if we take that $v=
X_{\mathcal{K}}\cdot \varepsilon \in \mathcal{K} (\subset \mathcal{F}), $
and for any $w \in \mathcal{F}, \; T\lambda(w)\neq 0, $ and
$\tau_{\mathcal{K}_\mu}\cdot T\pi_\mu \cdot w \neq 0 $,
$\tau_{\mathcal{K}_\mu}\cdot T\pi_\mu \cdot T\lambda(w) \neq 0, $
then we have that
\begin{align*}
& \omega_{\mathcal{K}_\mu}(T\bar{\gamma}_\mu \cdot \tilde{X}^\varepsilon, \;
\tau_{\mathcal{K}_\mu}\cdot T\pi_\mu \cdot w)=
\omega_{\mathcal{K}_\mu}(\tau_{\mathcal{K}_\mu}\cdot
T\bar{\gamma}_\mu \cdot \tilde{X}^\varepsilon, \; \tau_{\mathcal{K}_\mu}\cdot T\pi_\mu \cdot w)\\
& = \tau_{\mathcal{K}_\mu}\cdot \omega_{\mathcal{M}_\mu}(T(\pi_\mu
\cdot\gamma )\cdot \tilde{X}^\varepsilon, \; T\pi_\mu \cdot w ) =
\tau_{\mathcal{K}_\mu}\cdot i_{\mathcal{M}_\mu}^* \cdot
\omega_\mu(T\pi_\mu \cdot T\gamma \cdot T\pi_Q \cdot \tilde{X} \cdot \varepsilon,
\; T\pi_\mu \cdot w )\\
& = \tau_{\mathcal{K}_\mu}\cdot i_{\mathcal{M}_\mu}^*
\cdot \pi_\mu^*\omega_\mu(T\gamma \cdot T\pi_Q \cdot X_{\mathcal{K}}\cdot
\varepsilon, \; w) = \tau_{\mathcal{K}_\mu}\cdot i_{\mathcal{M}_\mu}^*
\cdot \omega(T(\gamma \cdot \pi_Q) \cdot X_{\mathcal{K}}\cdot \varepsilon, \; w)\\
& =\tau_{\mathcal{K}_\mu}\cdot i_{\mathcal{M}_\mu}^* \cdot
(\omega(X_{\mathcal{K}}\cdot \varepsilon, \; w-T(\gamma
\cdot \pi_Q)\cdot w) -\mathbf{d}\gamma(T\pi_{Q}(X_{\mathcal{K}}\cdot \varepsilon), \; T\pi_{Q}(w)))\\
& = \tau_{\mathcal{K}_\mu}\cdot i_{\mathcal{M}_\mu}^* \cdot
\omega(X_{\mathcal{K}}\cdot \varepsilon, \; w) - \tau_{\mathcal{K}_\mu}\cdot
i_{\mathcal{M}_\mu}^* \cdot \omega(X_{\mathcal{K}}\cdot
\varepsilon, \; T\lambda \cdot w)\\
& \;\;\;\;\;\; -\tau_{\mathcal{K}_\mu}\cdot i_{\mathcal{M}_\mu}^*
\cdot\mathbf{d}\gamma(T\pi_{Q}(X_{\mathcal{K}}\cdot \varepsilon), \; T\pi_{Q}(w))\\
& = \tau_{\mathcal{K}_\mu}\cdot i_{\mathcal{M}_\mu}^*
\cdot\pi_\mu^*\omega_\mu(X_{\mathcal{K}}\cdot \varepsilon, \; w) -
\tau_{\mathcal{K}_\mu}\cdot i_{\mathcal{M}_\mu}^*
\cdot\pi_\mu^*\omega_\mu(X_{\mathcal{K}}\cdot
\varepsilon, \; T\lambda \cdot w)\\
& \;\;\;\;\;\; +\tau_{\mathcal{K}_\mu}\cdot i_{\mathcal{M}_\mu}^* \cdot \lambda^* \omega
(X_{\mathcal{K}}\cdot \varepsilon, \; w)\\
& = \tau_{\mathcal{K}_\mu}\cdot i_{\mathcal{M}_\mu}^*
\cdot\omega_\mu(T\pi_\mu(X_{\mathcal{K}}\cdot \varepsilon), \; T\pi_\mu \cdot w) -
\tau_{\mathcal{K}_\mu}\cdot i_{\mathcal{M}_\mu}^* \cdot
\omega_\mu(T\pi_\mu\cdot(X_{\mathcal{K}}\cdot \varepsilon), \; T(\pi_\mu
\cdot\lambda) \cdot w)\\
& \;\;\;\;\;\; +\tau_{\mathcal{K}_\mu}\cdot
i_{\mathcal{M}_\mu}^* \cdot
\pi_\mu^*\omega_\mu(T\lambda\cdot X_{\mathcal{K}}\cdot
\varepsilon, \; T\lambda \cdot w)\\
& = \omega_{\mathcal{K}_\mu}\cdot
\tau_{\mathcal{K}_\mu}(T\pi_\mu(X_{\mathcal{K}}\cdot
\varepsilon), \; T\pi_\mu \cdot w) - \omega_{\mathcal{K}_\mu}\cdot
\tau_{\mathcal{K}_\mu}(T\pi_\mu(X_{\mathcal{K}}\cdot
\varepsilon), \; T\bar{\lambda}_\mu \cdot w)\\
& \;\;\;\;\;\; +\tau_{\mathcal{K}_\mu}\cdot
i_{\mathcal{M}_\mu}^* \cdot
\omega_\mu(T\pi_\mu\cdot T\lambda\cdot X_{\mathcal{K}}\cdot
\varepsilon, \; T\pi_\mu\cdot T\lambda \cdot w)\\
& = \omega_{\mathcal{K}_\mu}( \tau_{\mathcal{K}_\mu}\cdot T\pi_\mu\cdot(X_{\mathcal{K}}\cdot \varepsilon), \; \tau_{\mathcal{K}_\mu}\cdot T\pi_\mu
\cdot w) - \omega_{\mathcal{K}_\mu}(\tau_{\mathcal{K}_\mu}\cdot
T\pi_\mu\cdot(X_{\mathcal{K}}\cdot \varepsilon), \;
\tau_{\mathcal{K}_\mu}\cdot T\bar{\lambda}_\mu \cdot w)\\
& \;\;\;\;\;\; +\omega_{\mathcal{K}_\mu}(\tau_{\mathcal{K}_\mu}\cdot
T\bar{\lambda}_\mu \cdot X_{\mathcal{K}} \cdot
\varepsilon, \; \tau_{\mathcal{K}_\mu}\cdot T\bar{\lambda}_\mu \cdot w)\\
& = \omega_{\mathcal{K}_\mu}(X_{\mathcal{K}_\mu}\cdot
\bar{\varepsilon}_\mu, \; \tau_{\mathcal{K}_\mu} \cdot T\pi_\mu \cdot w)
- \omega_{\mathcal{K}_\mu}(X_{\mathcal {K}_\mu}\cdot
\bar{\varepsilon}_\mu, \; \tau_{\mathcal{K}_\mu} \cdot T\bar{\lambda}_\mu \cdot w)\\
& \;\;\;\;\;\; +\omega_{\mathcal{K}_\mu}( T\bar{\lambda}_\mu
\cdot X_{\mathcal{K}}\cdot \varepsilon, \; \tau_{\mathcal{K}_\mu} \cdot T\bar{\lambda}_\mu \cdot w),
\end{align*}
where we have used that $ \tau_{\mathcal{K}_\mu}\cdot T\bar{\gamma}_\mu=
T\bar{\gamma}_\mu $, $ \tau_{\mathcal{K}_\mu}\cdot T\bar{\lambda}_\mu=
T\bar{\lambda}_\mu, $ since $\textmd{Im}(T\bar{\gamma}_\mu)\subset \mathcal{K}_\mu, $
and $\tau_{\mathcal{K}_\mu}\cdot T\pi_\mu\cdot(X_{\mathcal{K}}\cdot \varepsilon)
= X_{\mathcal{K}_\mu}\cdot \bar{\varepsilon}_\mu. $
From the $\mathbf{J}$-nonholonomic $R_p$-reduced
distributional Hamiltonian equation (5.3),
$\mathbf{i}_{X_{\mathcal{K_\mu}}}\omega_{\mathcal{K_\mu}}
=\mathbf{d}h_{\mathcal{K_\mu}}$, we have that
$X_{\mathcal{K_\mu}}
=\tau_{\mathcal{K_\mu}}\cdot X_{h_{\mathcal{K_\mu}}},$
where $ X_{h_{\mathcal{K_\mu}}}$ is the Hamiltonian vector field of
the function $h_{\mathcal{K_\mu}}.$
Note that $\varepsilon: T^* Q \rightarrow T^* Q $ is symplectic, and
$\pi_\mu^*\omega_\mu= i_\mu^*\omega= \omega, $ along
$\textmd{Im}(\gamma)$, and hence $\bar{\varepsilon}_\mu=
\pi_\mu(\varepsilon): T^* Q \rightarrow (T^* Q)_\mu $ is also symplectic
along $\bar{\varepsilon}_\mu$, and hence $X_{h_{\mathcal {K}_\mu}}\cdot
\bar{\varepsilon}_\mu= T\bar{\varepsilon}_\mu \cdot X_{h_{\mathcal {K}_\mu} \cdot
\bar{\varepsilon}_\mu}, $ along $\bar{\varepsilon}_\mu$, and hence
$X_{\mathcal{K}_\mu}\cdot \bar{\varepsilon}_\mu
=\tau_{\mathcal{K}_\mu}\cdot X_{h_{\mathcal {K}_\mu}} \cdot \bar{\varepsilon}_\mu
= \tau_{\mathcal{K}_\mu}T\bar{\varepsilon}_\mu \cdot X_{h_{\mathcal {K}_\mu} \cdot
\bar{\varepsilon}_\mu}, $ along $\bar{\varepsilon}_\mu.$
Note that
$T\bar{\lambda}_\mu \cdot X_\mathcal{K}\cdot \varepsilon
=T\pi_{\mu}\cdot T\lambda \cdot X_\mathcal{K}\cdot \varepsilon
=T\pi_{\mu}\cdot T\gamma \cdot T\pi_Q\cdot X_\mathcal{K}\cdot \varepsilon
=T\pi_{\mu}\cdot T\gamma \cdot T\pi_Q\cdot \tilde{X}\cdot \varepsilon
=T\pi_{\mu}\cdot T\lambda \cdot \tilde{X}\cdot \varepsilon
=T\bar{\lambda}_\mu \cdot \tilde{X}\cdot \varepsilon.$
Then we have that
\begin{align*}
& \omega_{\mathcal{K}_\mu}(T\bar{\gamma}_\mu \cdot \tilde{X}^\varepsilon, \;
\tau_{\mathcal{K}_\mu}\cdot T\pi_\mu \cdot w)-
\omega_{\mathcal{K}_\mu}(X_{\mathcal{K}_\mu}\cdot \bar{\varepsilon}_\mu,
\; \tau_{\mathcal{K}_\mu} \cdot T\pi_\mu \cdot w) \nonumber \\
& = - \omega_{\mathcal{K}_\mu}(X_{\mathcal{K}_\mu}\cdot \bar{\varepsilon}_\mu,
\; \tau_{\mathcal{K}_\mu}\cdot T\bar{\lambda}_\mu \cdot w)
+\omega_{\mathcal{K}_\mu}(T\bar{\lambda}_\mu \cdot X_{\mathcal{K}}\cdot
\varepsilon, \; \tau_{\mathcal{K}_\mu}\cdot T\bar{\lambda}_\mu \cdot w)\\
&= - \omega_{\mathcal{K}_\mu}(\tau_{\mathcal{K}_\mu}
T\bar{\varepsilon}_\mu \cdot X_{h_{\mathcal {K}_\mu} \cdot \bar{\varepsilon}_\mu},
\; \tau_{\mathcal{K}_\mu}\cdot T\bar{\lambda}_\mu \cdot w)
+\omega_{\mathcal{K}_\mu}(T\bar{\lambda}_\mu \cdot \tilde{X}\cdot
\varepsilon, \; \tau_{\mathcal{K}_\mu}\cdot T\bar{\lambda}_\mu \cdot w)\\
& = \omega_{\mathcal{K}_\mu}(T\bar{\lambda}_\mu \cdot \tilde{X}\cdot
\varepsilon- \tau_{\mathcal{K}_\mu}\cdot T\bar{\varepsilon}_\mu
\cdot X_{h_{\mathcal {K}_\mu} \cdot \bar{\varepsilon}_\mu}, \;
\tau_{\mathcal{K}_\mu}\cdot T\bar{\lambda}_\mu \cdot w).
\end{align*}
Because the $\mathbf{J}$-nonholonomic $R_p$-reduced distributional two-form
$\omega_{\mathcal{K}_\mu}$ is non-degenerate, it follows that the equation
$T\bar{\gamma}_\mu\cdot \tilde{X}^\varepsilon = X_{\mathcal{K}_\mu}\cdot
\bar{\varepsilon}_\mu, $ is equivalent to the equation $T\bar{\lambda}_\mu\cdot \tilde{X}\cdot \varepsilon
= \tau_{\mathcal{K}_\mu}\cdot T\bar{\varepsilon}_\mu \cdot X_{h_{\mathcal {K}_\mu}
\cdot \bar{\varepsilon}_\mu}. $
Thus, $\varepsilon$ and $\bar{\varepsilon}_\mu$ satisfy the equation
$T\bar{\lambda}_\mu\cdot \tilde{X}\cdot \varepsilon
= \tau_{\mathcal{K}_\mu}\cdot T\bar{\varepsilon}_\mu \cdot X_{h_{\mathcal {K}_\mu}
\cdot \bar{\varepsilon}_\mu}, $ if and only if they satisfy
the Type II of Hamilton-Jacobi equation $T\bar{\gamma}_\mu\cdot \tilde{X}^\varepsilon
= X_{\mathcal{K}_\mu}\cdot \bar{\varepsilon}_\mu .$
\hskip 0.3cm $\blacksquare$\\

For a given $\mathbf{J}$-nonholonomic regular point reducible RCH system
$(T^*Q,G,\omega,\mathbf{J},\mathcal{D},H, F, W)$
with the associated distributional RCH system
$(\mathcal{K},\omega_{\mathcal {K}},H_{\mathcal{K}}, F_{\mathcal{K}}, u_{\mathcal{K}})$
and the $\mathbf{J}$-nonholonomic $R_p$-reduced distributional RCH system
$(\mathcal{K}_\mu,\omega_{\mathcal{K}_\mu},h_{\mathcal {K}_\mu},
f_{\mathcal {K}_\mu}, u_{\mathcal {K}_\mu})$, we know that the
nonholonomic dynamical vector field $X_{\mathcal{K}}$ and the $\mathbf{J}$-
nonholonomic $R_p$-reduced dynamical vector field $X_{\mathcal{K_\mu}}$ are $\pi_\mu$-related,
that is, $X_{\mathcal{K_\mu}}\cdot \pi_\mu=T\pi_\mu\cdot X_{\mathcal{K}}\cdot i_\mu.$ Then
we can prove the following Theorem 5.5, which states the
relationship between the solutions of Type II of Hamilton-Jacobi equations and
$\mathbf{J}$-nonholonomic regular point reduction.

\begin{theo}
For a given $\mathbf{J}$-nonholonomic regular point reducible RCH system
$(T^*Q,G,\omega, \mathbf{J}, \mathcal{D}, \\ H, F, u)$
with the associated distributional RCH system
$(\mathcal{K},\omega_{\mathcal {K}},H_{\mathcal{K}}, F_{\mathcal{K}}, u_{\mathcal{K}})$
and the $\mathbf{J}$-nonholonomic $R_p$-reduced distributional RCH system
$(\mathcal{K}_\mu,\omega_{\mathcal{K}_\mu}, h_{\mathcal {K}_\mu},
f_{\mathcal {K}_\mu}, u_{\mathcal {K}_\mu})$, assume that $\gamma:
Q \rightarrow T^*Q$ is an one-form on $Q$, and
$\varepsilon: T^* Q \rightarrow T^* Q $ is a symplectic map,
and $\bar{\gamma}_\mu=\pi_\mu(\gamma): Q \rightarrow \mathcal{M}_\mu $,
and $\bar{\varepsilon}_\mu=\pi_\mu(\varepsilon): \mathcal{M}\cap \mathbf{J}^{-1}(\mu)
(\subset T^*Q) \rightarrow (T^* Q)_\mu $.
Under the hypotheses and notations of Theorem 5.4, then we have that
$\varepsilon$ is a solution of the Type II of Hamilton-Jacobi equation $T\gamma\cdot
\tilde{X}^\varepsilon= X_{\mathcal{K}}\cdot \varepsilon, $ for the distributional
RCH system $(\mathcal{K},\omega_{\mathcal {K}},H_{\mathcal {K}},
F_{\mathcal {K}}, u_{\mathcal {K}})$, if and
only if $\varepsilon$ and $\bar{\varepsilon}_\mu $ satisfy the Type II of Hamilton-Jacobi
equation $T\bar{\gamma}_\mu\cdot \tilde{X}^\varepsilon=
X_{\mathcal{K}_\mu}\cdot \bar{\varepsilon}_\mu, $ for the
$\mathbf{J}$-nonholonomic $R_p$-reduced distributional RCH system
$(\mathcal{K}_\mu,\omega_{\mathcal{K}_\mu},h_{\mathcal {K}_\mu},
f_{\mathcal {K}_\mu}, u_{\mathcal {K}_\mu})$.
\end{theo}

\noindent{\bf Proof: } For the dynamical vector field of the distributional RCH system
$(\mathcal{K},\omega_{\mathcal {K}}, H_{\mathcal{K}}, F_{\mathcal{K}}, u_{\mathcal{K}})$,
$\tilde{X}=X_{(\mathcal{K},\omega_{\mathcal{K}},
H_{\mathcal{K}}, F_{\mathcal{K}}, u_{\mathcal{K}})}
=X_\mathcal {K}+ F_{\mathcal{K}}+u_{\mathcal{K}}$,
and $F_{\mathcal{K}}=\tau_{\mathcal{K}}\cdot \textnormal{vlift}(F_{\mathcal{M}})X_H$,
and $u_{\mathcal{K}}=\tau_{\mathcal{K}}\cdot \textnormal{vlift}(u_{\mathcal{M}})X_H$,
note that $T\pi_{Q}\cdot \textnormal{vlift}(F_{\mathcal{M}})X_H=T\pi_{Q}\cdot \textnormal{vlift}(u_{\mathcal{M}})X_H=0, $
then we have that $T\pi_{Q}\cdot F_{\mathcal{K}}=T\pi_{Q}\cdot u_{\mathcal{K}}=0,$
and hence $T\pi_{Q}\cdot \tilde{X}\cdot \varepsilon
=T\pi_{Q}\cdot X_{\mathcal{K}}\cdot \varepsilon. $
Next, under the hypotheses and notations of Theorem 5.4,
$\textmd{Im}(\gamma)\subset \mathcal{M}\cap \mathbf{J}^{-1}(\mu), $
and it is $G_\mu$-invariant, in this case,
$\pi_\mu^*\omega_\mu=
i_\mu^*\omega= \omega, $ along $\textmd{Im}(\gamma)$.
On the other hand, because
$\textmd{Im}(T\bar{\gamma}_\mu)\subset \mathcal{K}_\mu, $
then $\omega_{\mathcal{K}_\mu}\cdot
\tau_{\mathcal{K}_\mu}=\tau_{\mathcal{K}_\mu}\cdot
\omega_{\mathcal{M}_\mu}= \tau_{\mathcal{K}_\mu}\cdot
i_{\mathcal{M}_\mu}^* \cdot \omega_\mu, $ along
$\textmd{Im}(T\bar{\gamma}_\mu)$.
In addition, from the distributional Hamiltonian equation (5.2),
we have that $X_{\mathcal{K}}=\tau_{\mathcal{K}}\cdot X_H, $
and from the $\mathbf{J}$-nonholonomic $R_p$-reduced
distributional Hamiltonian equation (5.3), we have that
$X_{\mathcal{K_\mu}}
=\tau_{\mathcal{K_\mu}}\cdot X_{h_{\mathcal{K_\mu}}},$
and the vector fields $X_{\mathcal{K}}$
and $X_{\mathcal{K_\mu}}$ are $\pi_{\mu}$-related,
that is, $X_{\mathcal{K_\mu}}\cdot \pi_{\mu}=T\pi_{\mu}\cdot X_{\mathcal{K}}.$
Note that $\varepsilon(\mathcal{M})\subset \mathcal{M},$
and hence $X_H\cdot \varepsilon \in \mathcal{F}$ along $\varepsilon$,
and $\tau_{\mathcal{K}}\cdot X_{H}\cdot \varepsilon
= X_{\mathcal{K}}\cdot \varepsilon \in \mathcal{K}$.
Then $\tau_{\mathcal{K}_\mu}\cdot T\pi_{\mu}(X_{\mathcal{K}}\cdot \varepsilon)
=\tau_{\mathcal{K}_\mu}\cdot (T\pi_{\mu}(X_{\mathcal{K}}))\cdot (\varepsilon)
= \tau_{\mathcal{K}_\mu}\cdot (X_{\mathcal{K}_\mu}\cdot \pi_{\mu})\cdot (\varepsilon)
= \tau_{\mathcal{K}_\mu}\cdot X_{\mathcal{K}_\mu}\cdot \pi_{\mu}(\varepsilon)
= X_{\mathcal{K}_\mu}\cdot \bar{\varepsilon}.$
Thus, using the $\mathbf{J}$-nonholonomic $R_p$-reduced distributional two-form
$\omega_{\mathcal{K}_\mu}$, note that
$\tau_{\mathcal{K}_\mu}\cdot T\bar{\gamma}_\mu =T\bar{\gamma}_\mu,$
for any $w \in \mathcal{F}, \; \tau_{\mathcal{K}}\cdot w\neq 0, $ and
$\tau_{\mathcal{K}_\mu}\cdot T\pi_\mu \cdot w \neq 0 $,
then we have that
\begin{align*}
& \omega_{\mathcal{K}_\mu}(T\bar{\gamma}_\mu \cdot \tilde{X}^\varepsilon
- X_{\mathcal{K}_\mu}\cdot \bar{\varepsilon}_\mu, \;
\tau_{\mathcal{K}_\mu}\cdot T\pi_\mu \cdot w) \\
& = \omega_{\mathcal{K}_\mu}(T\bar{\gamma}_\mu \cdot \tilde{X}^\varepsilon, \;
\tau_{\mathcal{K}_\mu}\cdot T\pi_\mu \cdot w)-
\omega_{\mathcal{K}_\mu}(X_{\mathcal{K}_\mu}\cdot \bar{\varepsilon}_\mu,
\; \tau_{\mathcal{K}_\mu} \cdot T\pi_\mu \cdot w) \\
& = \omega_{\mathcal{K}_\mu}(\tau_{\mathcal{K}_\mu}\cdot T\bar{\gamma}_\mu \cdot \tilde{X}^\varepsilon, \;
\tau_{\mathcal{K}_\mu}\cdot T\pi_\mu \cdot w)-
\omega_{\mathcal{K}_\mu}(\tau_{\mathcal{K}_\mu}\cdot T\pi_{\mu}(X_{\mathcal{K}}\cdot \varepsilon), \;
\tau_{\mathcal{K}_\mu} \cdot T\pi_\mu \cdot w) \\
& = \omega_{\mathcal{K}_\mu}\cdot \tau_{\mathcal{K}_\mu}
(T\pi_\mu \cdot T\gamma \cdot \tilde{X}^\varepsilon, \; T\pi_\mu \cdot w)
-\omega_{\mathcal{K}_\mu}\cdot \tau_{\mathcal{K}_\mu}
(T\pi_\mu \cdot X_{\mathcal{K}}\cdot \varepsilon, \; T\pi_\mu \cdot w)\\
& = \tau_{\mathcal{K}_\mu}\cdot
i_{\mathcal{M}_\mu}^* \cdot \omega_\mu(T\pi_\mu \cdot T\gamma \cdot \tilde{X}^\varepsilon, \;
T\pi_\mu \cdot w) -\tau_{\mathcal{K}_\mu}\cdot
i_{\mathcal{M}_\mu}^* \cdot \omega_\mu
(T\pi_\mu \cdot X_{\mathcal{K}}\cdot \varepsilon, \; T\pi_\mu \cdot w)\\
& = \tau_{\mathcal{K}_\mu}\cdot
i_{\mathcal{M}_\mu}^* \cdot \pi_\mu^*\omega_\mu(T\gamma \cdot \tilde{X}^\varepsilon, \; w)
-\tau_{\mathcal{K}_\mu}\cdot
i_{\mathcal{M}_\mu}^* \cdot \pi_\mu^*\omega_\mu(X_{\mathcal{K}}\cdot \varepsilon, \; w).
\end{align*}
In the case we note that $\tau_{\mathcal{K}_\mu}\cdot
i_{\mathcal{M}_\mu}^* \cdot \pi_\mu^*\omega_\mu
=\tau_{\mathcal{K}}\cdot i_{\mathcal{M}}^* \cdot
\omega= \omega_{\mathcal{K}}\cdot \tau_{\mathcal{K}}, $
and
$\tau_{\mathcal{K}}\cdot T\gamma =T\gamma, \; \tau_{\mathcal{K}} \cdot \tilde{X}= X_{\mathcal{K}}$,
since $\textmd{Im}(\gamma)\subset
\mathcal{M}, $ and $\textmd{Im}(T\gamma)\subset \mathcal{K}. $
Thus, we have that
\begin{align*}
& \omega_{\mathcal{K}_\mu}(T\bar{\gamma}_\mu \cdot \tilde{X}^\varepsilon
- X_{\mathcal{K}_\mu}\cdot \bar{\varepsilon}_\mu, \;
\tau_{\mathcal{K}_\mu}\cdot T\pi_\mu \cdot w) \\
& = \omega_{\mathcal{K}}\cdot \tau_{\mathcal{K}}(T\gamma \cdot \tilde{X}^\varepsilon, \; w)
-\omega_{\mathcal{K}}\cdot \tau_{\mathcal{K}}(X_{\mathcal{K}}\cdot \varepsilon, \; w)\\
& = \omega_{\mathcal{K}}(\tau_{\mathcal{K}}\cdot T\gamma \cdot \tilde{X}^\varepsilon, \; \tau_{\mathcal{K}}\cdot w)
-\omega_{\mathcal{K}}(\tau_{\mathcal{K}}\cdot X_{\mathcal{K}}\cdot \varepsilon, \; \tau_{\mathcal{K}}\cdot w)\\
& = \omega_{\mathcal{K}}(T\gamma \cdot \tilde{X}^\varepsilon, \; \tau_{\mathcal{K}}\cdot w)
-\omega_{\mathcal{K}}(X_{\mathcal{K}}\cdot \varepsilon, \; \tau_{\mathcal{K}}\cdot w)\\
& = \omega_{\mathcal{K}}(T\gamma \cdot \tilde{X}^\varepsilon- X_{\mathcal{K}}\cdot \varepsilon, \; \tau_{\mathcal{K}}\cdot w).
\end{align*}
Because the distributional two-form $\omega_{\mathcal{K}}$ and
the $\mathbf{J}$-nonholonomic $R_p$-reduced distributional
two-form $\omega_{\mathcal{K}_\mu}$ are both non-degenerate,
it follows that the equation
$T\bar{\gamma}_\mu\cdot \tilde{X}^\varepsilon=
X_{\mathcal{K}_\mu}\cdot \bar{\varepsilon}_\mu, $ is equivalent to the equation
$T\gamma\cdot \tilde{X}^\varepsilon= X_{\mathcal{K}}\cdot \varepsilon$. Thus,
$\varepsilon$ is a solution of the Type II of Hamilton-Jacobi equation
$T\gamma\cdot \tilde{X}^\varepsilon= X_{\mathcal{K}}\cdot \varepsilon, $ for the distributional
RCH system $(\mathcal{K},\omega_{\mathcal {K}},H_{\mathcal {K}},
F_{\mathcal {K}}, u_{\mathcal {K}})$, if and only if
$\varepsilon$ and $\bar{\varepsilon}_\mu $ satisfy the Type II of Hamilton-Jacobi
equation $T\bar{\gamma}_\mu\cdot \tilde{X}^\varepsilon=
X_{\mathcal{K}_\mu}\cdot \bar{\varepsilon}_\mu, $ for the
$\mathbf{J}$-nonholonomic $R_p$-reduced distributional RCH system
$(\mathcal{K}_\mu,\omega_{\mathcal{K}_\mu},h_{\mathcal{K}_\mu},
f_{\mathcal {K}_\mu}, u_{\mathcal {K}_\mu})$.  \hskip 0.3cm
$\blacksquare$ \\

\begin{rema}
It is worthy of noting that,
the Type I of Hamilton-Jacobi equation
$T\bar{\gamma}_\mu\cdot \tilde{X}^\gamma= X_{\mathcal{K}_\mu}\cdot
\bar{\gamma}_\mu. $ is the equation of the $\mathbf{J}$-nonholonomic
$R_p$-reduced differential one-form $\bar{\gamma}_\mu$; and
the Type II of Hamilton-Jacobi equation $T\bar{\gamma}_\mu\cdot \tilde{X}^\varepsilon
= X_{\mathcal{K}_\mu}\cdot \bar{\varepsilon}_\mu .$ is the equation of the symplectic
diffeomorphism map $\varepsilon$ and the $\mathbf{J}$-nonholonomic $R_p$-reduced
symplectic diffeomorphism map $\bar{\varepsilon}_\mu. $
If a $\mathbf{J}$-nonholonomic regular point reducible RCH system we considered
$(T^*Q,G,\omega, \mathbf{J}, \mathcal{D}, H, F, u)$ has not any constrains,
in this case, the $\mathbf{J}$-nonholonomic $R_p$-reduced distributional
RCH system is just the $R_p$-reduced RCH system itself.
From the above Type I and Type II of Hamilton-Jacobi theorems, that is,
Theorem 5.3 and Theorem 5.4, we can get the Theorem 3.2
and Theorem 3.3 in Wang \cite{wa13d}.
It shows that Theorem 5.3 and Theorem 5.4 can be regarded as an extension of two types of
Hamilton-Jacobi theorem for the $R_p$-reduced RCH system
to that for the system with nonholonomic context.
If the $\mathbf{J}$-nonholonomic regular point reducible RCH system we considered
$(T^*Q,G,\omega, \mathbf{J}, \mathcal{D}, H, F, u)$
has not any the external force and control, that is, $F=0 $ and $u=0$,
in this case, the $\mathbf{J}$-nonholonomic regular point reducible RCH system
is just the $\mathbf{J}$-nonholonomic regular point reducible
Hamiltonian system $(T^*Q,G,\omega,\mathbf{J},\mathcal{D},H)$.
and with the canonical symplectic form $\omega$ on $T^*Q$.
From the above Type I and Type II of Hamilton-Jacobi theorems, that is,
Theorem 5.3 and Theorem 5.4, we can get the Theorem 5.2
and Theorem 5.3 in Le\'{o}n and Wang \cite{lewa15}.
It shows that Theorem 5.3 and Theorem 5.4 can be regarded as an extension of two types of
Hamilton-Jacobi theorem for the $\mathbf{J}$-nonholonomic regular point reducible
Hamiltonian system to that for the system with the external force and control.
In particular, in this case, if the $\mathbf{J}$-nonholonomic
regular point reducible RCH system we considered has not any constrains,
that is, $F=0, \; u=0 $ and $\mathcal{D}=\emptyset$, then
the $\mathbf{J}$-nonholonomic regular point reducible RCH system
is just a Marsden-Weinstein reducible Hamiltonian system $(T^*Q,G,\omega,H)$
with the canonical symplectic form $\omega$ on $T^*Q$,
we can obtain two types of Hamilton-Jacobi
equation for the associated Marsden-Weinstein reduced Hamiltonian system,
which is given in Wang \cite{wa17}.
Thus, Theorem 5.3 and Theorem 5.4 can be regarded as an extension of two types of Hamilton-Jacobi
theorem for a Marsden-Weinstein reducible Hamiltonian system to that for the system with external force,
control and nonholonomic constrain.
\end{rema}

\begin{rema}
If $(T^\ast Q, \omega)$ is a connected symplectic manifold, and
$\mathbf{J}:T^\ast Q\rightarrow \mathfrak{g}^\ast$ is a
non-equivariant momentum map with a non-equivariance group
one-cocycle $\sigma: G\rightarrow \mathfrak{g}^\ast$,
which is defined by $\sigma(g):=\mathbf{J}(g\cdot
z)-\operatorname{Ad}^\ast_{g^{-1}}\mathbf{J}(z)$, where $g\in G$ and
$z\in T^\ast Q$. Then we know that $\sigma$ produces a new affine
action $\Theta: G\times \mathfrak{g}^\ast \rightarrow
\mathfrak{g}^\ast $ defined by
$\Theta(g,\mu):=\operatorname{Ad}^\ast_{g^{-1}}\mu + \sigma(g)$,
where $\mu \in \mathfrak{g}^\ast$, with respect to which the given
momentum map $\mathbf{J}$ is equivariant. Assume that $G$ acts
freely and properly on $T^\ast Q$, and $\tilde{G}_\mu$ denotes the
isotropy subgroup of $\mu \in \mathfrak{g}^\ast$ relative to this
affine action $\Theta$ and $\mu$ is a regular value of $\mathbf{J}$.
Then the quotient space $(T^\ast
Q)_\mu=\mathbf{J}^{-1}(\mu)/\tilde{G}_\mu$ is also a symplectic
manifold with symplectic form $\omega_\mu$ uniquely characterized by
$(5.1)$, see Ortega and Ratiu \cite{orra04}. In this case,
we can also define the $\mathbf{J}$-nonholonomic regular point reducible
RCH system $(T^*Q,G,\omega,\mathbf{J},\mathcal{D},H,F,W)$
with the associated distributional RCH system
$(\mathcal{K},\omega_{\mathcal {K}},H_{\mathcal{K}}, F_{\mathcal{K}}, u_{\mathcal{K}})$
and the $\mathbf{J}$-nonholonomic $R_p$-reduced distributional RCH system
$(\mathcal{K}_\mu,\omega_{\mathcal{K}_\mu},h_{\mathcal{K}_\mu},
f_{\mathcal {K}_\mu}, u_{\mathcal {K}_\mu})$,
and prove the Type I and Type II of
the Hamilton-Jacobi theorem for the $\mathbf{J}$-nonholonomic $R_p$-reduced
distributional RCH system
$(\mathcal{K}_\mu,\omega_{\mathcal{K}_\mu}, h_{\mathcal{K}_\mu},
f_{\mathcal {K}_\mu}, u_{\mathcal {K}_\mu})$ by using the above similar way,
in which the $\mathbf{J}$-nonholonomic $R_p$-reduced space
$(\mathcal{K}_\mu,\omega_{\mathcal{K}_\mu})$
is determined by the affine action and $\mathbf{J}$-nonholonomic regular point reduction.
\end{rema}

\subsection{Hamilton-Jacobi equations in the case compatible with regular orbit reduction}

It is worthy of noting that the regular orbit reduction of a Hamiltonian system
with symmetry and momentum map is an alternative
approach to symplectic reduction given by Marle \cite{ma76}
and Kazhdan, Kostant and Sternberg \cite{kakost78}, which is different from the
Marsden-Weinstein reduction, because the regular orbit reduced symplectic form
is different from the Marsden-Weinstein reduced symplectic form,
and it is also very important subject in the research of geometric mechanics\\

In this subsection, at first, we consider the regular orbit reducible RCH system,
which is given by Marsden et al \cite{mawazh10}.
Assume that the cotangent lifted
left action $\Phi^{T^\ast}:G\times T^\ast Q\rightarrow T^\ast Q$ is
symplectic, free and proper, and admits an
$\operatorname{Ad}^\ast$-equivariant momentum map $\mathbf{J}:T^\ast
Q\rightarrow \mathfrak{g}^\ast$. Let $\mu\in \mathfrak{g}^\ast$ be a
regular value of the momentum map $\mathbf{J}$ and
$\mathcal{O}_\mu=G\cdot \mu\subset \mathfrak{g}^\ast$ be the
$G$-orbit of the coadjoint $G$-action through the point $\mu$. Since
$G$ acts freely, properly and symplectically on $T^\ast Q$, then the
quotient space $(T^\ast Q)_{\mathcal{O}_\mu}=
\mathbf{J}^{-1}(\mathcal{O}_\mu)/G$ is a regular quotient symplectic
manifold with the $R_o$-reduced symplectic form $\omega_{\mathcal{O}_\mu}$
uniquely characterized by the relation
\begin{equation}i_{\mathcal{O}_\mu}^\ast \omega=\pi_{\mathcal{O}_{\mu}}^\ast
\omega_{\mathcal{O}
_\mu}+\mathbf{J}_{\mathcal{O}_\mu}^\ast\omega_{\mathcal{O}_\mu}^+,
\label{5.7}
\end{equation} where $\mathbf{J}_{\mathcal{O}_\mu}$ is
the restriction of the momentum map $\mathbf{J}$ to
$\mathbf{J}^{-1}(\mathcal{O}_\mu)$, that is,
$\mathbf{J}_{\mathcal{O}_\mu}=\mathbf{J}\cdot i_{\mathcal{O}_\mu}$
and $\omega_{\mathcal{O}_\mu}^+$ is the $+$-symplectic structure on
the orbit $\mathcal{O}_\mu$ given by
\begin{equation}
\omega_{\mathcal{O}_\mu}^
+(\nu)(\xi_{\mathfrak{g}^\ast}(\nu),\eta_{\mathfrak{g}^\ast}(\nu))
=<\nu,[\xi,\eta]>,\;\; \forall\;\nu\in\mathcal{O}_\mu, \;
\xi,\eta\in \mathfrak{g}.
\label{5.8} \end{equation}
The maps
$i_{\mathcal{O}_\mu}:\mathbf{J}^{-1}(\mathcal{O}_\mu)\rightarrow
T^\ast Q$ and
$\pi_{\mathcal{O}_\mu}:\mathbf{J}^{-1}(\mathcal{O}_\mu)\rightarrow
(T^\ast Q)_{\mathcal{O}_\mu}$ are natural injection and the
projection, respectively. The pair $((T^\ast
Q)_{\mathcal{O}_\mu},\omega_{\mathcal{O}_\mu})$ is called
the $R_o$-reduced symplectic space of $(T^\ast Q,\omega)$ at $\mu$.
In general case, we
maybe thought that the structure of the $R_o$-reduced symplectic
space $((T^\ast Q)_{\mathcal{O}_\mu},\omega_{\mathcal{O}_\mu})$ is
more complex than that of the $R_p$-reduced symplectic space
$((T^\ast Q)_\mu,\omega_\mu)$, but, from the regular reduction
diagram, see Ortega and Ratiu \cite{orra04},
we know that the $R_o$-reduced space $((T^\ast
Q)_{\mathcal{O}_\mu},\omega_{\mathcal{O}_\mu})$ is symplectic
diffeomorphic to the $R_p$-reduced space $((T^*Q)_\mu,
\omega_\mu)$, and hence is also symplectic diffeomorphic to a
symplectic fiber bundle.\\

Let $H:T^\ast Q\rightarrow \mathbb{R}$ be a $G$-invariant
Hamiltonian, the flow $F_t$ of the Hamiltonian vector field $X_H$
leaves the connected components of
$\mathbf{J}^{-1}(\mathcal{O}_\mu)$ invariant and commutes with the
$G$-action, so it induces a flow $f_t^{\mathcal{O}_\mu}$ on $(T^\ast
Q)_{\mathcal{O}_\mu}$, defined by $f_t^{\mathcal{O}_\mu}\cdot
\pi_{\mathcal{O}_\mu}=\pi_{\mathcal{O}_\mu} \cdot F_t\cdot
i_{\mathcal{O}_\mu}$, and the vector field $X_{h_{\mathcal{O}_\mu}}$
generated by the flow $f_t^{\mathcal{O}_\mu}$ on $((T^\ast
Q)_{\mathcal{O}_\mu},\omega_{\mathcal{O}_\mu})$ is Hamiltonian with
the associated $R_o$-reduced Hamiltonian function
$h_{\mathcal{O}_\mu}:(T^\ast Q)_{\mathcal{O}_\mu}\rightarrow
\mathbb{R}$ defined by $h_{\mathcal{O}_\mu}\cdot
\pi_{\mathcal{O}_\mu}= H\cdot i_{\mathcal{O}_\mu}$, and the
Hamiltonian vector fields $X_H$ and $X_{h_{\mathcal{O}_\mu}}$ are
$\pi_{\mathcal{O}_\mu}$-related.
Moreover, assume that the fiber-preserving map $F:T^\ast Q\rightarrow T^\ast
Q$ and the control subset $W$ of\; $T^\ast Q$ are both $G$-invariant.
In order to get the $R_o$-reduced RCH system, we also assume that
$F(\mathbf{J}^{-1}(\mathcal{O}_\mu))\subset
\mathbf{J}^{-1}(\mathcal{O}_\mu)$, and $W \cap
\mathbf{J}^{-1}(\mathcal{O}_\mu)\neq \emptyset $.
Thus, we can introduce a regular orbit
reducible RCH systems as follows, see Marsden et al \cite{mawazh10}
and Wang \cite{wa18}.

\begin{defi}
(Regular Orbit Reducible RCH System) A 6-tuple $(T^\ast Q, G,
\omega,H,F,W)$ with
the canonical symplectic form $\omega$ on $T^*Q$, where the Hamiltonian $H: T^\ast Q\rightarrow
\mathbb{R}$, the fiber-preserving map $F: T^\ast Q\rightarrow T^\ast
Q$ and the fiber submanifold $W$ of $T^\ast Q$ are all
$G$-invariant, is called a regular orbit reducible RCH system, if
there exists an orbit $\mathcal{O}_\mu, \; \mu\in\mathfrak{g}^\ast$,
where $\mu$ is a regular value of the momentum map $\mathbf{J}$,
such that the regular orbit reduced system, that is, the 5-tuple
$((T^\ast
Q)_{\mathcal{O}_\mu},\omega_{\mathcal{O}_\mu},
h_{\mathcal{O}_\mu},f_{\mathcal{O}_\mu},
W_{\mathcal{O}_\mu})$, where $(T^\ast
Q)_{\mathcal{O}_\mu}=\mathbf{J}^{-1}(\mathcal{O}_\mu)/G$,
$\pi_{\mathcal{O}_\mu}^\ast \omega_{\mathcal{O}_\mu}
=i_{\mathcal{O}_\mu}^\ast\omega
-\mathbf{J}_{\mathcal{O}_\mu}^\ast\omega_{\mathcal{O}_\mu}^+$,
$h_{\mathcal{O}_\mu}\cdot \pi_{\mathcal{O}_\mu} =H\cdot
i_{\mathcal{O}_\mu}$, $F(\mathbf{J}^{-1}(\mathcal{O}_\mu))\subset
\mathbf{J}^{-1}(\mathcal{O}_\mu)$, $f_{\mathcal{O}_\mu}\cdot
\pi_{\mathcal{O}_\mu}=\pi_{\mathcal{O}_\mu}\cdot F\cdot
i_{\mathcal{O}_\mu}$, and $W \cap
\mathbf{J}^{-1}(\mathcal{O}_\mu)\neq \emptyset $,
$W_{\mathcal{O}_\mu}=\pi_{\mathcal{O}_\mu}(W \cap
\mathbf{J}^{-1}(\mathcal{O}_\mu))$, is an RCH system,
which is simply written as $R_o$-reduced RCH system. Where $((T^\ast
Q)_{\mathcal{O}_\mu},\omega_{\mathcal{O}_\mu})$ is the $R_o$-reduced
space, the function $h_{\mathcal{O}_\mu}:(T^\ast
Q)_{\mathcal{O}_\mu}\rightarrow \mathbb{R}$ is called the $R_o$-reduced
Hamiltonian, the fiber-preserving map $f_{\mathcal{O}_\mu}:(T^\ast
Q)_{\mathcal{O}_\mu} \rightarrow (T^\ast Q)_{\mathcal{O}_\mu}$ is
called the $R_o$-reduced (external) force map, $W_{\mathcal{O}_\mu}$ is a
fiber submanifold of $(T^\ast Q)_{\mathcal{O}_\mu}$, and is called
the $R_o$-reduced control subset.
\end{defi}

In the following we consider that a nonholonomic RCH system
with symmetry and momentum map
is 8-tuple $(T^*Q,G, \omega,\mathbf{J},\mathcal{D},H,F,W)$,
where $\omega$ is the canonical symplectic form on $T^* Q$,
and the Lie group $G$, which may not be Abelian, acts smoothly by the left on $Q$,
its tangent lifted action on $TQ$ and its cotangent lifted action on $T^\ast Q$,
and $\mathcal{D}\subset TQ$ is a
$\mathcal{D}$-completely and $\mathcal{D}$-regularly nonholonomic
constraint of the system, and $\mathcal{D}$, $H, F$ and $W$ are all
$G$-invariant. Thus, the nonholonomic RCH system with symmetry and momentum map
is a regular orbit reducible RCH system with $G$-invariant
nonholonomic constraint $\mathcal{D}$.
Moreover, in the following we shall give
carefully a geometric formulation of the $\mathbf{J}$-nonholonomic
$R_o$-reduced distributional RCH system, by using momentum map and the
nonholonomic reduction compatible with regular orbit reduction.\\

Note that the Legendre transformation $\mathcal{F}L: TQ \rightarrow T^*Q $
is a fiber-preserving map, and $\mathcal{D}\subset TQ$ is $G$-invariant
for the tangent lifted left action $\Phi^{T}: G\times TQ\rightarrow TQ, $
then the constraint submanifold
$\mathcal{M}=\mathcal{F}L(\mathcal{D})\subset T^*Q$ is
$G$-invariant for the cotangent lifted left action $\Phi^{T^\ast}:
G\times T^\ast Q\rightarrow T^\ast Q$.
For the nonholonomic RCH system with symmetry
and momentum map  $(T^*Q,G, \omega,\mathbf{J},\mathcal{D},H,F,W)$,
in the same way, we define the distribution $\mathcal{F}$, which is the pre-image of the
nonholonomic constraints $\mathcal{D}$ for the map $T\pi_Q: TT^* Q
\rightarrow TQ$, that is, $\mathcal{F}=(T\pi_Q)^{-1}(\mathcal{D})$,
and the distribution $\mathcal{K}=\mathcal{F} \cap T\mathcal{M}$.
Moreover, we can also define the distributional two-form $\omega_\mathcal{K}$,
which is induced from the canonical symplectic form $\omega$ on $T^* Q$, that is,
$\omega_\mathcal{K}= \tau_{\mathcal{K}}\cdot \omega_{\mathcal{M}},$ and
$\omega_{\mathcal{M}}= i_{\mathcal{M}}^* \omega $.
If the admissibility condition $\mathrm{dim}\mathcal{M}=
\mathrm{rank}\mathcal{F}$ and the compatibility condition
$T\mathcal{M}\cap \mathcal{F}^\bot= \{0\}$ hold, then
$\omega_\mathcal{K}$ is non-degenerate as a
bilinear form on each fibre of $\mathcal{K}$, there exists a vector
field $X_\mathcal{K}$ on $\mathcal{M}$ which takes values in the
constraint distribution $\mathcal{K}$, such that for the function $H_\mathcal{K}$,
the following distributional Hamiltonian equation holds, that is,
\begin{align}
\mathbf{i}_{X_\mathcal{K}}\omega_\mathcal{K}
=\mathbf{d}H_\mathcal{K},
\label{5.9} \end{align}
where the function $H_{\mathcal{K}}$ satisfies
$\mathbf{d}H_{\mathcal{K}}= \tau_{\mathcal{K}}\cdot \mathbf{d}H_{\mathcal {M}}$,
and $H_\mathcal{M}= \tau_{\mathcal{M}}\cdot H$
is the restriction of $H$ to $\mathcal{M}$, and
from the equation (5.9), we have that
$X_{\mathcal{K}}=\tau_{\mathcal{K}}\cdot X_H $.\\

Since the nonholonomic RCH system with symmetry and momentum map
is a regular orbit reducible RCH system with $G$-invariant
nonholonomic constraint $\mathcal{D}$,
for a regular value $\mu\in\mathfrak{g}^\ast$ of the
momentum map $\mathbf{J}:T^\ast Q\rightarrow \mathfrak{g}^\ast$,
$\mathcal{O}_\mu=G\cdot \mu\subset \mathfrak{g}^\ast$ is the
$G$-orbit of the coadjoint $G$-action through the point $\mu$,
we assume that the constraint submanifold $\mathcal{M}$
is clean intersection with $\mathbf{J}^{-1}(\mathcal{O}_\mu)$, that is,
$\mathcal{M} \cap \mathbf{J}^{-1}(\mathcal{O}_\mu)\neq \emptyset$.
It follows that the quotient
space $\mathcal{M}_{\mathcal{O}_\mu} =(\mathcal{M}\cap \mathbf{J}^{-1}(\mathcal{O}_\mu))
/G \subset (T^\ast Q)_{\mathcal{O}_\mu}$ of the $G$-orbit in
$\mathcal{M}\cap \mathbf{J}^{-1}(\mathcal{O}_\mu)$, is a smooth manifold with
projection $\pi_{\mathcal{O}_\mu}: \mathcal{M}\cap \mathbf{J}^{-1}(\mathcal{O}_\mu)
\rightarrow \mathcal{M}_{\mathcal{O}_\mu}$ which is a surjective submersion.
Denote $i_{\mathcal{M}_{\mathcal{O}_\mu}}: \mathcal{M}_{\mathcal{O}_\mu}\rightarrow
(T^*Q)_{\mathcal{O}_\mu}, $ and $\omega_{\mathcal{M}_{\mathcal{O}_\mu}}= i_{\mathcal{M}_{\mathcal{O}_\mu}}^*
\omega_{\mathcal{O}_\mu} $, that is, the symplectic form
$\omega_{\mathcal{M}_{\mathcal{O}_\mu}}$ is induced from the $R_o$-reduced symplectic
form $\omega_{\mathcal{O}_\mu}$ on $(T^* Q)_{\mathcal{O}_\mu}$ given in (5.7),
where $i_{\mathcal{M}_{\mathcal{O}_\mu}}^*:
T^*(T^*Q)_{\mathcal{O}_\mu} \rightarrow T^*\mathcal{M}_{\mathcal{O}_\mu}.$
Moreover, the distribution $\mathcal{F}$ is pushed down to a distribution
$\mathcal{F}_{\mathcal{O}_\mu}= T\pi_{\mathcal{O}_\mu}\cdot \mathcal{F}$
on $(T^\ast Q)_{\mathcal{O}_\mu}$,
and we define $\mathcal{K}_{\mathcal{O}_\mu}=\mathcal{F}_{\mathcal{O}_\mu} \cap
T\mathcal{M}_{\mathcal{O}_\mu}$. Assume that $\omega_{\mathcal{K}_{\mathcal{O}_\mu}}=
\tau_{\mathcal{K}_{\mathcal{O}_\mu}}\cdot \omega_{\mathcal{M}_{\mathcal{O}_\mu}}$ is the
restriction of the symplectic form $\omega_{\mathcal{M}_{\mathcal{O}_\mu}}$ on
$T^*\mathcal{M}_{\mathcal{O}_\mu}$ fibrewise to
the distribution $\mathcal{K}_{\mathcal{O}_\mu}$,
where $\tau_{\mathcal{K}_{\mathcal{O}_\mu}}$ is the restriction map to distribution
$\mathcal{K}_{\mathcal{O}_\mu}$. The $\omega_{\mathcal{K}_{\mathcal{O}_\mu}}$ is not
a true two-form on a manifold, which is called
as a $\mathbf{J}$-nonholonomic $R_o$-reduced
distributional two-form to avoid any confusion.\\

From the above construction we know that,
if the admissibility condition $\mathrm{dim}\mathcal{M}_{\mathcal{O}_\mu}=
\mathrm{rank}\mathcal{F}_{\mathcal{O}_\mu}$ and the compatibility condition
$T\mathcal{M}_{\mathcal{O}_\mu}\cap \mathcal{F}_{\mathcal{O}_\mu}^\bot= \{0\}$ hold, where
$\mathcal{F}_{\mathcal{O}_\mu}^\bot$ denotes the symplectic orthogonal of
$\mathcal{F}_{\mathcal{O}_\mu}$ with respect to the $R_o$-reduced symplectic form
$\omega_{\mathcal{O}_\mu}$, then
$\omega_{\mathcal{K}_{\mathcal{O}_\mu}}$ is
non-degenerate as a bilinear form on each fibre of
$\mathcal{K}_{\mathcal{O}_\mu}$, and hence there exists a vector field
$X_{\mathcal{K}_{\mathcal{O}_\mu}}$
on $\mathcal{M}_{\mathcal{O}_\mu}$, which takes values in the constraint
distribution $\mathcal{K}_{\mathcal{O}_\mu}$,
such that for the function $h_{\mathcal{K}_{\mathcal{O}_\mu}}$,
the $\mathbf{J}$-nonholonomic $R_o$-reduced distributional
Hamiltonian equation holds, that is,
\begin{align}
\mathbf{i}_{X_{\mathcal{K}_{\mathcal{O}_\mu}}}\omega_{\mathcal{K}_{\mathcal{O}_\mu}}
=\mathbf{d}h_{\mathcal{K}_{\mathcal{O}_\mu}},
\label{5.10} \end{align}
where $\mathbf{d}h_{\mathcal{K}_{\mathcal{O}_\mu}}$ is the restriction
of $\mathbf{d}h_{\mathcal{M}_{\mathcal{O}_\mu}}$ to $\mathcal{K}_{\mathcal{O}_\mu}$, and
the function $h_{\mathcal{K}_{\mathcal{O}_\mu}}$ satisfies
$\mathbf{d}h_{\mathcal{K}_{\mathcal{O}_\mu}}= \tau_{\mathcal{K}_{\mathcal{O}_\mu}}\cdot \mathbf{d}h_{\mathcal{M}_{\mathcal{O}_\mu}} $,
and $h_{\mathcal{M}_{\mathcal{O}_\mu}}
= \tau_{\mathcal{M}_{\mathcal{O}_\mu}}\cdot h_{\mathcal{O}_\mu}$ is the
restriction of $h_{\mathcal{O}_\mu}$ to $\mathcal{M}_{\mathcal{O}_\mu}$,
and $h_{\mathcal{O}_\mu}$ is the $R_o$-reduced
Hamiltonian function $h_{\mathcal{O}_\mu}: (T^* Q)_{\mathcal{O}_\mu} \rightarrow \mathbb{R}$ defined
by $h_{\mathcal{O}_\mu}\cdot \pi_{\mathcal{O}_\mu}= H\cdot i_{\mathcal{O}_\mu}$.
In addition, from the distributional Hamiltonian equation (5.9)
$\mathbf{i}_{X_\mathcal{K}}\omega_\mathcal{K}=\mathbf{d}H_\mathcal
{K},$ we have that $X_{\mathcal{K}}=\tau_{\mathcal{K}}\cdot X_H, $
and from the $\mathbf{J}$-nonholonomic $R_o$-reduced distributional Hamiltonian equation (5.10)
$\mathbf{i}_{X_{\mathcal{K_{\mathcal{O}_\mu}}}}\omega_{\mathcal{K_{\mathcal{O}_\mu}}}
=\mathbf{d}h_{\mathcal{K_{\mathcal{O}_\mu}}}$, we have that
$X_{\mathcal{K_{\mathcal{O}_\mu}}}
=\tau_{\mathcal{K_{\mathcal{O}_\mu}}}\cdot X_{h_{\mathcal{K_{\mathcal{O}_\mu}}}},$
where $ X_{h_{\mathcal{K_{\mathcal{O}_\mu}}}}$ is the Hamiltonian vector field of
the function $h_{\mathcal{K_{\mathcal{O}_\mu}}},$
and the vector fields $X_{\mathcal{K}}$
and $X_{\mathcal{K_{\mathcal{O}_\mu}}}$ are $\pi_{\mathcal{O}_\mu}$-related,
that is, $X_{\mathcal{K_{\mathcal{O}_\mu}}}\cdot \pi_{\mathcal{O}_\mu}
=T\pi_{\mathcal{O}_\mu}\cdot X_{\mathcal{K}}.$ \\

Moreover, if considering the external force $F$ and control subset $W$,
and we define the vector field $F_\mathcal{K}
=\tau_{\mathcal{K}}\cdot \textnormal{vlift}(F_{\mathcal{M}})X_H,$
and for a control law $u\in W$, the vector field
$u_\mathcal{K}= \tau_{\mathcal{K}}\cdot  \textnormal{vlift}(u_{\mathcal{M}})X_H,$
where $F_\mathcal{M}= \tau_{\mathcal{M}}\cdot F$ and
$u_\mathcal{M}= \tau_{\mathcal{M}}\cdot u$ are the restrictions of
$F$ and $u$ to $\mathcal{M}$, that is, $F_\mathcal{K}$ and $u_\mathcal{K}$
are the restrictions of the changes of Hamiltonian vector field $X_H$
under the actions of $F_\mathcal{M}$ and $u_\mathcal{M}$ to $\mathcal{K}$,
then the 5-tuple $(\mathcal{K},\omega_{\mathcal{K}},
H_\mathcal{K}, F_\mathcal{K}, u_\mathcal{K})$
is a distributional RCH system corresponding to the nonholonomic RCH system with symmetry
and momentum map $(T^*Q,G, \omega,\mathbf{J},\mathcal{D},H,F,u)$,
and the dynamical vector field of the distributional RCH system
can be expressed by
\begin{align}
\tilde{X}=X_{(\mathcal{K},\omega_{\mathcal{K}},
H_{\mathcal{K}}, F_{\mathcal{K}}, u_{\mathcal{K}})}
=X_\mathcal {K}+ F_{\mathcal{K}}+u_{\mathcal{K}},
\label{5.11} \end{align}
which is the synthetic
of the nonholonomic dynamical vector field $X_{\mathcal{K}}$ and
the vector fields $F_{\mathcal{K}}$ and $u_{\mathcal{K}}$.
Assume that the vector fields $F_\mathcal{K}$ and $u_\mathcal{K}$
on $\mathcal{M}$ are pushed down to the vector fields
$f_{\mathcal{M}_{\mathcal{O}_\mu}}=T\pi_{\mathcal{O}_\mu} \cdot F_\mathcal{K}$ and
$u_{\mathcal{M}_{\mathcal{O}_\mu}}=T\pi_{\mathcal{O}_\mu} \cdot u_\mathcal{K}$ on $\mathcal{M}_{\mathcal{O}_\mu}$.
Then we define that $f_{\mathcal{K}_{\mathcal{O}_\mu}}
=T\pi_{\mathcal{K}_{\mathcal{O}_\mu}}\cdot f_{\mathcal{M}_{\mathcal{O}_\mu}}$ and
$u_{\mathcal{K}_{\mathcal{O}_\mu}}
=T\pi_{\mathcal{K}_{\mathcal{O}_\mu}}\cdot u_{\mathcal{M}_{\mathcal{O}_\mu}},$
that is, $f_{\mathcal{K}_{\mathcal{O}_\mu}}$ and
$u_{\mathcal{K}_{\mathcal{O}_\mu}}$ are the restrictions of
$f_{\mathcal{M}_{\mathcal{O}_\mu}}$ and $u_{\mathcal{M}_{\mathcal{O}_\mu}}$ to $\mathcal{K}_{\mathcal{O}_\mu}$.
In consequence, the 5-tuple $(\mathcal{K}_{\mathcal{O}_\mu},\omega_{\mathcal{K}_{\mathcal{O}_\mu}},
h_{\mathcal{K}_{\mathcal{O}_\mu}}, f_{\mathcal{K}_{\mathcal{O}_\mu}},
u_{\mathcal{K}_{\mathcal{O}_\mu}})$
is a $\mathbf{J}$-nonholonomic $R_o$-reduced distributional RCH system of the nonholonomic
RCH system with symmetry and momentum map
$(T^*Q,G,\omega,\mathbf{J},\mathcal{D},H,F,W)$,
as well as with a control law $u \in W$.
Thus, the geometrical formulation
of the $\mathbf{J}$-nonholonomic $R_o$-reduced distributional RCH
system may be summarized as follows.

\begin{defi} ($\mathbf{J}$-Nonholonomic $R_o$-reduced Distributional RCH System)
Assume that the 8-tuple $(T^*Q,G,\omega,\mathbf{J},\mathcal{D},H, F, W)$
is a nonholonomic RCH system with symmetry and momentum map,
where $\omega$ is the canonical symplectic form on $T^* Q$,
and $\mathcal{D}\subset TQ$ is a $\mathcal{D}$-completely and
$\mathcal{D}$-regularly nonholonomic constraint of the system, and
$\mathcal{D}$, $H, F$ and $W$ are all $G$-invariant. For a regular value
$\mu\in\mathfrak{g}^\ast$ of the momentum map $\mathbf{J}$,
$\mathcal{O}_\mu=G\cdot \mu\subset \mathfrak{g}^\ast$ is the
$G$-orbit of the coadjoint $G$-action through the point $\mu$,
assume that there exists a
$\mathbf{J}$-nonholonomic $R_o$-reduced distribution
$\mathcal{K}_{\mathcal{O}_\mu}$, an associated non-degenerate
and $\mathbf{J}$-nonholonomic $R_o$-reduced distributional two-form
$\omega_{\mathcal{K}_{\mathcal{O}_\mu}}$
and a vector field $X_{\mathcal {K}_{\mathcal{O}_\mu}}$
on the $\mathbf{J}$-nonholonomic $R_o$-reduced constraint submanifold
$\mathcal{M}_{\mathcal{O}_\mu}=(\mathcal{M}\cap \mathbf{J}^{-1}(\mathcal{O}_\mu)) /G, $
where $\mathcal{M}=\mathcal{F}L(\mathcal{D}),$ and $\mathcal{M}\cap
\mathbf{J}^{-1}({\mathcal{O}_\mu})\neq \emptyset, $ such that the
$\mathbf{J}$-nonholonomic $R_o$-reduced distributional Hamiltonian
equation (5.10) holds, that is, $\mathbf{i}_{X_{\mathcal{K}_{\mathcal{O}_\mu}}}\omega_{\mathcal{K}_{\mathcal{O}_\mu}} =
\mathbf{d}h_{\mathcal{K}_{\mathcal{O}_\mu}}$, where
$\mathbf{d}h_{\mathcal{K}_{\mathcal{O}_\mu}}$ is the restriction of
$\mathbf{d}h_{\mathcal{M}_{\mathcal{O}_\mu}}$ to $\mathcal{K}_{\mathcal{O}_\mu}$,
and the function $h_{\mathcal{K}_{\mathcal{O}_\mu}}$,
and the vector fields $f_{\mathcal{K}_{\mathcal{O}_\mu}}$ and
$u_{\mathcal{K}_{\mathcal{O}_\mu}}$ are defined above.
Then the 5-tuple
$(\mathcal{K}_{\mathcal{O}_\mu},\omega_{\mathcal {K}_{\mathcal{O}_\mu}},
h_{\mathcal {K}_{\mathcal{O}_\mu}}, f_{\mathcal {K}_{\mathcal{O}_\mu}},
u_{\mathcal {K}_{\mathcal{O}_\mu}})$ is called a
$\mathbf{J}$-nonholonomic $R_o$-reduced distributional RCH system
of the nonholonomic RCH system with symmetry and momentum map
$(T^*Q,G,\omega,\mathbf{J},\mathcal{D},H, F, u)$
with a control law $u \in W$, and $X_{\mathcal
{K}_{\mathcal{O}_\mu}}$ is called the $\mathbf{J}$-nonholonomic
$R_o$-reduced dynamical vector field.
Denote by
\begin{align}
\hat{X}_{\mathcal{O}_\mu}=X_{(\mathcal{K}_{\mathcal{O}_\mu},
\omega_{\mathcal{K}_{\mathcal{O}_\mu}}, h_{\mathcal {K}_{\mathcal{O}_\mu}},
f_{\mathcal {K}_{\mathcal{O}_\mu}}, u_{\mathcal {K}_{\mathcal{O}_\mu}})}
=X_{\mathcal{K}_{\mathcal{O}_\mu}}+ f_{\mathcal{K}_{\mathcal{O}_\mu}}+u_{\mathcal{K}_{\mathcal{O}_\mu}}
\label{5.12} \end{align}
is the dynamical vector field of the
$\mathbf{J}$-nonholonomic $R_o$-reduced distributional RCH system
$(\mathcal{K}_{\mathcal{O}_\mu},\omega_{\mathcal{K}_{\mathcal{O}_\mu}},\\
h_{\mathcal {K}_{\mathcal{O}_\mu}},
f_{\mathcal {K}_{\mathcal{O}_\mu}}, u_{\mathcal {K}_{\mathcal{O}_\mu}})$,
which is the synthetic of the $\mathbf{J}$-nonholonomic
$R_o$-reduced dynamical vector field $X_{\mathcal {K}_{\mathcal{O}_\mu}}$ and
the vector fields $F_{\mathcal {K}_{\mathcal{O}_\mu}}$
and $u_{\mathcal {K}_{\mathcal{O}_\mu}}$.
Under the above circumstances, we refer to
$(T^*Q,G,\omega,\mathbf{J},\mathcal{D},H, F, u)$ as
a $\mathbf{J}$-nonholonomic regular orbit reducible RCH system
with the associated distributional RCH system
$(\mathcal{K},\omega_{\mathcal {K}},H_{\mathcal{K}}, F_{\mathcal{K}}, u_{\mathcal{K}})$
and the $\mathbf{J}$-nonholonomic $R_o$-reduced
distributional RCH system
$(\mathcal{K}_{\mathcal{O}_\mu},\omega_{\mathcal{K}_{\mathcal{O}_\mu}},
h_{\mathcal{K}_{\mathcal{O}_\mu}},
f_{\mathcal {K}_{\mathcal{O}_\mu}}, u_{\mathcal {K}_{\mathcal{O}_\mu}})$.
\end{defi}

Since the non-degenerate and $\mathbf{J}$-nonholonomic $R_o$-reduced
distributional two-form $\omega_{\mathcal{K}_{\mathcal{O}_\mu}}$ is not a "true two-form"
on a manifold, and it is not symplectic, and hence
the $\mathbf{J}$-nonholonomic $R_o$-reduced distributional RCH system
$(\mathcal{K}_{\mathcal{O}_\mu},\omega_{\mathcal {K}_{\mathcal{O}_\mu}},
h_{\mathcal{K}_{\mathcal{O}_\mu}},
f_{\mathcal {K}_{\mathcal{O}_\mu}}, u_{\mathcal {K}_{\mathcal{O}_\mu}})$
is not a Hamiltonian system, and has no yet generating function,
and hence we can not describe the Hamilton-Jacobi equation for a
$\mathbf{J}$-nonholonomic $R_o$-reduced
distributional RCH system just like as in Theorem 1.1.
But, for a given $\mathbf{J}$-nonholonomic regular orbit reducible RCH system $(T^*Q,G,\omega,\mathbf{J},\mathcal{D},H, F, u)$ with the associated distributional RCH system
$(\mathcal{K},\omega_{\mathcal {K}},H_{\mathcal{K}}, F_{\mathcal{K}}, u_{\mathcal{K}})$
and the $\mathbf{J}$-nonholonomic $R_o$-reduced distributional RCH system
$(\mathcal{K}_{\mathcal{O}_\mu},\omega_{\mathcal {K}_{\mathcal{O}_\mu}},
h_{\mathcal{K}_{\mathcal{O}_\mu}}, f_{\mathcal {K}_{\mathcal{O}_\mu}},
u_{\mathcal {K}_{\mathcal{O}_\mu}})$,
by using Lemma 3.3, we can derive precisely
the geometric constraint conditions of the $\mathbf{J}$-nonholonomic $R_o$-reduced distributional two-form
$\omega_{\mathcal{K}_{\mathcal{O}_\mu}}$ for the nonholonomic reducible dynamical vector field,
that is, the two types of Hamilton-Jacobi equation for the
$\mathbf{J}$-nonholonomic $R_o$-reduced distributional RCH system
$(\mathcal{K}_{\mathcal{O}_\mu},\omega_{\mathcal
{K}_{\mathcal{O}_\mu}},h_{\mathcal{K}_{\mathcal{O}_\mu}},
f_{\mathcal {K}_{\mathcal{O}_\mu}}, u_{\mathcal {K}_{\mathcal{O}_\mu}})$.
At first, using the fact that the one-form $\gamma: Q
\rightarrow T^*Q $ is closed on $\mathcal{D}$ with respect to
$T\pi_Q: TT^* Q \rightarrow TQ, $
and $\textmd{Im}(\gamma)\subset \mathcal{M} \cap
\mathbf{J}^{-1}(\mathcal{O}_\mu), $ and it is $G$-invariant,
as well as $ \textmd{Im}(T\bar{\gamma}_{\mathcal{O}_\mu})\subset \mathcal{K}_{\mathcal{O}_\mu},$
we can prove the Type I of
Hamilton-Jacobi theorem for the $\mathbf{J}$-nonholonomic $R_o$-reduced distributional
RCH system. For convenience, the maps involved in the
following theorem and its proof are shown in Diagram-7.
\begin{center}
\hskip 0cm \xymatrix{ \mathbf{J}^{-1}(\mathcal{O}_\mu) \ar[r]^{i_{\mathcal{O}_\mu}}
& T^* Q  \ar[r]^{\pi_Q}
& Q \ar[d]_{\tilde{X}^\gamma} \ar[r]^{\gamma}
& T^*Q \ar[d]_{\tilde{X}} \ar[r]^{\pi_{\mathcal{O}_\mu}} & (T^* Q)_{\mathcal{O}_\mu} \ar[d]_{\hat{X}_{\mathcal{O}_\mu}}
& \mathcal{M}_{\mathcal{O}_\mu}  \ar[l]_{i_{\mathcal{M}_{\mathcal{O}_\mu}}} \ar[d]_{X_{\mathcal{K}_{\mathcal{O}_\mu}}}\\
& T(T^*Q)  & TQ \ar[l]_{T\gamma} & T(T^*Q) \ar[l]^{T\pi_Q} \ar[r]_{T\pi_{\mathcal{O}_\mu}}
& T(T^* Q)_{\mathcal{O}_\mu} \ar[r]^{\tau_{\mathcal{K}_{\mathcal{O}_\mu}}} & \mathcal{K}_{\mathcal{O}_\mu} }
\end{center}
$$\mbox{Diagram-7}$$

\begin{theo} (Type I of Hamilton-Jacobi Theorem for a $\mathbf{J}$-Nonholonomic
$R_o$-reduced Distributional RCH System) For a given
$\mathbf{J}$-nonholonomic regular orbit reducible RCH system
$(T^*Q,G,\omega,\mathbf{J},\\ \mathcal{D},H, F, u)$ with the associated distributional RCH system
$(\mathcal{K},\omega_{\mathcal {K}},H_{\mathcal{K}}, F_{\mathcal{K}}, u_{\mathcal{K}})$
and the $\mathbf{J}$-nonholonomic $R_o$-reduced distributional RCH system
$(\mathcal{K}_{\mathcal{O}_\mu},\omega_{\mathcal{K}_{\mathcal{O}_\mu}},
h_{\mathcal{K}_{\mathcal{O}_\mu}},
f_{\mathcal {K}_{\mathcal{O}_\mu}},
u_{\mathcal {K}_{\mathcal{O}_\mu}})$, assume that $\gamma:
Q \rightarrow T^*Q$ is an one-form on $Q$, and
$\tilde{X}^\gamma = T\pi_{Q}\cdot \tilde{X}\cdot \gamma$,
where $\tilde{X}=X_{(\mathcal{K},\omega_{\mathcal{K}},
H_{\mathcal{K}}, F_{\mathcal{K}}, u_{\mathcal{K}})}
=X_\mathcal {K}+ F_{\mathcal{K}}+u_{\mathcal{K}}$
is the dynamical vector field of the distributional RCH system
$(\mathcal{K},\omega_{\mathcal{K}},
H_{\mathcal{K}}, F_{\mathcal{K}}, u_{\mathcal{K}})$
corresponding to the $\mathbf{J}$-nonholonomic regular orbit reducible RCH
system with symmetry and momentum map $(T^*Q,G,\omega,\mathbf{J},\mathcal{D},H, F,u)$.
Moreover, assume that $\mu\in\mathfrak{g}^\ast$ is a regular value of the momentum
map $\mathbf{J}$, and $\textmd{Im}(\gamma)\subset \mathcal{M} \cap
\mathbf{J}^{-1}(\mu), $ and it is $G$-invariant, and
$\bar{\gamma}_{\mathcal{O}_\mu}
=\pi_{\mathcal{O}_\mu}(\gamma): Q \rightarrow \mathcal{M}_{\mathcal{O}_\mu} $,
and $ \textmd{Im}(T\bar{\gamma}_{\mathcal{O}_\mu})\subset \mathcal{K}_{\mathcal{O}_\mu}. $
If the one-form $\gamma: Q \rightarrow T^*Q $ is closed on $\mathcal{D}$ with respect to
$T\pi_Q: TT^* Q \rightarrow TQ, $ then
$\bar{\gamma}_{\mathcal{O}_\mu}$ is a solution of the
equation $T\bar{\gamma}_{\mathcal{O}_\mu}\cdot
\tilde{X}^\gamma= X_{\mathcal{K}_{\mathcal{O}_\mu}}\cdot \bar{\gamma}_{\mathcal{O}_\mu}. $
Here $X_{\mathcal{K}_{\mathcal{O}_\mu}}$ is the $\mathbf{J}$-nonholonomic $R_o$-reduced
dynamical vector field. The equation
$T\bar{\gamma}_{\mathcal{O}_\mu}\cdot \tilde{X}^\gamma
= X_{\mathcal{K}_{\mathcal{O}_\mu}}\cdot
\bar{\gamma}_{\mathcal{O}_\mu},$ is called the Type I of Hamilton-Jacobi equation for the
$\mathbf{J}$-nonholonomic $R_o$-reduced distributional RCH system
$(\mathcal{K}_{\mathcal{O}_\mu},\omega_{\mathcal{K}_{\mathcal{O}_\mu}},
h_{\mathcal{K}_{\mathcal{O}_\mu}},
f_{\mathcal {K}_{\mathcal{O}_\mu}}, u_{\mathcal {K}_{\mathcal{O}_\mu}})$.
\end{theo}

\noindent{\bf Proof: } At first, for the dynamical vector field of the distributional RCH system
$(\mathcal{K},\omega_{\mathcal {K}}, H_{\mathcal{K}}, F_{\mathcal{K}}, u_{\mathcal{K}})$,
$\tilde{X}=X_{(\mathcal{K},\omega_{\mathcal{K}},
H_{\mathcal{K}}, F_{\mathcal{K}}, u_{\mathcal{K}})}
=X_\mathcal {K}+ F_{\mathcal{K}}+u_{\mathcal{K}}$,
and $F_{\mathcal{K}}=\tau_{\mathcal{K}}\cdot \textnormal{vlift}(F_{\mathcal{M}})X_H$,
and $u_{\mathcal{K}}=\tau_{\mathcal{K}}\cdot \textnormal{vlift}(u_{\mathcal{M}})X_H$,
note that $T\pi_{Q}\cdot \textnormal{vlift}(F_{\mathcal{M}})X_H=T\pi_{Q}\cdot \textnormal{vlift}(u_{\mathcal{M}})X_H=0, $
then we have that $T\pi_{Q}\cdot F_{\mathcal{K}}=T\pi_{Q}\cdot u_{\mathcal{K}}=0,$
and hence $T\pi_{Q}\cdot \tilde{X}\cdot \gamma=T\pi_{Q}\cdot X_{\mathcal{K}}\cdot \gamma. $
Moreover, from Theorem 3.4, we know that
$\gamma$ is a solution of the Hamilton-Jacobi equation
$T\gamma\cdot \tilde{X}^\gamma= X_{\mathcal{K}}\cdot \gamma .$ Next, we note that
the $R_o$-reduced symplectic space
$(T^\ast Q)_{\mathcal{O}_\mu}= \mathbf{J}^{-1}(\mathcal{O}_\mu)/G
\cong \mathbf{J}^{-1}(\mu)/G \times \mathcal{O}_\mu, $ with the
symplectic form $\omega_{\mathcal{O}_\mu}$ uniquely characterized by
the relation $i_{\mathcal{O}_\mu}^\ast
\omega=\pi_{\mathcal{O}_{\mu}}^\ast \omega_{\mathcal{O}
_\mu}+\mathbf{J}_{\mathcal{O}_\mu}^\ast\omega_{\mathcal{O}_\mu}^+. $
Since
$\textmd{Im}(\gamma)\subset \mathcal{M} \cap \mathbf{J}^{-1}(\mu), $
and it is $G$-invariant,
in this case for any $V\in TQ, $ and $w\in TT^*Q, $
we have that
$\mathbf{J}_{\mathcal{O}_\mu}^\ast\omega_{\mathcal{O}_\mu}^+(T\gamma
\cdot V, \; w)=0, $ and hence
$\pi_{\mathcal{O}_\mu}^*\omega_{\mathcal{O}_\mu}=
i_{\mathcal{O}_\mu}^*\omega= \omega, $ along $\textmd{Im}(\gamma)$.
On the other hand, because
$\textmd{Im}(T\bar{\gamma}_{\mathcal{O}_\mu})\subset \mathcal{K}_{\mathcal{O}_\mu}, $
then
$\omega_{\mathcal{K}_{\mathcal{O}_\mu}}\cdot
\tau_{\mathcal{K}_{\mathcal{O}_\mu}}=\tau_{\mathcal{K}_{\mathcal{O}_\mu}}\cdot
\omega_{\mathcal{M}_{\mathcal{O}_\mu}}= \tau_{\mathcal{K}_{\mathcal{O}_\mu}}\cdot
i_{\mathcal{M}_{\mathcal{O}_\mu}}^* \cdot \omega_{\mathcal{O}_\mu}, $ along
$\textmd{Im}(T\bar{\gamma}_{\mathcal{O}_\mu})$.
From the distributional Hamiltonian equation (5.9),
we have that $X_{\mathcal{K}}= \tau_{\mathcal{K}}\cdot X_H,$
and $\tau_{\mathcal{K}}\cdot X_{H}\cdot \gamma
= X_{\mathcal{K}}\cdot \gamma \in \mathcal{K}$.
Because the vector fields $X_{\mathcal{K}}$
and $X_{\mathcal{K}_{\mathcal{O}_\mu}}$ are $\pi_{\mathcal{O}_\mu}$-related,
$T\pi_{\mathcal{O}_\mu}(X_{\mathcal{K}})
=X_{\mathcal{K}_{\mathcal{O}_\mu}}\cdot \pi_{\mathcal{O}_\mu}$,
and hence $\tau_{\mathcal{K}_{\mathcal{O}_\mu}}\cdot T\pi_{\mathcal{O}_\mu}(X_{\mathcal{K}}\cdot \gamma)
=\tau_{\mathcal{K}_{\mathcal{O}_\mu}}\cdot (T\pi_{\mathcal{O}_\mu}(X_{\mathcal{K}}))\cdot (\gamma)
= \tau_{\mathcal{K}_{\mathcal{O}_\mu}}\cdot (X_{\mathcal{K}_{\mathcal{O}_\mu}}\cdot \pi_{\mathcal{O}_\mu})\cdot (\gamma)
= \tau_{\mathcal{K}_{\mathcal{O}_\mu}}\cdot X_{\mathcal{K}_{\mathcal{O}_\mu}}\cdot \pi_{\mathcal{O}_\mu}(\gamma)
= X_{\mathcal{K}_{\mathcal{O}_\mu}}\cdot \bar{\gamma}.$
Thus, using the $\mathbf{J}$-nonholonomic $R_o$-reduced distributional two-form
$\omega_{\mathcal{K}_{\mathcal{O}_\mu}}$, from Lemma 3.3(ii) and (iii), if we take that $v=
X_{\mathcal{K}}\cdot \gamma \in \mathcal{K} (\subset \mathcal{F}), $ and
for any $w \in \mathcal{F}, \; T\lambda(w)\neq 0, $ and
$\tau_{\mathcal{K}_{\mathcal{O}_\mu}}\cdot T\pi_{\mathcal{O}_\mu} \cdot w \neq 0, $ then we have that
\begin{align*}
& \omega_{\mathcal{K}_{\mathcal{O}_\mu}}(T\bar{\gamma}_{\mathcal{O}_\mu} \cdot \tilde{X}^\gamma, \;
\tau_{\mathcal{K}_{\mathcal{O}_\mu}}\cdot T\pi_{\mathcal{O}_\mu} \cdot w)=
\omega_{\mathcal{K}_{\mathcal{O}_\mu}}(\tau_{\mathcal{K}_{\mathcal{O}_\mu}}\cdot
T\bar{\gamma}_{\mathcal{O}_\mu} \cdot \tilde{X}^\gamma, \; \tau_{\mathcal{K}_{\mathcal{O}_\mu}}\cdot T\pi_{\mathcal{O}_\mu} \cdot w)\\
& = \tau_{\mathcal{K}_{\mathcal{O}_\mu}}\cdot \omega_{\mathcal{M}_{\mathcal{O}_\mu}}(T(\pi_{\mathcal{O}_\mu}
\cdot\gamma )\cdot \tilde{X}^\gamma, \; T\pi_{\mathcal{O}_\mu} \cdot w ) =
\tau_{\mathcal{K}_{\mathcal{O}_\mu}}\cdot i_{\mathcal{M}_{\mathcal{O}_\mu}}^* \cdot
\omega_{\mathcal{O}_\mu}(T\pi_{\mathcal{O}_\mu} \cdot T\gamma \cdot \tilde{X}^\gamma,
\; T\pi_{\mathcal{O}_\mu} \cdot w )\\
& = \tau_{\mathcal{K}_{\mathcal{O}_\mu}}\cdot i_{\mathcal{M}_{\mathcal{O}_\mu}}^*
\cdot \pi_{\mathcal{O}_\mu}^*\omega_{\mathcal{O}_\mu}(T\gamma \cdot T\pi_Q \cdot \tilde{X} \cdot
\gamma, \; w) = \tau_{\mathcal{K}_{\mathcal{O}_\mu}}\cdot i_{\mathcal{M}_{\mathcal{O}_\mu}}^*
\cdot \omega(T(\gamma \cdot \pi_Q) \cdot X_{\mathcal{K}}\cdot \gamma, \; w)\\
& =\tau_{\mathcal{K}_{\mathcal{O}_\mu}}\cdot i_{\mathcal{M}_{\mathcal{O}_\mu}}^* \cdot
(\omega(X_{\mathcal{K}}\cdot \gamma, \; w-T(\gamma
\cdot \pi_Q)\cdot w) -\mathbf{d}\gamma(T\pi_{Q}(X_{\mathcal{K}}\cdot \gamma), \; T\pi_{Q}(w)))\\
& = \tau_{\mathcal{K}_{\mathcal{O}_\mu}}\cdot i_{\mathcal{M}_{\mathcal{O}_\mu}}^*
\cdot\pi_{\mathcal{O}_\mu}^*\omega_{\mathcal{O}_\mu}(X_{\mathcal{K}}\cdot \gamma, \; w) -
\tau_{\mathcal{K}_{\mathcal{O}_\mu}}\cdot i_{\mathcal{M}_{\mathcal{O}_\mu}}^*
\cdot\pi_{\mathcal{O}_\mu}^*\omega_{\mathcal{O}_\mu}(X_{\mathcal{K}}\cdot
\gamma, \; T(\gamma \cdot \pi_Q) \cdot w)\\
& \;\;\;\;\;\; -\tau_{\mathcal{K}_{\mathcal{O}_\mu}}\cdot i_{\mathcal{M}_{\mathcal{O}_\mu}}^* \cdot\mathbf{d}\gamma(T\pi_{Q}(X_{\mathcal{K}}\cdot \gamma), \; T\pi_{Q}(w))\\
& = \tau_{\mathcal{K}_{\mathcal{O}_\mu}}\cdot i_{\mathcal{M}_{\mathcal{O}_\mu}}^*
\cdot\omega_{\mathcal{O}_\mu}(T\pi_{\mathcal{O}_\mu}(X_{\mathcal{K}}\cdot \gamma),
\; T\pi_{\mathcal{O}_\mu} \cdot w)\\
& \;\;\;\;\;\; -\tau_{\mathcal{K}_{\mathcal{O}_\mu}}\cdot i_{\mathcal{M}_{\mathcal{O}_\mu}}^* \cdot
\omega_{\mathcal{O}_\mu}(T\pi_{\mathcal{O}_\mu}\cdot(X_{\mathcal{K}}\cdot \gamma), \; T(\pi_{\mathcal{O}_\mu}
\cdot\gamma) \cdot T\pi_{Q}(w))\\
& \;\;\;\;\;\; -\tau_{\mathcal{K}_{\mathcal{O}_\mu}}\cdot
i_{\mathcal{M}_{\mathcal{O}_\mu}}^* \cdot\mathbf{d}\gamma(T\pi_{Q}(X_{\mathcal{K}}\cdot
\gamma), \; T\pi_{Q}(w))\\
& = \omega_{\mathcal{K}_{\mathcal{O}_\mu}}\cdot
\tau_{\mathcal{K}_{\mathcal{O}_\mu}}(T\pi_{\mathcal{O}_\mu}\cdot (X_{\mathcal{K}}\cdot \gamma), \;
T\pi_{\mathcal{O}_\mu}\cdot w) \\
& \;\;\;\;\;\; -\omega_{\mathcal{K}_{\mathcal{O}_\mu}}\cdot \tau_{\mathcal{K}_{\mathcal{O}_\mu}}
(T\pi_{\mathcal{O}_\mu}\cdot (X_{\mathcal{K}}\cdot \gamma), \;
T\bar{\gamma}_{\mathcal{O}_\mu}\cdot T\pi_{Q}(w))\\
& \;\;\;\;\;\; -\tau_{\mathcal{K}_{\mathcal{O}_\mu}}\cdot i_{\mathcal{M}_{\mathcal{O}_\mu}}^* \cdot\mathbf{d}\gamma (T\pi_{Q}(X_{\mathcal{K}}\cdot \gamma), \; T\pi_{Q}(w))\\
& = \omega_{\mathcal{K}_{\mathcal{O}_\mu}}( \tau_{\mathcal{K}_{\mathcal{O}_\mu}}
\cdot T\pi_{\mathcal{O}_\mu} \cdot (X_{\mathcal{K}}\cdot \gamma), \; \tau_{\mathcal{K}_{\mathcal{O}_\mu}}\cdot T\pi_{\mathcal{O}_\mu}\cdot w) \\
 & \;\;\;\;\;\; - \omega_{\mathcal{K}_{\mathcal{O}_\mu}}(\tau_{\mathcal{K}_{\mathcal{O}_\mu}}\cdot
T\pi_{\mathcal{O}_\mu}\cdot (X_{\mathcal{K}}\cdot \gamma), \;
\tau_{\mathcal{K}_{\mathcal{O}_\mu}}\cdot T\bar{\gamma}_{\mathcal{O}_\mu}\cdot T\pi_{Q}(w))\\
& \;\;\;\;\;\; -\tau_{\mathcal{K}_{\mathcal{O}_\mu}}\cdot i_{\mathcal{M}_{\mathcal{O}_\mu}}^* \cdot\mathbf{d}\gamma (T\pi_{Q}(X_{\mathcal{K}}\cdot \gamma), \; T\pi_{Q}(w))\\
& = \omega_{\mathcal{K}_{\mathcal{O}_\mu}}(X_{\mathcal{K}_{\mathcal{O}_\mu}}\cdot
\bar{\gamma}_{\mathcal{O}_\mu}, \; \tau_{\mathcal{K}_{\mathcal{O}_\mu}}\cdot T\pi_{\mathcal{O}_\mu}\cdot w)
- \omega_{\mathcal{K}_{\mathcal{O}_\mu}}(X_{\mathcal{K}_{\mathcal{O}_\mu}} \cdot
\bar{\gamma}_{\mathcal{O}_\mu}, \; \tau_{\mathcal{K}_{\mathcal{O}_\mu}}\cdot
T\bar{\gamma}_{\mathcal{O}_\mu}\cdot T\pi_{Q}(w))\\
& \;\;\;\;\;\; -\tau_{\mathcal{K}_{\mathcal{O}_\mu}}\cdot i_{\mathcal{M}_{\mathcal{O}_\mu}}^*
\cdot\mathbf{d}\gamma(T\pi_{Q}(X_{\mathcal{K}}\cdot \gamma), \; T\pi_{Q}(w)),
\end{align*}
where we have used that $ \tau_{\mathcal{K}_{\mathcal{O}_\mu}}\cdot T\bar{\gamma}_{\mathcal{O}_\mu}=
T\bar{\gamma}_{\mathcal{O}_\mu}, $ since $\textmd{Im}(T\bar{\gamma}_{\mathcal{O}_\mu})\subset \mathcal{K}_{\mathcal{O}_\mu}, $ and $\tau_{\mathcal{K}_{\mathcal{O}_\mu}}\cdot
T\pi_{\mathcal{O}_\mu} \cdot (X_{\mathcal{K}}\cdot \gamma )
= X_{\mathcal{K}_{\mathcal{O}_\mu}}\cdot \bar{\gamma}_{\mathcal{O}_\mu}. $
If the one-form $\gamma: Q \rightarrow T^*Q $ is closed on $\mathcal{D}$ with respect to
$T\pi_Q: TT^* Q \rightarrow TQ, $ then we have that
$\mathbf{d}\gamma(T\pi_{Q}(X_{\mathcal{K}}\cdot \gamma), \; T\pi_{Q}(w))=0, $
since $X_{\mathcal{K}}\cdot \gamma, \; w \in \mathcal{F},$ and
$T\pi_{Q}(X_{\mathcal{K}}\cdot \gamma), \; T\pi_{Q}(w) \in \mathcal{D}, $ and hence
$$
\tau_{\mathcal{K}_{\mathcal{O}_\mu}}\cdot
i_{\mathcal{M}_{\mathcal{O}_\mu}}^* \cdot\mathbf{d}\gamma(T\pi_{Q}(X_{\mathcal{K}}\cdot \gamma),
\; T\pi_{Q}(w))=0,
$$
and
\begin{align}
& \omega_{\mathcal{K}_{\mathcal{O}_\mu}}(T\bar{\gamma}_{\mathcal{O}_\mu} \cdot \tilde{X}^\gamma, \;
\tau_{\mathcal{K}_{\mathcal{O}_\mu}}\cdot T\pi_{\mathcal{O}_\mu} \cdot w)
- \omega_{\mathcal{K}_{\mathcal{O}_\mu}}(X_{\mathcal{K}_{\mathcal{O}_\mu}}\cdot
\bar{\gamma}_{\mathcal{O}_\mu}, \; \tau_{\mathcal{K}_{\mathcal{O}_\mu}} \cdot T\pi_{\mathcal{O}_\mu} \cdot w) \nonumber \\
& = -\omega_{\mathcal{K}_{\mathcal{O}_\mu}}(X_{\mathcal{K}_{\mathcal{O}_\mu}} \cdot
\bar{\gamma}_{\mathcal{O}_\mu}, \; \tau_{\mathcal{K}_{\mathcal{O}_\mu}}\cdot
T\bar{\gamma}_{\mathcal{O}_\mu}\cdot T\pi_{Q}(w)).
\label{5.13} \end{align}
If $\bar{\gamma}_{\mathcal{O}_\mu}$ satisfies the equation
$T\bar{\gamma}_{\mathcal{O}_\mu}\cdot \tilde{X}^\gamma= X_{\mathcal{K}_{\mathcal{O}_\mu}}\cdot
\bar{\gamma}_{\mathcal{O}_\mu},$
from Lemma 3.3(i) we know that the right side of (5.13) becomes that
\begin{align*}
& -\omega_{\mathcal{K}_{\mathcal{O}_\mu}}(X_{\mathcal{K}_{\mathcal{O}_\mu}} \cdot
\bar{\gamma}_{\mathcal{O}_\mu}, \; \tau_{\mathcal{K}_{\mathcal{O}_\mu}}\cdot T\bar{\gamma}_{\mathcal{O}_\mu}
\cdot T\pi_{Q}(w))\\
& = -\omega_{\mathcal{K}_{\mathcal{O}_\mu}}(T\bar{\gamma}_{\mathcal{O}_\mu}\cdot \tilde{X}^\gamma,
\; \tau_{\mathcal{K}_{\mathcal{O}_\mu}}\cdot T\bar{\gamma}_{\mathcal{O}_\mu} \cdot T\pi_{Q}(w))\\
& = - \omega_{\mathcal{K}_{\mathcal{O}_\mu}}
(\tau_{\mathcal{K}_{\mathcal{O}_\mu}}T\bar{\gamma}_{\mathcal{O}_\mu}\cdot \tilde{X}^\gamma,
\; \tau_{\mathcal{K}_{\mathcal{O}_\mu}}\cdot T\bar{\gamma}_{\mathcal{O}_\mu} \cdot T\pi_{Q}(w))\\
& = -\tau_{\mathcal{K}_{\mathcal{O}_\mu}}\cdot i_{\mathcal{M}_{\mathcal{O}_\mu}}^* \cdot
\omega_{\mathcal{O}_\mu}(T\bar{\gamma}_{\mathcal{O}_\mu}\cdot \tilde{X}^\gamma,
\; T\bar{\gamma}_{\mathcal{O}_\mu} \cdot T\pi_{Q}(w))\\
& = -\tau_{\mathcal{K}_{\mathcal{O}_\mu}}\cdot i_{\mathcal{M}_{\mathcal{O}_\mu}}^* \cdot \bar{\gamma}_{\mathcal{O}_\mu}^* \cdot
\omega_{\mathcal{O}_\mu}( T\pi_{Q} \cdot \tilde{X} \cdot
\gamma, \; T\pi_{Q}(w))\\
& = -\tau_{\mathcal{K}_{\mathcal{O}_\mu}}\cdot i_{\mathcal{M}_{\mathcal{O}_\mu}}^* \cdot
 \gamma^* \cdot \pi^*_{\mathcal{O}_\mu}\cdot \omega_{\mathcal{O}_\mu}(T\pi_{Q} \cdot X_{\mathcal{K}}\cdot
\gamma, \; T\pi_{Q}(w))\\
& = -\tau_{\mathcal{K}_{\mathcal{O}_\mu}}\cdot
i_{\mathcal{M}_{\mathcal{O}_\mu}}^* \cdot\gamma^*\omega( T\pi_{Q}(X_{\mathcal{K}}\cdot\gamma), \; T\pi_{Q}(w))\\
& = \tau_{\mathcal{K}_{\mathcal{O}_\mu}}\cdot i_{\mathcal{M}_{\mathcal{O}_\mu}}^* \cdot
\mathbf{d}\gamma(T\pi_{Q}(X_{\mathcal{K}}\cdot\gamma ), \; T\pi_{Q}(w))=0.
\end{align*}
But, because the $\mathbf{J}$-nonholonomic $R_o$-reduced distributional two-form
$\omega_{\mathcal{K}_{\mathcal{O}_\mu}}$ is non-degenerate,
the left side of (5.13) equals zero, only when
$\bar{\gamma}_{\mathcal{O}_\mu}$ satisfies the equation
$T\bar{\gamma}_{\mathcal{O}_\mu}\cdot \tilde{X}^\gamma= X_{\mathcal{K}_{\mathcal{O}_\mu}}\cdot
\bar{\gamma}_{\mathcal{O}_\mu}.$ Thus,
if the one-form $\gamma: Q \rightarrow T^*Q $ is closed on $\mathcal{D}$ with respect to
$T\pi_Q: TT^* Q \rightarrow TQ, $ then $\bar{\gamma}_{\mathcal{O}_\mu}$ must be a solution
of the Type I of Hamilton-Jacobi equation
$T\bar{\gamma}_{\mathcal{O}_\mu}\cdot \tilde{X}^\gamma
= X_{\mathcal{K}_{\mathcal{O}_\mu}}\cdot
\bar{\gamma}_{\mathcal{O}_\mu}. $
\hskip 0.3cm $\blacksquare$\\

Next, for any $G$-invariant symplectic map $\varepsilon: T^* Q \rightarrow T^* Q $,
we can prove the following Type II of Hamilton-Jacobi theorem for the
$\mathbf{J}$-nonholonomic $R_o$-reduced distributional RCH system.
For convenience, the maps involved in the following
theorem and its proof are shown in Diagram-8.
\begin{center}
\hskip 0cm \xymatrix{ \mathbf{J}^{-1}(\mathcal{O}_\mu) \ar[r]^{i_{\mathcal{O}_\mu}} & T^* Q
\ar[d]_{X_{H\cdot \varepsilon}} \ar[dr]^{\tilde{X}^\varepsilon} \ar[r]^{\pi_Q}
& Q \ar[r]^{\gamma} & T^*Q \ar[d]_{\tilde{X}} \ar[dr]^{X_{h_{\mathcal{K}_{\mathcal{O}_\mu}}}
\cdot\bar{\varepsilon}} \ar[r]^{\pi_{\mathcal{O}_\mu}}
& (T^* Q)_{\mathcal{O}_\mu} \ar[d]^{X_{h_{\mathcal{K}_{\mathcal{O}_\mu}}}}
& \mathcal{M}_{\mathcal{O}_\mu} \ar[l]_{i_{\mathcal{M}_{\mathcal{O}_\mu}}}
\ar[d]_{X_{\mathcal{K}_{\mathcal{O}_\mu}}}\\
& T(T^*Q)  & TQ \ar[l]^{T\gamma} & T(T^*Q) \ar[l]^{T\pi_Q} \ar[r]_{T\pi_{\mathcal{O}_\mu}}
& T(T^* Q)_{\mathcal{O}_\mu} \ar[r]^{\tau_{\mathcal{K}_{\mathcal{O}_\mu}}} & \mathcal{K}_{\mathcal{O}_\mu} }
\end{center}
$$\mbox{Diagram-8}$$

\begin{theo} (Type II of Hamilton-Jacobi Theorem for a $\mathbf{J}$-Nonholonomic
$R_o$-reduced Distributional RCH System) For a given
$\mathbf{J}$-nonholonomic regular orbit reducible RCH system
$(T^*Q,G,\omega,\mathbf{J},\\ \mathcal{D},H, F, u)$ with the associated distributional RCH system
$(\mathcal{K},\omega_{\mathcal {K}},H_{\mathcal{K}}, F_{\mathcal{K}}, u_{\mathcal{K}})$
and the $\mathbf{J}$-nonholonomic $R_o$-reduced distributional RCH system
$(\mathcal{K}_{\mathcal{O}_\mu},\omega_{\mathcal{K}_{\mathcal{O}_\mu}},
h_{\mathcal{K}_{\mathcal{O}_\mu}},
f_{\mathcal {K}_{\mathcal{O}_\mu}}, u_{\mathcal {K}_{\mathcal{O}_\mu}})$, assume that $\gamma:
Q \rightarrow T^*Q$ is an one-form on $Q$, and $\lambda=\gamma \cdot
\pi_{Q}: T^* Q \rightarrow T^* Q, $ and for any
symplectic map $\varepsilon:T^* Q \rightarrow T^* Q, $
denote $\tilde{X}^\varepsilon = T\pi_{Q}\cdot \tilde{X}\cdot \varepsilon$,
where $\tilde{X}=X_{(\mathcal{K},\omega_{\mathcal{K}},
H_{\mathcal{K}}, F_{\mathcal{K}}, u_{\mathcal{K}})}
=X_\mathcal {K}+ F_{\mathcal{K}}+u_{\mathcal{K}}$
is the dynamical vector field of the distributional RCH system
$(\mathcal{K},\omega_{\mathcal{K}},
H_{\mathcal{K}}, F_{\mathcal{K}}, u_{\mathcal{K}})$
corresponding to the $\mathbf{J}$-nonholonomic regular orbit reducible RCH
system with symmetry and momentum map $(T^*Q,G,\omega,\mathbf{J},\mathcal{D},H, F,u)$.
Moreover, assume that $\mu\in\mathfrak{g}^\ast$ is a regular value of the momentum
map $\mathbf{J}$, and $\textmd{Im}(\gamma)\subset \mathcal{M} \cap
\mathbf{J}^{-1}(\mu), $ and it is $G$-invariant,
and $\varepsilon$ is also $G$-invariant and
$\varepsilon(\mathcal{M}\cap \mathbf{J}^{-1}(\mathcal{O}_\mu))
\subset \mathcal{M}\cap \mathbf{J}^{-1}(\mathcal{O}_\mu). $ Denote
$\bar{\gamma}_{\mathcal{O}_\mu}
=\pi_{\mathcal{O}_\mu}(\gamma): Q \rightarrow \mathcal{M}_{\mathcal{O}_\mu} $,
and $ \textmd{Im}(T\bar{\gamma}_{\mathcal{O}_\mu})\subset \mathcal{K}_{\mathcal{O}_\mu}, $ and
$\bar{\lambda}_{\mathcal{O}_\mu}=\pi_{\mathcal{O}_\mu}(\lambda): \mathcal{M}\cap \mathbf{J}^{-1}(\mathcal{O}_\mu)
(\subset T^* Q) \rightarrow \mathcal{M}_{\mathcal{O}_\mu}, $ and
$\bar{\varepsilon}_{\mathcal{O}_\mu}=\pi_{\mathcal{O}_\mu}(\varepsilon):
\mathcal{M} \cap \mathbf{J}^{-1}(\mathcal{O}_\mu) (\subset T^* Q) \rightarrow
\mathcal{M}_{\mathcal{O}_\mu}. $ Then $\varepsilon$ and $\bar{\varepsilon}_{\mathcal{O}_\mu}$ satisfy the
equation $\tau_{\mathcal{K}_{\mathcal{O}_\mu}} \cdot T\bar{\varepsilon}(X_{h_{\mathcal{K}_{\mathcal{O}_\mu}}\cdot \bar{\varepsilon}_{\mathcal{O}_\mu}})
= T\bar{\lambda}_{\mathcal{O}_\mu} \cdot \tilde{X} \cdot\varepsilon, $ if and only if
they satisfy the equation $T\bar{\gamma}_{\mathcal{O}_\mu}\cdot
\tilde{X}^\varepsilon= X_{\mathcal{K}_{\mathcal{O}_\mu}}\cdot \bar{\varepsilon}_{\mathcal{O}_\mu}. $
Here $X_{h_{\mathcal{K}_{\mathcal{O}_\mu}} \cdot\bar{\varepsilon}_{\mathcal{O}_\mu}}$
is the Hamiltonian vector field of the
function $h_{\mathcal{K}_{\mathcal{O}_\mu}}\cdot \bar{\varepsilon}_{\mathcal{O}_\mu}: T^* Q\rightarrow \mathbb{R}, $
and $X_{\mathcal{K}_{\mathcal{O}_\mu}}$ is the $\mathbf{J}$-nonholonomic
$R_o$-reduced dynamical vector field. The equation
$T\bar{\gamma}_{\mathcal{O}_\mu}\cdot \tilde{X}^\varepsilon= X_{\mathcal{K}_{\mathcal{O}_\mu}}\cdot
\bar{\varepsilon}_{\mathcal{O}_\mu},$ is called the Type II of Hamilton-Jacobi equation for the
$\mathbf{J}$-nonholonomic $R_o$-reduced distributional RCH system
$(\mathcal{K}_{\mathcal{O}_\mu},\omega_{\mathcal{K}_{\mathcal{O}_\mu}},
h_{\mathcal{K}_{\mathcal{O}_\mu}},
f_{\mathcal {K}_{\mathcal{O}_\mu}}, u_{\mathcal {K}_{\mathcal{O}_\mu}})$.
\end{theo}

\noindent{\bf Proof: } In the same way, for the dynamical vector field of the distributional RCH system
$(\mathcal{K},\omega_{\mathcal {K}}, H_{\mathcal{K}}, F_{\mathcal{K}}, u_{\mathcal{K}})$,
$\tilde{X}=X_{(\mathcal{K},\omega_{\mathcal{K}},
H_{\mathcal{K}}, F_{\mathcal{K}}, u_{\mathcal{K}})}
=X_\mathcal {K}+ F_{\mathcal{K}}+u_{\mathcal{K}}$,
and $F_{\mathcal{K}}=\tau_{\mathcal{K}}\cdot \textnormal{vlift}(F_{\mathcal{M}})X_H$,
and $u_{\mathcal{K}}=\tau_{\mathcal{K}}\cdot \textnormal{vlift}(u_{\mathcal{M}})X_H$,
note that $T\pi_{Q}\cdot \textnormal{vlift}(F_{\mathcal{M}})X_H=T\pi_{Q}\cdot \textnormal{vlift}(u_{\mathcal{M}})X_H=0, $
then we have that $T\pi_{Q}\cdot F_{\mathcal{K}}=T\pi_{Q}\cdot u_{\mathcal{K}}=0,$
and hence $T\pi_{Q}\cdot \tilde{X}\cdot \varepsilon
=T\pi_{Q}\cdot X_{\mathcal{K}}\cdot \varepsilon. $
Next, we note that
the $R_o$-reduced symplectic space
$(T^\ast Q)_{\mathcal{O}_\mu}= \mathbf{J}^{-1}(\mathcal{O}_\mu)/G
\cong (\mathbf{J}^{-1}(\mu)/G) \times \mathcal{O}_\mu, $ with the $R_o$-reduced
symplectic form $\omega_{\mathcal{O}_\mu}$ uniquely characterized by
the relation $i_{\mathcal{O}_\mu}^\ast
\omega=\pi_{\mathcal{O}_{\mu}}^\ast \omega_{\mathcal{O}
_\mu}+\mathbf{J}_{\mathcal{O}_\mu}^\ast\omega_{\mathcal{O}_\mu}^+. $
Since
$\textmd{Im}(\gamma)\subset \mathcal{M}\cap \mathbf{J}^{-1}(\mu), $
and it is $G$-invariant,
in this case for any $V\in TQ, $ and $w\in TT^*Q, $
we have that
$\mathbf{J}_{\mathcal{O}_\mu}^\ast\omega_{\mathcal{O}_\mu}^+(T\gamma
\cdot V, \; w)=0, $ and hence
$\pi_{\mathcal{O}_\mu}^*\omega_{\mathcal{O}_\mu}=
i_{\mathcal{O}_\mu}^*\omega= \omega, $ along $\textmd{Im}(\gamma)$.
On the other hand, because
$\textmd{Im}(T\bar{\gamma}_{\mathcal{O}_\mu})\subset \mathcal{K}_{\mathcal{O}_\mu}, $
then $\omega_{\mathcal{K}_{\mathcal{O}_\mu}}\cdot
\tau_{\mathcal{K}_{\mathcal{O}_\mu}}=\tau_{\mathcal{K}_{\mathcal{O}_\mu}}\cdot
\omega_{\mathcal{M}_{\mathcal{O}_\mu}}= \tau_{\mathcal{K}_{\mathcal{O}_\mu}}\cdot
i_{\mathcal{M}_{\mathcal{O}_\mu}}^* \cdot \omega_{\mathcal{O}_\mu}, $ along
$\textmd{Im}(T\bar{\gamma}_{\mathcal{O}_\mu})$.
Moreover, from the distributional Hamiltonian equation (5.9),
we have that $X_{\mathcal{K}}= \tau_{\mathcal{K}}\cdot X_H.$
Note that $\varepsilon(\mathcal{M})\subset \mathcal{M},$ and
$T\pi_{Q}(X_H\cdot \varepsilon(q,p))\in
\mathcal{D}_{q}, \; \forall q \in Q, \; (q,p) \in \mathcal{M}(\subset T^* Q), $
and hence $X_H\cdot \varepsilon \in \mathcal{F}$ along $\varepsilon$,
and $\tau_{\mathcal{K}}\cdot X_{H}\cdot \varepsilon
= X_{\mathcal{K}}\cdot \varepsilon \in \mathcal{K}$.
Because the vector fields $X_{\mathcal{K}}$
and $X_{\mathcal{K}_{\mathcal{O}_\mu}}$ are $\pi_{\mathcal{O}_\mu}$-related, then
$T\pi_{\mathcal{O}_\mu}(X_{\mathcal{K}})=X_{\mathcal{K}_{\mathcal{O}_\mu}}\cdot \pi_{\mathcal{O}_\mu}$,
and hence $\tau_{\mathcal{K}_{\mathcal{O}_\mu}}\cdot T\pi_{\mathcal{O}_\mu}(X_{\mathcal{K}}\cdot \varepsilon)
=\tau_{\mathcal{K}_{\mathcal{O}_\mu}}\cdot (T\pi_{\mathcal{O}_\mu}(X_{\mathcal{K}}))\cdot (\varepsilon)
= \tau_{\mathcal{K}_{\mathcal{O}_\mu}}\cdot (X_{\mathcal{K}_{\mathcal{O}_\mu}}\cdot \pi_{\mathcal{O}_\mu})\cdot (\varepsilon)
= \tau_{\mathcal{K}_{\mathcal{O}_\mu}}\cdot X_{\mathcal{K}_{\mathcal{O}_\mu}}\cdot \pi_{\mathcal{O}_\mu}(\varepsilon)
= X_{\mathcal{K}_{\mathcal{O}_\mu}}\cdot \bar{\varepsilon}.$
Thus, using the $\mathbf{J}$-nonholonomic $R_o$-reduced distributional two-form
$\omega_{\mathcal{K}_{\mathcal{O}_\mu}}$, from Lemma 3.3(ii) and (iii), if we take that $v=
X_{\mathcal{K}}\cdot \varepsilon \in \mathcal{K} (\subset \mathcal{F}), $ and
for any $w \in \mathcal{F}, \; T\lambda(w)\neq 0, $ and
$\tau_{\mathcal{K}_{\mathcal{O}_\mu}}\cdot T\pi_{\mathcal{O}_\mu} \cdot w \neq 0, $
and $\tau_{\mathcal{K}_{\mathcal{O}_\mu}}\cdot T\pi_{\mathcal{O}_\mu} \cdot T\lambda(w) \neq 0, $
then we have that
\begin{align*}
& \omega_{\mathcal{K}_{\mathcal{O}_\mu}}(T\bar{\gamma}_{\mathcal{O}_\mu} \cdot \tilde{X}^\varepsilon, \;
\tau_{\mathcal{K}_{\mathcal{O}_\mu}}\cdot T\pi_{\mathcal{O}_\mu} \cdot w)=
\omega_{\mathcal{K}_{\mathcal{O}_\mu}}(\tau_{\mathcal{K}_{\mathcal{O}_\mu}}\cdot
T\bar{\gamma}_{\mathcal{O}_\mu} \cdot \tilde{X}^\varepsilon, \;
\tau_{\mathcal{K}_{\mathcal{O}_\mu}}\cdot T\pi_{\mathcal{O}_\mu} \cdot w)\\
& = \tau_{\mathcal{K}_{\mathcal{O}_\mu}}\cdot \omega_{\mathcal{M}_{\mathcal{O}_\mu}}(T(\pi_{\mathcal{O}_\mu}
\cdot\gamma )\cdot \tilde{X}^\varepsilon, \; T\pi_{\mathcal{O}_\mu} \cdot w ) =
\tau_{\mathcal{K}_{\mathcal{O}_\mu}}\cdot i_{\mathcal{M}_{\mathcal{O}_\mu}}^* \cdot
\omega_{\mathcal{O}_\mu}(T\pi_{\mathcal{O}_\mu} \cdot T\gamma \cdot \tilde{X}^\varepsilon,
\; T\pi_{\mathcal{O}_\mu} \cdot w )\\
& = \tau_{\mathcal{K}_{\mathcal{O}_\mu}}\cdot i_{\mathcal{M}_{\mathcal{O}_\mu}}^*
\cdot \pi_{\mathcal{O}_\mu}^*\omega_{\mathcal{O}_\mu}(T\gamma \cdot T\pi_Q \cdot \tilde{X} \cdot
\varepsilon, \; w) = \tau_{\mathcal{K}_{\mathcal{O}_\mu}}\cdot i_{\mathcal{M}_{\mathcal{O}_\mu}}^*
\cdot \omega(T(\gamma \cdot \pi_Q) \cdot X_{\mathcal{K}}\cdot \varepsilon, \; w)\\
& =\tau_{\mathcal{K}_{\mathcal{O}_\mu}}\cdot i_{\mathcal{M}_{\mathcal{O}_\mu}}^* \cdot
(\omega(X_{\mathcal{K}}\cdot \varepsilon, \; w-T(\gamma
\cdot \pi_Q)\cdot w) -\mathbf{d}\gamma(T\pi_{Q}(X_{\mathcal{K}}\cdot \varepsilon), \; T\pi_{Q}(w)))\\
& = \tau_{\mathcal{K}_{\mathcal{O}_\mu}}\cdot i_{\mathcal{M}_{\mathcal{O}_\mu}}^* \cdot
\omega(X_{\mathcal{K}}\cdot \varepsilon, \; w) - \tau_{\mathcal{K}_{\mathcal{O}_\mu}}\cdot
i_{\mathcal{M}_{\mathcal{O}_\mu}}^* \cdot \omega(X_{\mathcal{K}}\cdot
\varepsilon, \; T\lambda \cdot w)\\
& \;\;\;\;\;\; -\tau_{\mathcal{K}_{\mathcal{O}_\mu}}\cdot i_{\mathcal{M}_{\mathcal{O}_\mu}}^* \cdot\mathbf{d}\gamma(T\pi_{Q}(X_{\mathcal{K}}\cdot \varepsilon), \; T\pi_{Q}(w))\\
& = \tau_{\mathcal{K}_{\mathcal{O}_\mu}}\cdot i_{\mathcal{M}_{\mathcal{O}_\mu}}^*
\cdot\pi_{\mathcal{O}_\mu}^*\omega_{\mathcal{O}_\mu}(X_{\mathcal{K}}\cdot \varepsilon, \; w) -
\tau_{\mathcal{K}_{\mathcal{O}_\mu}}\cdot i_{\mathcal{M}_{\mathcal{O}_\mu}}^*
\cdot\pi_{\mathcal{O}_\mu}^*\omega_{\mathcal{O}_\mu}(X_{\mathcal{K}}\cdot
\varepsilon, \; T\lambda \cdot w)\\
& \;\;\;\;\;\; +\tau_{\mathcal{K}_{\mathcal{O}_\mu}}\cdot i_{\mathcal{M}_{\mathcal{O}_\mu}}^*
\cdot \lambda^* \omega(X_{\mathcal{K}}\cdot \varepsilon, \; w)\\
& = \tau_{\mathcal{K}_{\mathcal{O}_\mu}}\cdot i_{\mathcal{M}_{\mathcal{O}_\mu}}^*
\cdot\omega_{\mathcal{O}_\mu}(T\pi_{\mathcal{O}_\mu}(X_{\mathcal{K}}\cdot \varepsilon),
\; T\pi_{\mathcal{O}_\mu} \cdot w)\\
& \;\;\;\;\;\; -\tau_{\mathcal{K}_{\mathcal{O}_\mu}}\cdot i_{\mathcal{M}_{\mathcal{O}_\mu}}^* \cdot
\omega_{\mathcal{O}_\mu}(T\pi_{\mathcal{O}_\mu}\cdot(X_{\mathcal{K}}\cdot \varepsilon),
\; T(\pi_{\mathcal{O}_\mu}\cdot\lambda) \cdot w)\\
& \;\;\;\;\;\; +\tau_{\mathcal{K}_{\mathcal{O}_\mu}}\cdot
i_{\mathcal{M}_{\mathcal{O}_\mu}}^* \cdot
\pi_{\mathcal{O}_\mu}^*\omega_{\mathcal{O}_\mu}(T\lambda\cdot X_{\mathcal{K}}\cdot
\varepsilon, \; T\lambda \cdot w)\\
& = \omega_{\mathcal{K}_{\mathcal{O}_\mu}}\cdot
\tau_{\mathcal{K}_{\mathcal{O}_\mu}}(T\pi_{\mathcal{O}_\mu}(X_{\mathcal{K}}\cdot
\varepsilon), \; T\pi_{\mathcal{O}_\mu}\cdot w) - \omega_{\mathcal{K}_{\mathcal{O}_\mu}}\cdot
\tau_{\mathcal{K}_{\mathcal{O}_\mu}}(T\pi_{\mathcal{O}_\mu}(X_{\mathcal{K}}\cdot
\varepsilon), \; T\bar{\lambda}_{\mathcal{O}_\mu}\cdot w)\\
& \;\;\;\;\;\; +\tau_{\mathcal{K}_{\mathcal{O}_\mu}}\cdot
i_{\mathcal{M}_{\mathcal{O}_\mu}}^* \cdot
\omega_{\mathcal{O}_\mu}(T\pi_{\mathcal{O}_\mu}\cdot T\lambda\cdot X_{\mathcal{K}}\cdot
\varepsilon, \; T\pi_{\mathcal{O}_\mu}\cdot T\lambda \cdot w)\\
& = \omega_{\mathcal{K}_{\mathcal{O}_\mu}}( \tau_{\mathcal{K}_{\mathcal{O}_\mu}}\cdot T\pi_{\mathcal{O}_\mu}\cdot(X_{\mathcal{K}}\cdot \varepsilon), \;
\tau_{\mathcal{K}_{\mathcal{O}_\mu}}\cdot T\pi_{\mathcal{O}_\mu}\cdot w) \\
& \;\;\;\;\;\;  - \omega_{\mathcal{K}_{\mathcal{O}_\mu}}(\tau_{\mathcal{K}_{\mathcal{O}_\mu}}\cdot
T\pi_{\mathcal{O}_\mu}\cdot(X_{\mathcal{K}}\cdot \varepsilon), \;
\tau_{\mathcal{K}_{\mathcal{O}_\mu}}\cdot T\bar{\lambda}_{\mathcal{O}_\mu}\cdot w)\\
& \;\;\;\;\;\; +\omega_{\mathcal{K}_{\mathcal{O}_\mu}}(\tau_{\mathcal{K}_{\mathcal{O}_\mu}}\cdot
T\bar{\lambda}_{\mathcal{O}_\mu}\cdot X_{\mathcal{K}} \cdot
\varepsilon, \; \tau_{\mathcal{K}_{\mathcal{O}_\mu}}\cdot T\bar{\lambda}_{\mathcal{O}_\mu}\cdot w)\\
& = \omega_{\mathcal{K}_{\mathcal{O}_\mu}}(X_{\mathcal{K}_{\mathcal{O}_\mu}}\cdot
\bar{\varepsilon}_{\mathcal{O}_\mu}, \;
\tau_{\mathcal{K}_{\mathcal{O}_\mu}} \cdot T\pi_{\mathcal{O}_\mu} \cdot w)
- \omega_{\mathcal{K}_{\mathcal{O}_\mu}}(\tau_{\mathcal{K}_{\mathcal{O}_\mu}}\cdot X_{h_{\mathcal{K}_{\mathcal{O}_\mu}}} \cdot
\bar{\varepsilon}_{\mathcal{O}_\mu}, \; \tau_{\mathcal{K}_{\mathcal{O}_\mu}} \cdot T\bar{\lambda}_{\mathcal{O}_\mu} \cdot w)\\
& \;\;\;\;\;\; +\omega_{\mathcal{K}_{\mathcal{O}_\mu}}(T\bar{\lambda}_{\mathcal{O}_\mu}
\cdot X_{\mathcal{K}} \cdot
\varepsilon, \; \tau_{\mathcal{K}_{\mathcal{O}_\mu}} \cdot T\bar{\lambda}_{\mathcal{O}_\mu} \cdot w),
\end{align*}
where we have used that $ \tau_{\mathcal{K}_{\mathcal{O}_\mu}}\cdot T\bar{\gamma}_{\mathcal{O}_\mu}=
T\bar{\gamma}_{\mathcal{O}_\mu}, $ and $ \tau_{\mathcal{K}_{\mathcal{O}_\mu}}\cdot T\bar{\lambda}_{\mathcal{O}_\mu}=T\bar{\lambda}_{\mathcal{O}_\mu}, $
since $\textmd{Im}(T\bar{\gamma}_{\mathcal{O}_\mu})\subset \mathcal{K}_{\mathcal{O}_\mu}, $
and $\tau_{\mathcal{K}_{\mathcal{O}_\mu}}\cdot
T\pi_{\mathcal{O}_\mu}\cdot(X_{\mathcal{K}}\cdot \varepsilon)
= X_{\mathcal{K}_{\mathcal{O}_\mu}}\cdot \bar{\varepsilon}_{\mathcal{O}_\mu}. $
From the $\mathbf{J}$-nonholonomic $R_o$-reduced distributional Hamiltonian equation (5.10)
$\mathbf{i}_{X_{\mathcal{K_{\mathcal{O}_\mu}}}}\omega_{\mathcal{K_{\mathcal{O}_\mu}}}
=\mathbf{d}h_{\mathcal{K_{\mathcal{O}_\mu}}}$, we have that
$X_{\mathcal{K_{\mathcal{O}_\mu}}}
=\tau_{\mathcal{K_{\mathcal{O}_\mu}}}\cdot X_{h_{\mathcal{K_{\mathcal{O}_\mu}}}},$
where $ X_{h_{\mathcal{K_{\mathcal{O}_\mu}}}}$ is the Hamiltonian vector field of
the function $h_{\mathcal{K_{\mathcal{O}_\mu}}}.$
Note that $\varepsilon: T^* Q \rightarrow T^* Q $ is symplectic, and
$\pi_{\mathcal{O}_\mu}^*\omega_{\mathcal{O}_\mu}= i_{\mathcal{O}_\mu}^*\omega= \omega, $ along
$\textmd{Im}(\gamma)$, and hence $\bar{\varepsilon}_{\mathcal{O}_\mu}=
\pi_{\mathcal{O}_\mu}(\varepsilon): T^* Q \rightarrow (T^* Q)_{\mathcal{O}_\mu} $ is also symplectic
along $\bar{\varepsilon}_{\mathcal{O}_\mu}$, and hence $X_{h_{\mathcal{K}_{\mathcal{O}_\mu}}}\cdot
\bar{\varepsilon}_{\mathcal{O}_\mu}= T\bar{\varepsilon}_{\mathcal{O}_\mu} \cdot X_{h_{\mathcal{K}_{\mathcal{O}_\mu}} \cdot
\bar{\varepsilon}_{\mathcal{O}_\mu}}, $ along $\bar{\varepsilon}_{\mathcal{O}_\mu}$, and hence
$X_{\mathcal{K}_{\mathcal{O}_\mu}}\cdot
\bar{\varepsilon}_{\mathcal{O}_\mu}=
\tau_{\mathcal{K}_{\mathcal{O}_\mu}}\cdot X_{h_{\mathcal{K}_{\mathcal{O}_\mu}}} \cdot
\bar{\varepsilon}_{\mathcal{O}_\mu}= \tau_{\mathcal{K}_{\mathcal{O}_\mu}}\cdot T\bar{\varepsilon}_{\mathcal{O}_\mu}
\cdot X_{h_{\mathcal{K}_{\mathcal{O}_\mu}}\cdot \bar{\varepsilon}_{\mathcal{O}_\mu}}, $
along $\bar{\varepsilon}_{\mathcal{O}_\mu}. $
Note that
$T\bar{\lambda}_{\mathcal{O}_\mu}\cdot X_\mathcal{K}\cdot \varepsilon
=T\pi_{\mathcal{O}_\mu}\cdot T\lambda \cdot X_\mathcal{K}\cdot \varepsilon
=T\pi_{\mathcal{O}_\mu}\cdot T\gamma \cdot T\pi_Q\cdot X_\mathcal{K}\cdot \varepsilon
=T\pi_{\mathcal{O}_\mu}\cdot T\gamma \cdot T\pi_Q\cdot \tilde{X}\cdot \varepsilon
=T\pi_{\mathcal{O}_\mu}\cdot T\lambda \cdot \tilde{X}\cdot \varepsilon
=T\bar{\lambda}_{\mathcal{O}_\mu}\cdot \tilde{X}\cdot \varepsilon.$
Then we have that
\begin{align*}
& \omega_{\mathcal{K}_{\mathcal{O}_\mu}}(T\bar{\gamma}_{\mathcal{O}_\mu} \cdot \tilde{X}^\varepsilon, \;
\tau_{\mathcal{K}_{\mathcal{O}_\mu}}\cdot T\pi_{\mathcal{O}_\mu} \cdot w)-
\omega_{\mathcal{K}_{\mathcal{O}_\mu}}(X_{\mathcal{K}_{\mathcal{O}_\mu}}\cdot \bar{\varepsilon}_{\mathcal{O}_\mu},
\; \tau_{\mathcal{K}_{\mathcal{O}_\mu}} \cdot T\pi_{\mathcal{O}_\mu} \cdot w) \nonumber \\
& = - \omega_{\mathcal{K}_{\mathcal{O}_\mu}}( X_{\mathcal{K}_{\mathcal{O}_\mu}} \cdot
\bar{\varepsilon}_{\mathcal{O}_\mu}, \; \tau_{\mathcal{K}_{\mathcal{O}_\mu}} \cdot T\bar{\lambda}_{\mathcal{O}_\mu} \cdot w)
+\omega_{\mathcal{K}_{\mathcal{O}_\mu}}(T\bar{\lambda}_{\mathcal{O}_\mu} \cdot X_\mathcal{K} \cdot
\varepsilon, \; \tau_{\mathcal{K}_{\mathcal{O}_\mu}} \cdot T\bar{\lambda}_{\mathcal{O}_\mu} \cdot w)\\
& = - \omega_{\mathcal{K}_{\mathcal{O}_\mu}}(\tau_{\mathcal{K}_{\mathcal{O}_\mu}}\cdot X_{h_{\mathcal{K}_{\mathcal{O}_\mu}}} \cdot
\bar{\varepsilon}_{\mathcal{O}_\mu}, \; \tau_{\mathcal{K}_{\mathcal{O}_\mu}} \cdot T\bar{\lambda}_{\mathcal{O}_\mu} \cdot w)
+\omega_{\mathcal{K}_{\mathcal{O}_\mu}}(T\bar{\lambda}_{\mathcal{O}_\mu} \cdot \tilde{X} \cdot
\varepsilon, \; \tau_{\mathcal{K}_{\mathcal{O}_\mu}} \cdot T\bar{\lambda}_{\mathcal{O}_\mu} \cdot w)\\
& = - \omega_{\mathcal{K}_{\mathcal{O}_\mu}}(\tau_{\mathcal{K}_{\mathcal{O}_\mu}}\cdot T\bar{\varepsilon}_{\mathcal{O}_\mu}
\cdot X_{h_{\mathcal{K}_{\mathcal{O}_\mu}} \cdot \bar{\varepsilon}_{\mathcal{O}_\mu}}, \; \tau_{\mathcal{K}_{\mathcal{O}_\mu}} \cdot T\bar{\lambda}_{\mathcal{O}_\mu} \cdot w)
+\omega_{\mathcal{K}_{\mathcal{O}_\mu}}(T\bar{\lambda}_{\mathcal{O}_\mu} \cdot \tilde{X} \cdot
\varepsilon, \; \tau_{\mathcal{K}_{\mathcal{O}_\mu}} \cdot T\bar{\lambda}_{\mathcal{O}_\mu} \cdot w)\\
& = \omega_{\mathcal{K}_{\mathcal{O}_\mu}}(T\bar{\lambda}_{\mathcal{O}_\mu} \cdot \tilde{X} \cdot
\varepsilon- \tau_{\mathcal{K}_{\mathcal{O}_\mu}}\cdot T\bar{\varepsilon}_{\mathcal{O}_\mu}
\cdot X_{h_{\mathcal{K}_{\mathcal{O}_\mu}} \cdot \bar{\varepsilon}_{\mathcal{O}_\mu}}, \; \tau_{\mathcal{K}_{\mathcal{O}_\mu}} \cdot T\bar{\lambda}_{\mathcal{O}_\mu} \cdot w).
\end{align*}
Because the $\mathbf{J}$-nonholonomic $R_o$-reduced distributional two-form
$\omega_{\mathcal{K}_{\mathcal{O}_\mu}}$ is non-degenerate, it follows that the equation
$T\bar{\gamma}_{\mathcal{O}_\mu}\cdot \tilde{X}^\varepsilon = X_{\mathcal{K}_{\mathcal{O}_\mu}}\cdot
\bar{\varepsilon}_{\mathcal{O}_\mu}, $ is equivalent to the equation
$T\bar{\lambda}_{\mathcal{O}_\mu}\cdot \tilde{X} \cdot \varepsilon
= \tau_{\mathcal{K}_{\mathcal{O}_\mu}}\cdot T\bar{\varepsilon}_{\mathcal{O}_\mu} \cdot X_{h_{\mathcal{K}_{\mathcal{O}_\mu}}
\cdot \bar{\varepsilon}_{\mathcal{O}_\mu}}. $
Thus, $\varepsilon$ and $\bar{\varepsilon}_{\mathcal{O}_\mu}$ satisfy the equation
$T\bar{\lambda}_{\mathcal{O}_\mu}\cdot \tilde{X} \cdot \varepsilon
= \tau_{\mathcal{K}_{\mathcal{O}_\mu}}\cdot T\bar{\varepsilon}_{\mathcal{O}_\mu} \cdot X_{h_{\mathcal{K}_{\mathcal{O}_\mu}}
\cdot \bar{\varepsilon}_{\mathcal{O}_\mu}}, $ if and only if they satisfy
the Type II of Hamilton-Jacobi equation
$T\bar{\gamma}_{\mathcal{O}_\mu}\cdot \tilde{X}^\varepsilon = X_{\mathcal{K}_{\mathcal{O}_\mu}}\cdot
\bar{\varepsilon}_{\mathcal{O}_\mu}.$
\hskip 0.3cm $\blacksquare$\\

For a given $\mathbf{J}$-nonholonomic regular orbit reducible RCH system
$(T^*Q,G,\omega,\mathbf{J},\mathcal{D},H, F, W)$ with the associated distributional RCH system
$(\mathcal{K},\omega_{\mathcal {K}},H_{\mathcal{K}}, F_{\mathcal{K}}, u_{\mathcal{K}})$
and the $\mathbf{J}$-nonholonomic $R_o$-reduced distributional RCH system
$(\mathcal{K}_{\mathcal{O}_\mu},\omega_{\mathcal{K}_{\mathcal{O}_\mu}},
h_{\mathcal{K}_{\mathcal{O}_\mu}},
f_{\mathcal {K}_{\mathcal{O}_\mu}}, u_{\mathcal {K}_{\mathcal{O}_\mu}})$,
we know that the nonholonomic dynamical vector field
$X_{\mathcal{K}}$ and the $\mathbf{J}$-nonholonomic $R_o$-reduced dynamical vector field
$X_{\mathcal{K}_{\mathcal{O}_\mu}}$ are $\pi_{\mathcal{O}_\mu}$-related,
that is, $X_{\mathcal{K}_{\mathcal{O}_\mu}}\cdot \pi_{\mathcal{O}_\mu}
=T\pi_{\mathcal{O}_\mu}\cdot X_{\mathcal{K}}\cdot i_{\mathcal{O}_\mu}.$ Then
we can prove the following Theorem 5.12, which states the
relationship between the solutions of Type II of Hamilton-Jacobi equations and
$\mathbf{J}$-nonholonomic regular orbit reduction.

\begin{theo}
For a given $\mathbf{J}$-nonholonomic regular orbit reducible RCH system
$(T^*Q, G, \omega,\mathbf{J}, \mathcal{D}, \\ H, F, W)$ with the associated distributional RCH system
$(\mathcal{K},\omega_{\mathcal {K}},H_{\mathcal{K}}, F_{\mathcal{K}}, u_{\mathcal{K}})$
and the $\mathbf{J}$-nonholonomic $R_o$-reduced distributional RCH system
$(\mathcal{K}_{\mathcal{O}_\mu},\omega_{\mathcal{K}_{\mathcal{O}_\mu}},
h_{\mathcal{K}_{\mathcal{O}_\mu}},
f_{\mathcal {K}_{\mathcal{O}_\mu}}, u_{\mathcal {K}_{\mathcal{O}_\mu}})$,
assume that $\gamma: Q \rightarrow T^*Q$ is an one-form on $Q$, and
$\varepsilon: T^* Q \rightarrow T^* Q $ is a symplectic map,
and $\bar{\gamma}_{\mathcal{O}_\mu}=\pi_{\mathcal{O}_\mu}(\gamma): Q \rightarrow \mathcal{M}_{\mathcal{O}_\mu} $,
and $\bar{\varepsilon}_{\mathcal{O}_\mu}=\pi_{\mathcal{O}_\mu}(\varepsilon):
\mathcal{M}\cap \mathbf{J}^{-1}(\mathcal{O}_\mu) (\subset T^*Q) \rightarrow (T^* Q)_{\mathcal{O}_\mu} $.
Under the hypotheses and notations of Theorem 5.11, then we have that $\varepsilon$
is a solution of the Type II of Hamilton-Jacobi equation $T\gamma\cdot
\tilde{X}^\varepsilon= X_{\mathcal{K}}\cdot \varepsilon, $ for the distributional
RCH system $(\mathcal{K},\omega_{\mathcal {K}},H_{\mathcal {K}},
F_{\mathcal {K}}, u_{\mathcal {K}})$, if and
only if $\varepsilon$ and $\bar{\varepsilon}_{\mathcal{O}_\mu} $ satisfy the Type II of Hamilton-Jacobi
equation $T\bar{\gamma}_{\mathcal{O}_\mu}\cdot \tilde{X}^\varepsilon=
X_{\mathcal{K}_{\mathcal{O}_\mu}}\cdot \bar{\varepsilon}_{\mathcal{O}_\mu}, $ for the
$\mathbf{J}$-nonholonomic $R_o$-reduced distributional RCH system
$(\mathcal{K}_{\mathcal{O}_\mu},\omega_{\mathcal{K}_{\mathcal{O}_\mu}},
h_{\mathcal{K}_{\mathcal{O}_\mu}},
f_{\mathcal {K}_{\mathcal{O}_\mu}}, u_{\mathcal {K}_{\mathcal{O}_\mu}})$.
\end{theo}

\noindent{\bf Proof: } For the dynamical vector field of the distributional RCH system
$(\mathcal{K},\omega_{\mathcal {K}}, H_{\mathcal{K}}, F_{\mathcal{K}}, u_{\mathcal{K}})$,
$\tilde{X}=X_{(\mathcal{K},\omega_{\mathcal{K}},
H_{\mathcal{K}}, F_{\mathcal{K}}, u_{\mathcal{K}})}
=X_\mathcal {K}+ F_{\mathcal{K}}+u_{\mathcal{K}}$,
and $F_{\mathcal{K}}=\tau_{\mathcal{K}}\cdot \textnormal{vlift}(F_{\mathcal{M}})X_H$,
and $u_{\mathcal{K}}=\tau_{\mathcal{K}}\cdot \textnormal{vlift}(u_{\mathcal{M}})X_H$,
note that $T\pi_{Q}\cdot \textnormal{vlift}(F_{\mathcal{M}})X_H=T\pi_{Q}\cdot \textnormal{vlift}(u_{\mathcal{M}})X_H=0, $
then we have that $T\pi_{Q}\cdot F_{\mathcal{K}}=T\pi_{Q}\cdot u_{\mathcal{K}}=0,$
and hence $T\pi_{Q}\cdot \tilde{X}\cdot \varepsilon
=T\pi_{Q}\cdot X_{\mathcal{K}}\cdot \varepsilon. $
Next, under the hypotheses and notations of Theorem 5.12,
the $R_o$-reduced symplectic space
$(T^\ast Q)_{\mathcal{O}_\mu}= \mathbf{J}^{-1}(\mathcal{O}_\mu)/G
\cong (\mathbf{J}^{-1}(\mu)/G) \times \mathcal{O}_\mu, $ with the $R_o$-reduced
symplectic form $\omega_{\mathcal{O}_\mu}$ uniquely characterized by
the relation $i_{\mathcal{O}_\mu}^\ast
\omega=\pi_{\mathcal{O}_{\mu}}^\ast \omega_{\mathcal{O}
_\mu}+\mathbf{J}_{\mathcal{O}_\mu}^\ast\omega_{\mathcal{O}_\mu}^+. $
Since
$\textmd{Im}(\gamma)\subset \mathcal{M}\cap \mathbf{J}^{-1}(\mu), $
and it is $G$-invariant,
in this case for any $V\in TQ, $ and $w\in TT^*Q, $
we have that
$\mathbf{J}_{\mathcal{O}_\mu}^\ast\omega_{\mathcal{O}_\mu}^+(T\gamma
\cdot V, \; w)=0, $ and hence
$\pi_{\mathcal{O}_\mu}^*\omega_{\mathcal{O}_\mu}=
i_{\mathcal{O}_\mu}^*\omega= \omega, $ along $\textmd{Im}(\gamma)$.
On the other hand, because
$\textmd{Im}(T\bar{\gamma}_{\mathcal{O}_\mu})\subset \mathcal{K}_{\mathcal{O}_\mu}, $
then $\omega_{\mathcal{K}_{\mathcal{O}_\mu}}\cdot
\tau_{\mathcal{K}_{\mathcal{O}_\mu}}=\tau_{\mathcal{K}_{\mathcal{O}_\mu}}\cdot
\omega_{\mathcal{M}_{\mathcal{O}_\mu}}= \tau_{\mathcal{K}_{\mathcal{O}_\mu}}\cdot
i_{\mathcal{M}_{\mathcal{O}_\mu}}^* \cdot \omega_{\mathcal{O}_\mu}, $ along
$\textmd{Im}(T\bar{\gamma}_{\mathcal{O}_\mu})$.
In addition, from the distributional Hamiltonian equation (5.9)
we have that $X_{\mathcal{K}}=\tau_{\mathcal{K}}\cdot X_H, $
and from the $\mathbf{J}$-nonholonomic $R_o$-reduced distributional Hamiltonian equation (5.10)
we have that $X_{\mathcal{K_{\mathcal{O}_\mu}}}
=\tau_{\mathcal{K_{\mathcal{O}_\mu}}}\cdot X_{h_{\mathcal{K_{\mathcal{O}_\mu}}}},$
and the vector fields $X_{\mathcal{K}}$
and $X_{\mathcal{K_{\mathcal{O}_\mu}}}$ are $\pi_{\mathcal{O}_\mu}$-related,
that is, $X_{\mathcal{K_{\mathcal{O}_\mu}}}\cdot \pi_{\mathcal{O}_\mu}
=T\pi_{\mathcal{O}_\mu}\cdot X_{\mathcal{K}}.$
Note that $\varepsilon(\mathcal{M})\subset \mathcal{M},$
and hence $X_H\cdot \varepsilon \in \mathcal{F}$ along $\varepsilon$,
and $\tau_{\mathcal{K}}\cdot X_{H}\cdot \varepsilon
= X_{\mathcal{K}}\cdot \varepsilon \in \mathcal{K}$.
Then $\tau_{\mathcal{K}_{\mathcal{O}_\mu}}\cdot T\pi_{\mathcal{O}_\mu}(X_{\mathcal{K}}\cdot \varepsilon)
=\tau_{\mathcal{K}_{\mathcal{O}_\mu}}\cdot (T\pi_{\mathcal{O}_\mu}(X_{\mathcal{K}}))\cdot (\varepsilon)
= \tau_{\mathcal{K}_{\mathcal{O}_\mu}}\cdot (X_{\mathcal{K}_{\mathcal{O}_\mu}}\cdot \pi_{\mathcal{O}_\mu})\cdot (\varepsilon)
= \tau_{\mathcal{K}_{\mathcal{O}_\mu}}\cdot X_{\mathcal{K}_{\mathcal{O}_\mu}}\cdot \pi_{\mathcal{O}_\mu}(\varepsilon)
= X_{\mathcal{K}_{\mathcal{O}_\mu}}\cdot \bar{\varepsilon}.$
Thus, using the $\mathbf{J}$-nonholonomic $R_o$-reduced distributional two-form
$\omega_{\mathcal{K}_{\mathcal{O}_\mu}}$, note that
$\tau_{\mathcal{K}_{\mathcal{O}_\mu}}\cdot T\bar{\gamma}_{\mathcal{O}_\mu}
= T\bar{\gamma}_{\mathcal{O}_\mu}, $
for any $w \in \mathcal{F}, \; \tau_{\mathcal{K}}\cdot w\neq 0, $ and
$\tau_{\mathcal{K}_{\mathcal{O}_\mu}}\cdot T\pi_{\mathcal{O}_\mu} \cdot w \neq 0, $
then we have that
\begin{align*}
& \omega_{\mathcal{K}_{\mathcal{O}_\mu}}(T\bar{\gamma}_{\mathcal{O}_\mu}
\cdot \tilde{X}^\varepsilon- X_{\mathcal{K}_{\mathcal{O}_\mu}}\cdot \bar{\varepsilon}_{\mathcal{O}_\mu}, \;
\tau_{\mathcal{K}_{\mathcal{O}_\mu}}\cdot T\pi_{\mathcal{O}_\mu} \cdot w) \\
& = \omega_{\mathcal{K}_{\mathcal{O}_\mu}}(T\bar{\gamma}_{\mathcal{O}_\mu} \cdot \tilde{X}^\varepsilon, \;
\tau_{\mathcal{K}_{\mathcal{O}_\mu}}\cdot T\pi_{\mathcal{O}_\mu} \cdot w)-
\omega_{\mathcal{K}_{\mathcal{O}_\mu}}(X_{\mathcal{K}_{\mathcal{O}_\mu}}\cdot \bar{\varepsilon}_{\mathcal{O}_\mu},
\; \tau_{\mathcal{K}_{\mathcal{O}_\mu}} \cdot T\pi_{\mathcal{O}_\mu} \cdot w) \\
& = \omega_{\mathcal{K}_{\mathcal{O}_\mu}}(\tau_{\mathcal{K}_{\mathcal{O}_\mu}}
\cdot T\bar{\gamma}_{\mathcal{O}_\mu} \cdot \tilde{X}^\varepsilon, \;
\tau_{\mathcal{K}_{\mathcal{O}_\mu}}\cdot T\pi_{\mathcal{O}_\mu} \cdot w)-
\omega_{\mathcal{K}_{\mathcal{O}_\mu}}(\tau_{\mathcal{K}_{\mathcal{O}_\mu}}
\cdot T\pi_{\mathcal{O}_\mu}(X_{\mathcal{K}}\cdot \varepsilon), \;
\tau_{\mathcal{K}_{\mathcal{O}_\mu}} \cdot T\pi_{\mathcal{O}_\mu} \cdot w) \\
& = \omega_{\mathcal{K}_{\mathcal{O}_\mu}}\cdot \tau_{\mathcal{K}_{\mathcal{O}_\mu}}
(T\pi_{\mathcal{O}_\mu} \cdot T\gamma \cdot \tilde{X}^\varepsilon, \; T\pi_{\mathcal{O}_\mu} \cdot w)
-\omega_{\mathcal{K}_{\mathcal{O}_\mu}}\cdot \tau_{\mathcal{K}_{\mathcal{O}_\mu}}
(T\pi_{\mathcal{O}_\mu} \cdot X_{\mathcal{K}}\cdot \varepsilon, \; T\pi_{\mathcal{O}_\mu} \cdot w)\\
& = \tau_{\mathcal{K}_{\mathcal{O}_\mu}}\cdot
i_{\mathcal{M}_{\mathcal{O}_\mu}}^* \cdot \omega_{\mathcal{O}_\mu}
(T\pi_{\mathcal{O}_\mu} \cdot T\gamma \cdot \tilde{X}^\varepsilon, \; T\pi_{\mathcal{O}_\mu} \cdot w)
-\tau_{\mathcal{K}_{\mathcal{O}_\mu}}\cdot
i_{\mathcal{M}_{\mathcal{O}_\mu}}^* \cdot \omega_{\mathcal{O}_\mu}
(T\pi_{\mathcal{O}_\mu} \cdot X_{\mathcal{K}}\cdot \varepsilon, \; T\pi_{\mathcal{O}_\mu} \cdot w)\\
& = \tau_{\mathcal{K}_{\mathcal{O}_\mu}}\cdot
i_{\mathcal{M}_{\mathcal{O}_\mu}}^* \cdot \pi_{\mathcal{O}_\mu}^*\omega_{\mathcal{O}_\mu}
(T\gamma \cdot \tilde{X}^\varepsilon, \; w)
-\tau_{\mathcal{K}_{\mathcal{O}_\mu}}\cdot
i_{\mathcal{M}_{\mathcal{O}_\mu}}^* \cdot \pi_{\mathcal{O}_\mu}^*\omega_{\mathcal{O}_\mu}
(X_{\mathcal{K}}\cdot \varepsilon, \; w).
\end{align*}
In the case we note that $\tau_{\mathcal{K}_{\mathcal{O}_\mu}}\cdot
i_{\mathcal{M}_{\mathcal{O}_\mu}}^* \cdot \pi_{\mathcal{O}_\mu}^*\omega_{\mathcal{O}_\mu}
=\tau_{\mathcal{K}}\cdot i_{\mathcal{M}}^* \cdot \omega= \omega_{\mathcal{K}}\cdot \tau_{\mathcal{K}}, $
and
$\tau_{\mathcal{K}}\cdot T\gamma =T\gamma, \; \tau_{\mathcal{K}} \cdot \tilde{X}= X_{\mathcal{K}}$,
since $\textmd{Im}(\gamma)\subset
\mathcal{M}, $ and $\textmd{Im}(T\gamma)\subset \mathcal{K}. $
Thus, we have that
\begin{align*}
& \omega_{\mathcal{K}_{\mathcal{O}_\mu}}(T\bar{\gamma}_{\mathcal{O}_\mu}
\cdot \tilde{X}^\varepsilon- X_{\mathcal{K}_{\mathcal{O}_\mu}}\cdot \bar{\varepsilon}_{\mathcal{O}_\mu}, \;
\tau_{\mathcal{K}_{\mathcal{O}_\mu}}\cdot T\pi_{\mathcal{O}_\mu} \cdot w) \\
& = \omega_{\mathcal{K}}\cdot \tau_{\mathcal{K}}(T\gamma \cdot \tilde{X}^\varepsilon, \; w)
-\omega_{\mathcal{K}}\cdot \tau_{\mathcal{K}}(X_{\mathcal{K}}\cdot \varepsilon, \; w)\\
& = \omega_{\mathcal{K}}(\tau_{\mathcal{K}}\cdot T\gamma \cdot \tilde{X}^\varepsilon, \; \tau_{\mathcal{K}}\cdot w)
-\omega_{\mathcal{K}}(\tau_{\mathcal{K}}\cdot X_{\mathcal{K}}\cdot \varepsilon, \; \tau_{\mathcal{K}}\cdot w)\\
& = \omega_{\mathcal{K}}(T\gamma \cdot \tilde{X}^\varepsilon, \; \tau_{\mathcal{K}}\cdot w)
-\omega_{\mathcal{K}}(X_{\mathcal{K}}\cdot \varepsilon, \; \tau_{\mathcal{K}}\cdot w)\\
& = \omega_{\mathcal{K}}(T\gamma \cdot \tilde{X}^\varepsilon- X_{\mathcal{K}}\cdot \varepsilon, \; \tau_{\mathcal{K}}\cdot w).
\end{align*}
Because the distributional two-form $\omega_{\mathcal{K}}$ and
the $\mathbf{J}$-nonholonomic $R_o$-reduced distributional
two-form $\omega_{\mathcal{K}_{\mathcal{O}_\mu}}$ are both non-degenerate,
it follows that the equation
$T\bar{\gamma}_{\mathcal{O}_\mu}\cdot \tilde{X}^\varepsilon=
X_{\mathcal{K}_{\mathcal{O}_\mu}}\cdot \bar{\varepsilon}_{\mathcal{O}_\mu}, $ is equivalent to the equation
$T\gamma\cdot \tilde{X}^\varepsilon= X_{\mathcal{K}}\cdot \varepsilon$. Thus,
$\varepsilon$ is a solution of the Type II of Hamilton-Jacobi equation
$T\gamma\cdot \tilde{X}^\varepsilon= X_{\mathcal{K}}\cdot \varepsilon, $ for the distributional
RCH system $(\mathcal{K},\omega_{\mathcal {K}},H_{\mathcal {K}},
F_{\mathcal {K}}, u_{\mathcal {K}})$, if and only if
$\varepsilon$ and $\bar{\varepsilon}_{\mathcal{O}_\mu} $ satisfy the Type II of Hamilton-Jacobi
equation $T\bar{\gamma}_{\mathcal{O}_\mu}\cdot \tilde{X}^\varepsilon=
X_{\mathcal{K}_{\mathcal{O}_\mu}}\cdot \bar{\varepsilon}_{\mathcal{O}_\mu}, $ for the
$\mathbf{J}$-nonholonomic $R_o$-reduced distributional RCH system
$(\mathcal{K}_{\mathcal{O}_\mu},\omega_{\mathcal{K}_{\mathcal{O}_\mu}},
h_{\mathcal{K}_{\mathcal{O}_\mu}},
f_{\mathcal {K}_{\mathcal{O}_\mu}}, u_{\mathcal {K}_{\mathcal{O}_\mu}})$.  \hskip 0.3cm
$\blacksquare$ \\

\begin{rema}
It is worthy of noting that,
the Type I of Hamilton-Jacobi equation
$T\bar{\gamma}_{\mathcal{O}_\mu}\cdot \tilde{X}^\gamma
= X_{\mathcal{K}_{\mathcal{O}_\mu}}\cdot \bar{\gamma}_{\mathcal{O}_\mu}. $
is the equation of the $\mathbf{J}$-nonholonomic
$R_o$-reduced differential one-form $\bar{\gamma}_{\mathcal{O}_\mu}$; and
the Type II of Hamilton-Jacobi equation
$T\bar{\gamma}_{\mathcal{O}_\mu}\cdot \tilde{X}^\varepsilon
= X_{\mathcal{K}_{\mathcal{O}_\mu}}\cdot \bar{\varepsilon}_{\mathcal{O}_\mu} .$
is the equation of the symplectic diffeomorphism map $\varepsilon$
and the $\mathbf{J}$-nonholonomic $R_o$-reduced
symplectic diffeomorphism map $\bar{\varepsilon}_{\mathcal{O}_\mu}. $
If a $\mathbf{J}$-nonholonomic regular orbit reducible RCH system we considered
$(T^*Q,G,\omega, \mathbf{J}, \mathcal{D}, H, F, u)$ has not any constrains,
in this case, the $\mathbf{J}$-nonholonomic $R_o$-reduced distributional
RCH system is just the $R_o$-reduced RCH system itself.
From the above Type I and Type II of Hamilton-Jacobi theorems, that is,
Theorem 5.10 and Theorem 5.11, we can get the Theorem 4.2
and Theorem 4.3 in Wang \cite{wa13d}.
It shows that Theorem 5.10 and Theorem 5.11 can be regarded as an extension of two types of
Hamilton-Jacobi theorem for the $R_o$-reduced RCH system
to that for the system with nonholonomic context.
If the $\mathbf{J}$-nonholonomic regular orbit reducible RCH system we considered
$(T^*Q,G,\omega, \mathbf{J}, \mathcal{D}, H, F, u)$
has not any the external force and control, that is, $F=0 $ and $u=0$,
in this case, the $\mathbf{J}$-nonholonomic regular orbit reducible RCH system
is just the $\mathbf{J}$-nonholonomic regular orbit reducible
Hamiltonian system $(T^*Q,G,\omega,\mathbf{J},\mathcal{D},H)$.
and with the canonical symplectic form $\omega$ on $T^*Q$.
From the above Type I and Type II of Hamilton-Jacobi theorems, that is,
Theorem 5.10 and Theorem 5.11, we can get the Theorem 5.9
and Theorem 5.10 in Le\'{o}n and Wang \cite{lewa15}.
It shows that Theorem 5.10 and Theorem 5.11 can be regarded as an extension of two types of
Hamilton-Jacobi theorem for the $\mathbf{J}$-nonholonomic regular orbit reducible Hamiltonian system
to that for the system with the external force and control.
In particular, in this case,
if the $\mathbf{J}$-nonholonomic regular orbit reducible RCH system we considered has not any constrains,
that is, $F=0, \; u=0 $ and $\mathcal{D}=\emptyset$, then
the $\mathbf{J}$-nonholonomic regular orbit reducible RCH system
is just a regular orbit reducible Hamiltonian system $(T^*Q,G,\omega,H)$
with the canonical symplectic form $\omega$ on $T^*Q$,
we can obtain two types of Hamilton-Jacobi
equation for the associated $R_o$-reduced Hamiltonian system,
which is given in Wang \cite{wa17}.
Thus, Theorem 5.10 and Theorem 5.11 can be regarded as an extension of two types of Hamilton-Jacobi
theorem for a regular orbit reducible Hamiltonian system to that for the system with external force,
control and nonholonomic constrain.
\end{rema}

\begin{rema}
If $(T^\ast Q, \omega)$ is a connected symplectic manifold, and
$\mathbf{J}:T^\ast Q\rightarrow \mathfrak{g}^\ast$ is a
non-equivariant momentum map with a non-equivariance group
one-cocycle $\sigma: G\rightarrow \mathfrak{g}^\ast$,
which is defined by $\sigma(g):=\mathbf{J}(g\cdot
z)-\operatorname{Ad}^\ast_{g^{-1}}\mathbf{J}(z)$, where $g\in G$ and
$z\in T^\ast Q$. Then we know that $\sigma$ produces a new affine
action $\Theta: G\times \mathfrak{g}^\ast \rightarrow
\mathfrak{g}^\ast $ defined by
$\Theta(g,\mu):=\operatorname{Ad}^\ast_{g^{-1}}\mu + \sigma(g)$,
where $\mu \in \mathfrak{g}^\ast$, with respect to which the given
momentum map $\mathbf{J}$ is equivariant. Assume that $G$ acts
freely and properly on $T^\ast Q$, and $\mathcal{O}_\mu= G\cdot \mu
\subset \mathfrak{g}^\ast$ denotes the G-orbit of the point $\mu \in
\mathfrak{g}^\ast$ with respect to this affine action $\Theta$,
and $\mu$ is a regular value of $\mathbf{J}$.
Then the quotient space $(T^\ast
Q)_{\mathcal{O}_\mu}=\mathbf{J}^{-1}(\mathcal{O}_\mu)/G$ is also a symplectic
manifold with symplectic form $\omega_{\mathcal{O}_\mu}$ uniquely characterized by
$(5.7)$, see Ortega and Ratiu \cite{orra04}. In this case,
we can also define the $\mathbf{J}$-nonholonomic regular orbit reducible
RCH system $(T^*Q,G,\omega,\mathbf{J},\mathcal{D},H, F, W)$
with the associated distributional RCH system
$(\mathcal{K},\omega_{\mathcal {K}},H_{\mathcal{K}}, F_{\mathcal{K}}, u_{\mathcal{K}})$
and the $\mathbf{J}$-nonholonomic $R_o$-reduced distributional RCH system
$(\mathcal{K}_{\mathcal{O}_\mu},\omega_{\mathcal{K}_{\mathcal{O}_\mu}},
h_{\mathcal{K}_{\mathcal{O}_\mu}},
f_{\mathcal {K}_{\mathcal{O}_\mu}}, u_{\mathcal {K}_{\mathcal{O}_\mu}})$,
and prove the Type I and Type II of
the Hamilton-Jacobi theorem for the $\mathbf{J}$-nonholonomic $R_o$-reduced
distributional RCH system
$(\mathcal{K}_{\mathcal{O}_\mu},\omega_{\mathcal{K}_{\mathcal{O}_\mu}},
h_{\mathcal{K}_{\mathcal{O}_\mu}},
f_{\mathcal {K}_{\mathcal{O}_\mu}}, u_{\mathcal {K}_{\mathcal{O}_\mu}})$
by using the above similar way,
in which the $\mathbf{J}$-nonholonomic $R_o$-reduced space $((\mathcal{K}_{\mathcal{O}_\mu},\omega_{\mathcal{K}_{\mathcal{O}_\mu}})$
is determined by the affine action and $\mathbf{J}$-nonholonomic regular orbit reduction.
\end{rema}

In this paper, we study the
Hamilton-Jacobi theory for
the nonholonomic RCH system and
the nonholonomic reducible RCH system
on a cotangent bundle, by using the distributional RCH system
and the nonholonomic reduced distributional RCH system,
and by combining with the analysis
of dynamic systems and the geometric reduction theory of
RCH systems.These researches,
from the geometrical point of view, reveal the internal relationships of
nonholonomic constraints, distributional two forms and
nonholonomic dynamical vector fields of an
RCH system and its nonholonomic reduced RCH systems.
In particular, note that Marsden et al. in \cite{mawazh10} set up the
regular reduction theory of an RCH system on a symplectic fiber
bundle, by using momentum map and the associated reduced symplectic
form, and from the viewpoint of completeness of Marsden-Weinstein symplectic
reduction, and some developments around the above work are given in Wang and
Zhang \cite{wazh12}, Ratiu and Wang \cite{rawa12}, and
Wang \cite{wa15a, wa13d, wa20a, wa13e}.
Since the Hamilton-Jacobi theory
is developed based on the Hamiltonian picture of dynamics, it is
natural idea to extend the regular reduction and Hamilton-Jacobi theory to the
controlled magnetic Hamiltonian (CMH) system and its a variety of reduced CMH systems,
and it is also possible to describe the relationship between the
CMH-equivalence for the CMH systems and the solutions
of corresponding Hamilton-Jacobi equations, see Wang \cite{wa21b}
for more details. In addition, it is also
an important topic for us to explore and reveal the deeply internal
relationships between the geometrical structures of phase spaces and the dynamical
vector fields of the regular controlled Lagrangian systems.
All of topics are our goal of pursuing and inheriting.\\

\noindent{\bf Acknowledgments:}
It is an important responsibility for the scientific researchers
to guarantee the development of the scientific research work
on a right way. Now, I understand that it is not easy to do for Professor
Jerrold E. Marsden, also for us, to fight with wrong, ignorance
and naivete. In particular,
it is an important task for us to correct and develop well
the research work of Professor Jerrold E. Marsden,
such that we never feel sorry for his great cause.\\

\end{document}